\documentclass{amsart}
\usepackage{amssymb}
\usepackage{amsfonts}

\newtheorem{theorem}{Theorem}
\theoremstyle{plain}

\newtheorem{corollary}{Corollary}

\newtheorem{definition}{Definition}
\newtheorem{example}{Example}

\newtheorem{lemma}{Lemma}
\newtheorem{notation}{Notation}

\newtheorem{proposition}{Proposition}
\newtheorem{remark}{Remark}

\newtheorem{terminology}{Terminology}
\numberwithin{equation}{section}
\input{tcilatex}

\begin{document}
\title[Extendable Cohomologies]{Extendable Cohomologies for Complex Analytic
Varieties}
\author{Carlo Perrone}
\address{Dipartimento di Matematica della II Universit\`{a} di Roma "Tor
Vergata" \\
Via della Ricerca Scientifica, 1 00133 - Roma (Italia)}
\date{October 16, 2008}
\dedicatory{Citro Cucurbit\ae que maximis}

\begin{abstract}
We introduce a cohomology, called extendable cohomology, for abstract
complex singular varieties based on suitable differential forms. Beside a
study of the general properties of such a cohomology, we show that, given a
complex vector bundle, one can compute its topological Chern classes using
the extendable Chern classes, defined via a Chern-Weil type theory. We also
prove that the localizations of the extendable Chern classes represent the
localizations of the respective topological Chern classes, thus obtaining an
abstract residue theorem for compact singular complex analytic varieties. As
an application of our theory, we prove a Camacho-Sad type index theorem for
holomorphic foliations of singular complex varieties.
\end{abstract}

\maketitle

\textbf{Introduction.} One of the more important contributions to the study
of complex vector bundles over differentiable manifolds has been given by
the Chern-Weil theory. Thanks to such a theory it is possible to describe
the topological Chern classes of a complex vector bundle on a manifold
(which lie in the topological cohomology groups of the manifold) by means of
the differentiable Chern classes of the bundle (which belong to the de Rham
cohomology groups of the manifold).

By their very definition, the differentiable Chern classes of a complex
vector bundle are built starting from suitable differentiable differential
forms on the manifold. This is the reason why, until now, it was impossible
to achieve a generalization of the Chern-Weil theory allowing to study
complex vector bundles over singular varieties. In fact, the hurdles for
having such a theory are tied to the difficulties of giving an appropriate
definition of differential forms on singular spaces.

In this paper we solve the problem of extending the Chern-Weil theory to the
case of abstract complex analytic varieties. Namely, we introduce a suitable
notion of differential forms, the extendable differentiable differential
forms, we develop a cohomology theory based on such forms, we define the
extendable Chern classes for differentiable complex vector bundles over
complex analytic varieties and we prove that these classes represent the
topological Chern classes of the bundle.

The starting point is the following. In the case of complex analytic
varieties, it can be given several natural definitions of holomorphic
differential forms. Nevertheless, although remarkable results have been
obtained, the development of the theories based on such holomorphic forms
did not carry on, because of the failure of the Poincar\'{e}\ lemma. Namely,
the cohomologies associated with these holomorphic forms are not, in
general, locally trivial (cp. \cite{Ferrari 1}, \cite{Ferrari 2}, \cite%
{Herrera Tesi}, \cite{Herrera Articolo}, \cite{Bloom-Herrera}). Anyway, all
these definitions of holomorphic differential form are such that the sheaves
of germs of such forms share the following property: it is a sheaf of
modules (over the ring of holomorphic functions on the variety) that, even
if it is not necessarily locally free, it is always coherent (cp. \cite%
{Ferrari 2}).

On the other hand, we need differential forms that are differentiable but
not necessarily holomorphic. So, we define the sheaf of extendable
differentiable differential forms we are interested in by tensorizing one of
the sheaves of holomorphic differential forms with the sheaf of rings of
differentiable functions on the variety. In view of our aim of extending the
Chern-Weil theory, the choice of the sheaf of holomorphic differential forms
is not important, even if, of course, different choices generally lead to
different results. In fact, we only need the coherence of the sheaf of
extendable differentiable differential forms.

Extendable differentiable differential forms enjoy many properties of
differential forms on smooth manifolds (for example, they always have
bounded coefficients), even if the proof of some of these properties is not
trivial, because of the presence of singularities (see, for example, Lemma %
\ref{Esiste Part. d'Unita' "Liscia" copy(1)}, Proposition \ref{Integrale per
Triangolazione(Omega)}, Theorem \ref{Teorema di Stokes}). On the other hand,
the cohomology groups associated with extendable forms (the extendable
cohomology groups) are not locally trivial (cp. Example \ref{Es. Bloom
Herrera Complesso}).

Let $E\rightarrow X$ be a differentiable (holomorphic) $\mathbb{%
\mathbb{C}
}$-vector bundle over an abstract finite dimensional complex analytic
variety $X$. By using the theory of extendable forms, we introduce the
notions of extendable linear connections and extendable curvatures for $E$.
Let $\nabla $ be an extendable linear connection for $E$. We define the
extendable Chern forms $c_{ext}^{\bullet }(\nabla )$ associated with $\nabla 
$ and, arguing as in the smooth case, we show that $c_{ext}^{\bullet
}(\nabla )$ are closed and only depend on $E$. Then, we define extendable
Chern classes of $E$ as the cohomology classes $c_{ext}^{\bullet
}(E)=[c_{ext}^{\bullet }(\nabla )]$.

Next, we may define an operator of integration of extendable forms on
simplices (recall that any complex analytic variety is triangulable). Let $%
H_{ext}^{\bullet }(X)$ denote the extendable cohomology groups of $X$. Then
the operator of integration induces a homomorphism $H^{\bullet
}:H_{ext}^{\bullet }(X)\rightarrow H^{\bullet }(X)$ between extendable and
topological cohomology groups (cp. Section \ref{Integration (Sezione)}). We
prove the following theorem (cp. Theorem \ref{Classi Estendibili e
Topologiche coincidono copy(1)}).

\begin{enumerate}
\item[\textbf{Theorem}] \textit{Let }$X$\textit{\ be an abstract complex
analytic variety of complex dimension }$n$\textit{\ and }$E\rightarrow X$%
\textit{\ a differentiable (holomorphic) complex vector bundle of rank }$e$%
\textit{. Take }$q\in \{1,...,n\}$\textit{\ with }$q\leq e$\textit{. Then }$%
c_{top}^{q}(E)=H^{2q}\left( c_{ext}^{q}(E)\right) $\textit{, where }$%
c_{top}^{\bullet }(E)$\textit{\ denote the topological Chern classes of }$E$%
\textit{.}
\end{enumerate}

One of the topics we deal with in this paper is to prove residue theorems
for holomorphic complex vector bundles over compact irreducible abstract
complex analytic varieties. Let $X$ be a complex analytic variety of
dimension $n$ and $\chi \in H^{\bullet }(X)$ an element in the topological
cohomology groups of $X$. It can happen that the class $\chi $ could
represent the first order obstruction to the existence of a certain global
object $\mathfrak{o}$\ on $X$ (for example, topological Chern classes
represent the first order obstruction to the existence of global frames for
complex vector bundles). Generally and very roughly speaking, one asks where
on the variety $X$ the existence of $\mathfrak{o}$\ is obstructed. Namely,
one asks where on $X$ the class $\chi $ vanishes. Let $S$ denote the loci
where $\mathfrak{o}$\ exists ($\chi $ vanishes). It could be possible to
make a clever choice of $S$, even if such loci, in general, are not unique.
Then, on $X\setminus S$, that is outside $S$, the object $\mathfrak{o}$
exists and the class $\chi $ vanishes.

Now, assume that $X$ is compact and let $\boldsymbol{P}_{\bullet }^{\ast
}:H^{\bullet }\left( X\right) \rightarrow H_{2n-\bullet }\left( X\right) $
be the Poincar\'{e} homomorphism. Suppose that $S$ is an analytic subvariety
of $X$ and consider the exact sequence%
\begin{equation*}
\cdots \rightarrow H^{\bullet }\left( X,X\setminus S\right) \rightarrow
H^{\bullet }\left( X\right) \rightarrow H^{\bullet }\left( X\setminus
S\right) \rightarrow \cdots
\end{equation*}%
If the image of $\chi \in H^{\bullet }\left( X\right) $ in $H^{\bullet
}\left( X\setminus S\right) $ is $0$, then there exist $\kappa \in
H^{\bullet }\left( X,X\setminus S\right) $ whose image in $H^{\bullet
}\left( X\right) $\ is $\chi $. Such a $\kappa $ is the localization of $%
\chi $ at $S$ and, in general, it is not unique. Nevertheless, if $S$ is
compact, by taking into account the Alexander-Lefschetz homomorphism $%
\boldsymbol{A}_{S,\bullet }^{\ast }:H^{\bullet }\left( X,X\setminus S\right)
\rightarrow H_{2n-\bullet }\left( S\right) $ and the commutative diagram%
\begin{equation*}
\begin{array}[b]{ccccc}
H^{\bullet }\left( X,X\setminus S\right) &  & \rightarrow &  & H^{\bullet
}\left( X\right) \\ 
&  &  &  &  \\ 
\downarrow _{\text{ }^{\boldsymbol{A}_{S,\bullet }^{\ast }}} &  &  &  & 
\downarrow _{\text{ }^{\boldsymbol{P}_{\bullet }^{\ast }}} \\ 
&  &  &  &  \\ 
H_{2n-\bullet }\left( S\right) &  & \overset{i_{\ast }}{\rightarrow } &  & 
H_{2n-\bullet }\left( X\right)%
\end{array}%
,
\end{equation*}%
we get the formula $\boldsymbol{P}_{\bullet }^{\ast }(\chi )=(i_{\ast }\circ 
\boldsymbol{A}_{S,\bullet }^{\ast })(\kappa )$. This is an "index theorem"
(see \cite{SuwaLibro}). If $\bullet =2n$ and $S$ is a finite set of points $%
\{p_{\nu }\}$, then $H_{0}(S)=\oplus _{\nu }H_{0}(p_{\nu })$ and $%
\boldsymbol{A}_{S,\bullet }^{\ast }(\kappa )=\tsum\nolimits_{\nu }Res(\kappa
,p_{\nu })$, where $Res(\kappa ,p_{\nu })\in H_{0}(p_{\nu })$ is "the
residue of $\kappa $ at $p_{\nu }$". So, the index theorem can be written as%
\begin{equation*}
\boldsymbol{P}_{2n}^{\ast }(\chi )=\tsum\nolimits_{\nu }i_{\ast }(Res(\kappa
,p_{\nu })).
\end{equation*}%
Next, taking into account the homomorphism $H^{2n}:H_{ext}^{2n}(X)%
\rightarrow H^{2n}(X)$ induced by integration on simplices, we have $%
\boldsymbol{P}_{2n}^{\ast }\circ H^{2n}=\tint\nolimits_{[X]}$, where $[X]$
is the fundamental class of $X$. So, if $\chi _{ext}\in H_{ext}^{2n}(X)$ is
such that $\chi =H_{ext}^{2n}(\chi _{ext})$, then we get%
\begin{equation*}
\tint\nolimits_{\lbrack X]}\chi _{ext}=\tsum\nolimits_{\nu }i_{\ast
}(Res(\kappa ,p_{\nu }))
\end{equation*}%
Namely, a "residue theorem". We prove the following theorem (cp. Theorem \ref%
{Teor. dei Res.}).

\begin{enumerate}
\item[\textbf{Theorem}] \textit{Let }$X$\textit{\ be a compact and
irreducible complex analytic variety of complex dimension }$n$\textit{\ and }%
$E\rightarrow X$\textit{\ a holomorphic complex vector bundle of rank }$e$%
\textit{. Take }$q\in \{0,...,n\}$\textit{\ with }$q\leq e$\textit{\ and set 
}$r=e-q+1$\textit{. Let }$s^{\left( r\right) }$\textit{\ be a holomorphic }$%
r $\textit{-section of }$E$\textit{\ and }$S$\textit{\ the singular locus of 
}$s^{(r)}$\textit{\ and }$c_{top}^{q}(E,s^{\left( r\right) })$\textit{\ the
localization at }$S$\textit{\ of }$c_{top}^{q}(E)$\textit{\ determined by }$%
s^{\left( r\right) }$\textit{. Set }$TopRes_{c_{top}^{q}}(E,s^{\left(
r\right) },S)=\boldsymbol{A}_{S,2q}^{\ast }(c_{top}^{q}(E,s^{\left( r\right)
}))$\textit{. Then}%
\begin{equation*}
\boldsymbol{P}_{2q}^{\ast }\circ H^{2q}(c_{ext}^{q}\left( E\right) )=i_{\ast
}(TopRes_{c_{top}^{q}}(E,s^{\left( r\right) },S)).
\end{equation*}%
\textit{If }$q=n$\textit{, then}%
\begin{equation*}
\tint\nolimits_{\lbrack X]}c_{ext}^{n}(E)=i_{\ast
}(TopRes_{c_{top}^{n}}(E,s^{\left( r\right) },S)).
\end{equation*}
\end{enumerate}

Actually, a residue theorem becomes really useful only if the residue can be
explicitly and easily computed. In order to solve this problem, several
people widely and successfully used the theory of \v{C}ech-de Rham. Among
other authors, we mention J. P. Brasselet, D. Lehmann and T. Suwa, who,
furthermore, greatly developed such a theory (see \cite{Lehmann-Suwa}, \cite%
{Bracci-Suwa}, \cite{SuwaLibro} and references therein). Indeed, \v{C}ech-de
Rham theory provides with very handleable tools to explicitly compute the
residue, at least in the case of isolated singularities. We should note
that, until now, only the cases of manifolds, submanifolds embedded in a
manifold and, at most, subvarieties embedded in a manifold were studied (see 
\cite{Lehmann-Suwa}, \cite{SuwaLibro}, \cite{Suwa Sao Carlos Versione II}).
In fact, in these cases, in order to have a \v{C}ech-de Rham type theory,
the differentiable differential forms of the ambient are used. On the other
hand, our theory allows to take into account the case of abstract complex
analytic varieties.

In this paper we also develop a \v{C}ech-de Rham type theory for extendable
forms. By means of such a theory, we can compute the residue of the
localizations of several characteristic classes, at least if the
singularities are isolated. For instance, we obtain a generalization of
Camacho-Sad index theorem (see Theorem \ref{Teorema Camacho-Sad} for a more
general statement).

\begin{enumerate}
\item[\textbf{Theorem}] \textit{Let }$X$\textit{\ be an abstract complex
analytic variety of complex dimension }$2$\textit{, }$\mathcal{F}$\textit{\
a holomorphic foliation of }$X$\textit{\ and }$Y$\textit{\ an }$\mathcal{F}$%
\textit{-invariant globally irreducible Cartier divisor of }$X$\textit{\
such that }$Y\nsubseteq Sing(X)$\textit{. Set }$S=(Sing(\mathcal{F})\cap
Y)\cup Sing(Y)$\textit{\ and let }$N_{Y}\rightarrow Y$\textit{\ be the line
bundle }$\mathcal{O}([Y])$\textit{. Then}%
\begin{equation*}
\tint\nolimits_{\lbrack Y]}c_{ext}^{1}(N_{Y})=i_{\ast
}(Res_{c_{ext}^{1}}(N_{Y},\mathcal{F},S)),
\end{equation*}%
\textit{where }$c_{ext}^{1}(N_{Y},\mathcal{F},S)$\textit{\ is the
localization of }$c_{ext}^{1}(N_{Y})$\textit{\ at }$S$\textit{\ determined
by }$\mathcal{F}$\textit{\ and }$Res_{c_{ext}^{1}}(N_{Y},\mathcal{F},S)$%
\textit{\ are, in fact, complex numbers which only depend on the behaviour }$%
\mathcal{F}$\textit{\ of around }$S$.

\item[ ] \qquad \textit{Suppose that }$S$\textit{\ only contains an isolated
singular point }$p\in Sing(Y)\cap Sing(\mathcal{F})\cap Sing(X)$\textit{\
and that the stalk }$\mathcal{F}_{p}$\textit{\ is generated on }$\mathcal{O}%
_{X,p}$\textit{\ by a single element of }$\mathcal{T}X_{p}$\textit{. Let }$%
(h,y)$\textit{\ be local coordinates on }$X^{Reg}$\textit{\ near }$p$\textit{%
\ such that }$y$\textit{\ is a local coordinate on }$Y^{\prime }=Y\setminus
((Sing(X)\cap Y)\cup Sing(Y))$\textit{\ near each point of }$Y^{\prime
}\setminus \{p\}$\textit{. If the holomorphic vector field }$\digamma \in 
\mathcal{T}X$\textit{\ generating }$\mathcal{F}$\textit{\ is locally given
by }$\digamma =a(h,y)h\frac{\partial }{\partial h}+b(h,y)\frac{\partial }{%
\partial y}$\textit{, with }$a$\textit{\ and }$b$\textit{\ holomorphic
functions such that }$b(0,y)$\textit{\ is not identically equal to zero, then%
}%
\begin{equation*}
i_{\ast }(Res_{c_{ext}^{1}}(N_{Y},\mathcal{F},p))=\tfrac{1}{2\pi \sqrt{-1}}%
\tint\nolimits_{Lk(p)}\tfrac{a(0,y)}{b(0,y)}dy,
\end{equation*}%
\textit{where }$Lk(p)$\textit{\ is the link of the singularity.}
\end{enumerate}

The work is organized as follows. In Section \ref{Notazioni (Sezione)}\ we
fix some important notations. In Section \ref{Extendable bundles (Sezione)}
we define extendable vector bundles on complex analytic varieties and
extendable sections of extendable bundles. Then, we study in a deeper way
the important case of extendable differentiable differential forms. In
Section \ref{Extendable cohomologies (Sezione)} we define the extendable
cohomology groups and we prove several important results concerning these
groups. In Section \ref{Integration (Sezione)} we define a homomorphism
between extendable and topological cohomology groups of complex analytic
varieties. In Section \ref{Vector bundles (Sezione)} lie the main results of
our work. We use the notion of extendable sections to introduce extendable
connections and extendable Chern classes for complex vector bundles over
complex analytic varieties. Then, we show that these classes represent the
topological Chern classes (defined by means of obstruction theory) via the
homomorphism of integration described in Section \ref{Integration (Sezione)}
(cp. Theorem \ref{Classi Estendibili e Topologiche coincidono copy(1)}).
More precisely, we represent the localizations of the topological Chern
classes by means of the respective localizations of the extendable Chern
classes (cp. Theorem \ref{Classi Estendibili e Topologiche coincidono
copy(2)}). Furthermore, in the compact case, we prove an abstract residue
theorem (cp. Theorem\ \ref{Teor. dei Res.}). Finally, we prove a Camacho-Sad
type index theorem for holomorphic foliations of singular complex varieties
and, under suitable hypotheses, we explicitly compute the residue at
isolated singularities (cp. Theorem \ref{Teorema Camacho-Sad}).

It is in our opinion that the extendable objects we introduced can be
successfully used in order to solve problems of continuous and discrete
holomorphic dynamics in the setting of singular varieties, avoiding the
desingularization processes. Indeed, we think that the theory of extendable
differentiable forms which we have developed lends itself to many uses and
applications. In fact, generalizations similar to the ones of residue
theorem (cp. Theorem \ref{Teor. dei Res.}) and Camacho-Sad index theorem
(cp. Theorem \ref{Teorema Camacho-Sad}) can be achieved in several contexts.

I would like to thank very much Professors F. Bracci and T. Suwa for their
generous help and for their precious advices. I wish to thank also Professor
J. V. Pereira for his useful remarks improving this work. Finally, I want to
thank the anonymous referee for the large amount of mathematical suggestions
that improved this paper.

\section{Main notations\label{Notazioni (Sezione)}}

Let $M$ be a complex differentiable manifold. The sheaves of germs of
differentiable and holomorphic functions on $M$ are denoted by $\mathcal{C}%
_{M}^{\infty }$\ and, respectively, by $\mathcal{O}_{M}$. The real
(holomorphic, antiholomorphic, complexified real) cotangent and tangent
bundles of $M$ are denoted by $T^{\ast }M$ ($\mathbf{T}^{\ast }M$, $\mathbf{%
\bar{T}}^{\ast }M$, $T^{%
\mathbb{C}
\ast }M$) and, respectively, by $TM$ ($\mathbf{T}M$, $\mathbf{\bar{T}}M$, $%
T^{%
\mathbb{C}
}M$). For each $p\in 
\mathbb{N}
$ we denote by $\mathcal{E}_{M}^{p}$ and by $\Omega _{M}^{p}$ the sheaves of
germs of differentiable and, respectively, holomorphic differential $p$%
-forms on $M$.

Recall that an abstract complex analytic variety $X$\ of complex dimension $%
n $ is a second countable, Hausdorff topological space for which there exist
an open covering $\mathcal{C}=\{A_{l}\}_{l\in L}$ and homeomorphisms%
\begin{equation}
F_{l}:A_{l}\rightarrow W_{l}  \label{Notazione Biolomorfismi 1}
\end{equation}%
between the subsets $A_{l}\subseteq X$ and holomorphic subvarieties $%
W_{l}\subseteq U_{l}$ of open sets $U_{l}\subseteq 
\mathbb{C}
^{n_{l}}$ such that for each nonempty intersection $%
A_{(l_{1},l_{2})}=A_{l_{1}}\cap A_{l_{2}}$ the map%
\begin{equation}
F_{(l_{1},l_{2})}:F_{l_{2}}(A_{(l_{1},l_{2})})\rightarrow
F_{l_{1}}(A_{(l_{1},l_{2})})  \label{Notazione Biolomorfismi 2}
\end{equation}%
defined by $F_{(l_{1},l_{2})}=F_{l_{1}}\circ
F_{l_{2}}^{-1}|_{F_{l_{2}}(A_{(l_{1},l_{2})})}$ is a biholomorphism such
that the regular part of $X$ is endowed with a structure of a complex
manifold of complex dimension $n$. A covering as $\mathcal{C}$\ is a \emph{%
coordinate open covering of }$X$ or an \emph{atlas of }$X$. Sometimes, to
make explicit all the data carried by an atlas $\mathcal{C}$, we write%
\begin{equation}
\mathcal{C}=\{(A_{l},n_{l},U_{l},W_{l},F_{l})\}_{l\in L}.
\label{Notazione Atlante}
\end{equation}

Let $X$ be an abstract finite dimensional complex analytic variety. The
singular locus and the regular part of $X$ will be denoted by $Sing(X)$ and,
respectively, by either $X^{Reg}$ or $X^{\prime }$. Recall that $Sing(X)$\
is a complex analytic subvariety of $X$ and that $X^{Reg}$\ is an open and
dense subset of $X$. The maximal atlas of $X$ will be denoted by $\mathcal{A}%
=\{A_{i}\}_{i\in I}$ and, given any $x\in X$, the set $\{i\in I:A_{i}\ni
x\}\subset I$ will be denoted by $I(x)$. Finally, the sheaves of germs of
differentiable and holomorphic functions on $X$ will be denoted by $\mathcal{%
C}_{X}^{\infty }$ and, respectively, by $\mathcal{O}_{X}$.

A\ finite dimensional complex analytic variety $X$ is a locally compact and
paracompact topological space.

\begin{lemma}
\label{Lemma su Topologia e Ricoprimenti}Let $X$ be a finite dimensional
complex analytic variety and $\mathcal{V}=\{V_{j}\}_{j\in J}$ an open
covering of $X$. Then

\begin{enumerate}
\item There exists an open covering $\mathcal{V}^{\ast }=\{V_{j}^{\ast
}\}_{j\in J}$\ of $X$ whose set of indices is still $J$ and such that for
any $j\in J$ it holds $\overline{V_{j}^{\ast }}\subseteq V_{j}$.

\item There exists an open covering $\mathcal{V}^{\bullet }=\{V_{\lambda
}^{\bullet }\}_{\lambda \in \Lambda }$ of $X$ refining $\mathcal{V}$ and
such that for each $\tilde{\lambda}\in \Lambda $ there is a finite subset of
indices $\Lambda (\tilde{\lambda})\subset \Lambda $ such that $\overline{%
V_{\lambda }^{\bullet }}\cap V_{\tilde{\lambda}}^{\bullet }\neq \emptyset $
if and only if $\lambda \in \Lambda (\tilde{\lambda})$.
\end{enumerate}
\end{lemma}

\begin{proof}
$X$ is a locally compact and paracompact topological space. So, the results
follow from General Topology (cp. \cite{Checcucci Tognoli Vesentini}).
\end{proof}

Finally, for a general reference on complex analytic varieties, see \cite%
{Gunning}, Vol. II.

\section{Extendable bundles\label{Extendable bundles (Sezione)}}

\subsection{Extendable vector bundles\label{Extendable vector bundles
(Sottosezione)}}

We begin with the following definition.

\begin{definition}
\label{Def. Fibrato Estendibile}Let $X$ be an abstract finite dimensional
complex analytic variety and $E^{\prime }\rightarrow X^{\prime }$ a
differentiable real (complex) vector bundle over $X^{\prime }$. We say that $%
E^{\prime }$ is $\mathcal{S}_{E^{\prime }}$\emph{-extendable} if there
exists a coherent sheaf $\mathcal{S}_{E^{\prime }}$ of $\mathcal{C}%
_{X}^{\infty }$-modules over $X$ such that $\mathcal{S}_{E^{\prime
}}|_{X^{\prime }}=\mathcal{C}_{X}^{\infty }(E^{\prime })$. A sheaf as $%
\mathcal{S}_{E^{\prime }}$\ is \emph{associated with} $E^{\prime }$.
\end{definition}

An other definition will be also necessary.

\begin{definition}
\label{Def. Sezioni Estendibili}Let $X$ be an abstract finite dimensional
complex analytic variety and $E^{\prime }\rightarrow X^{\prime }$ an $%
\mathcal{S}_{E^{\prime }}$-extendable differentiable real (complex) vector
bundle. A section $s^{\prime }\in \mathcal{S}_{E^{\prime }}|_{X^{\prime
}}(X^{\prime })$ of $E^{\prime }$ is $\mathcal{S}_{E^{\prime }}$\emph{%
-extendable} if there exists a section $s\in \mathcal{S}_{E^{\prime }}(X)$
such that $s^{\prime }=s|_{X^{\prime }}$. A section as $s$ is an $\mathcal{S}%
_{E^{\prime }}$\emph{-extension of }$s^{\prime }$.
\end{definition}

As a matter of notations and terminologies, the set%
\begin{equation}
\Gamma _{ext}(X^{\prime },E^{\prime })=\{s^{\prime }\in \Gamma (X^{\prime
},E^{\prime }):s^{\prime }\emph{is\ }\mathcal{S}_{E^{\prime }}\emph{%
-extendable}\},  \label{Spazio delle Sezioni Est.}
\end{equation}%
also denoted by $\Gamma _{ext}(E^{\prime })=\Gamma _{ext}(X^{\prime
},E^{\prime })$, is called \emph{space of }$\mathcal{S}_{E^{\prime }}$\emph{%
-extendable differentiable sections of }$E^{\prime }$. The sheaf of germs of 
$\mathcal{S}_{E^{\prime }}$-extendable differentiable sections of $E^{\prime
}$ is denoted by $_{ext}\mathcal{E}^{\prime }$.

We present a simple but fundamental example.

\begin{example}
\label{Esempio E-->X}Let $X$ be an abstract finite dimensional complex
analytic variety, $E\rightarrow X$ a differentiable real (complex) vector
bundle over the whole of $X$ and $\mathcal{E}$ the sheaf of germs of
differentiable sections of $E$. Then the bundle $E^{\prime }=E|_{X^{\prime
}}\rightarrow X^{\prime }$ is $\mathcal{E}$-extendable.

Let $s:X\rightarrow E$ be differentiable global section of $E$. Then $%
s^{\prime }=s|_{X^{^{\prime }}}$ is an $\mathcal{E}$-extendable section of $%
E^{\prime }$.
\end{example}

Less trivial examples of extendable vector bundles will be discussed in
Subsection \ref{Extendable forms (Sottosezione)}.

\begin{notation}
\label{Fasci (Notazioni)}Let $X$ be an abstract finite dimensional complex
analytic variety and $\mathcal{S}$ a sheaf\ over $X$. Given an atlas $%
\{(A_{k},n_{k},U_{k},W_{k},F_{k})\}_{k\in K}$ of $X$ (cp. Section \ref%
{Notazioni (Sezione)} and, in particular, (\ref{Notazione Atlante})), we
simply denote by $\mathcal{S}_{k}$ the sheaf $(F_{k})_{\ast }(\mathcal{S}%
|_{A_{k}})$ over $F_{k}(A_{k})$.
\end{notation}

Let $X$ be an abstract finite dimensional complex analytic variety and $%
\mathcal{A}=\{A_{i}\}_{i\in I}$\ the maximal atlas of $X$ (cp. Section \ref%
{Notazioni (Sezione)}). For any $i\in I$ write $A_{i}^{\prime
}=A_{i}\setminus Sing(X)$. Let $E^{\prime }\rightarrow X^{\prime }$ be an $%
\mathcal{S}_{E^{\prime }}$-extendable differentiable real (complex) vector
bundle over $X^{\prime }$. By the very definition of extendable bundle (cp.
Definition\ \ref{Def. Fibrato Estendibile}), the sheaf $\mathcal{S}%
_{E^{\prime }}$ is coherent. So, for any $x\in X$ there exists an index $%
l\in I(x)$ such that the restriction $\mathcal{S}_{E^{\prime }}|_{A_{l}}$ of 
$\mathcal{S}_{E^{\prime }}$\ at $A_{l}$ is generated by a finite number of
sections and such that the sequence%
\begin{equation*}
(\mathcal{C}_{X}^{\infty })^{\nu _{l}}|_{A_{l}}\rightarrow \mathcal{S}%
_{E^{\prime }}|_{A_{l}}\rightarrow 0,
\end{equation*}%
with $\nu _{l}\in 
\mathbb{N}
$, is exact. As a note, the number $\nu _{l}\in 
\mathbb{N}
$ is, in general, bigger than the Zariski dimension of the germ $X_{x}$ of $%
X $ at $x$. Then the sequence%
\begin{equation*}
((\mathcal{C}_{X}^{\infty })^{\nu _{l}}|_{A_{l}})_{l}\rightarrow (\mathcal{S}%
_{E^{\prime }}|_{A_{l}})_{l}\rightarrow 0
\end{equation*}%
is also exact. Now, taking into account the following exact sequence of
sheaves $(\mathcal{C}_{%
\mathbb{C}
^{n_{l}}}^{\infty })^{\nu _{l}}|_{U_{l}}\rightarrow ((\mathcal{C}%
_{X}^{\infty })^{\nu _{l}}|_{A_{l}})_{l}\rightarrow 0$, we get a diagram of
surjective maps%
\begin{equation}
\begin{array}{c}
(\mathcal{C}_{%
\mathbb{C}
^{n_{l}}}^{\infty })^{\nu _{l}}|_{U_{l}} \\ 
\downarrow \\ 
((\mathcal{C}_{X}^{\infty })^{\nu _{l}}|_{A_{l}})_{l} \\ 
\downarrow \\ 
0%
\end{array}%
\rightarrow (\mathcal{S}_{E^{\prime }}|_{A_{l}})_{l}\rightarrow 0
\label{Trivializing Extension}
\end{equation}%
In particular, the map $\zeta _{l}:(\mathcal{C}_{%
\mathbb{C}
^{n_{l}}}^{\infty })^{\nu _{l}}|_{U_{l}}\rightarrow (\mathcal{S}_{E^{\prime
}}|_{A_{l}})_{l}$ is surjective.

Next, consider the restriction map $\epsilon _{l}:\mathcal{S}_{E^{\prime
}}|_{A_{l}}\rightarrow \mathcal{S}_{E^{\prime }}|_{A_{l}^{\prime }}$ and the
map%
\begin{equation}
\varepsilon _{l}:(\mathcal{S}_{E^{\prime }}|_{A_{l}})_{l}\rightarrow (%
\mathcal{S}_{E^{\prime }}|_{A_{l}^{\prime }})_{l}
\label{Restrizione tra fasci}
\end{equation}%
induced by $\epsilon _{l}$. These maps are not surjective in general.
However, if $a^{\prime }\in \func{Im}(\epsilon _{l})$, then there exists $%
\tilde{a}\in (\mathcal{C}_{%
\mathbb{C}
^{n_{l}}}^{\infty })^{\nu _{l}}|_{U_{l}}$ such that $a_{l}^{\prime
}=\varepsilon _{l}\circ \zeta _{l}(\tilde{a})$, with $a_{l}^{\prime
}=[F_{l}]_{\ast }(a^{\prime })$.

This happens, for example, in the case of extendable sections. Indeed, if $%
s^{\prime }\in \mathcal{S}_{E^{\prime }}|_{X^{\prime }}(X^{\prime })$ is an $%
\mathcal{S}_{E^{\prime }}$-extendable section of $E^{\prime }$ and if $s\in 
\mathcal{S}_{E^{\prime }}(X)$ is an $\mathcal{S}_{E^{\prime }}$-extension of 
$s^{\prime }$, then $s^{\prime }$ gives rise to an element $s_{l}^{\prime }$
of $(\mathcal{S}_{E^{\prime }}|_{A_{l}^{\prime }})_{l}$ which lies in the
image of $\varepsilon _{l}:(\mathcal{S}_{E^{\prime
}}|_{A_{l}})_{l}\rightarrow (\mathcal{S}_{E^{\prime }}|_{A_{l}^{\prime
}})_{l}$. Denoting by $s_{l}\in (\mathcal{S}_{E^{\prime }}|_{A_{l}})_{l}$
the element determined by $s$, we have $s_{l}^{\prime }=\varepsilon
_{l}(s_{l})$. So, there exists $\tilde{s}_{l}\in (\mathcal{C}_{%
\mathbb{C}
^{n_{l}}}^{\infty })^{\nu _{l}}|_{U_{l}}$ such that%
\begin{equation}
s_{l}^{\prime }=\varepsilon _{l}\circ \zeta _{l}(\tilde{s}_{l}),
\label{Estensione Locale di una Sezione}
\end{equation}%
because of the surjectivity of $\zeta _{l}:(\mathcal{C}_{%
\mathbb{C}
^{n_{l}}}^{\infty })^{\nu _{l}}|_{U_{l}}\rightarrow (\mathcal{S}_{E^{\prime
}}|_{A_{l}})_{l}$.

\begin{remark}
\label{Atlante Associato}Let $X$, $\mathcal{A}$, $E^{\prime }\rightarrow
X^{\prime }$ and $\mathcal{S}_{E^{\prime }}$be as in the above discussion.
We wish to stress the fact that the sheaf $\mathcal{S}_{E^{\prime }}$
determines an atlas $\mathcal{C}_{E^{\prime }}$ of $X$. Namely, the atlas
that, using the above notations (cp. (\ref{Notazione Atlante})), is given by 
$\mathcal{C}_{E^{\prime }}=\{(A_{l},n_{l},U_{l},W_{l},F_{l})\}$. An atlas as 
$\mathcal{C}_{E^{\prime }}$ is an \emph{atlas\ associated with} $E^{\prime }$
or an \emph{atlas of trivializing extensions for }$E^{\prime }$.
\end{remark}

We need to introduce some terminology.

\begin{terminology}
\label{Terminologia}Let $X$, $E^{\prime }\rightarrow X^{\prime }$, $\mathcal{%
S}_{E^{\prime }}$be as above and\emph{\ }$s^{\prime }\in \Gamma
_{ext}(E^{\prime })$ an $\mathcal{S}_{E^{\prime }}$-extendable section of $%
E^{\prime }$. Let $\mathcal{C}_{E^{\prime }}=\{A_{l}\}$\ be an atlas
associated with $E^{\prime }$.

Let $Y$ be a subset of $X$ and $x\in Y$. We say that $s^{\prime }$\ is \emph{%
extended by} $\tilde{s}_{l}$\emph{\ on }$Y$\emph{\ around }$x\in $\emph{\ }$%
Y $ if there exist $A_{l}\in \mathcal{C}_{E^{\prime }}$ such that $A_{l}\ni
x $ and $\tilde{s}_{l}\in (\mathcal{C}_{%
\mathbb{C}
^{n_{l}}}^{\infty })^{\nu _{l}}|_{U_{l}}$ such that $(F_{l}|_{A_{l}^{\prime
}})_{\ast }(s^{\prime }|_{A_{l}^{\prime }})=\varepsilon _{l}\circ \zeta _{l}(%
\tilde{s}_{l})$. Let $Y=A$ be an open subset of $X$. We say that $s^{\prime
} $ is \emph{completely extendable on }$A$ if for each $x\in A$ the open set 
$A_{l}$ contains $A$.
\end{terminology}

We have the following proposition.

\begin{proposition}
\label{Fibr. e Mappe}Let $X_{1}$ and $X_{2}$ be finite dimensional complex
analytic varieties and $h:X_{1}\rightarrow X_{2}$ an analytic map. Let $%
E\rightarrow X_{2}$ be a differentiable complex (real) vector bundle defined
over the whole of $X_{2}$, $s^{\prime }:X_{2}^{\prime }\rightarrow
E|_{X_{2}^{\prime }}$ an extendable section of $E|_{X_{2}^{\prime
}}\rightarrow X_{2}^{\prime }$ and $s\in \mathcal{E}(X_{2})$ an extension of 
$s^{\prime }$. Then $(h|_{X_{1}^{\prime }})^{\ast }(s):X_{1}^{^{\prime
}}\rightarrow h^{\ast }(E)|_{X_{1}^{\prime }}$ is an extendable
differentiable section of $h^{\ast }(E)|_{X_{1}^{\prime }}\rightarrow
X_{1}^{\prime }$.
\end{proposition}

\begin{proof}
$h^{\ast }(\mathcal{E})$ is a locally free sheaf of $\mathcal{C}%
_{X_{2}}^{\infty }$-modules over $X_{1}$, because $\mathcal{E}$, the sheaf
over $X_{2}$ of germs of differentiable sections of $E\rightarrow X_{2}$, is
a locally free sheaf of $\mathcal{C}_{X_{2}}^{\infty }$-modules. Moreover,
since the map $h$ induces a morphism $h^{\ast }:\mathcal{C}_{X_{2}}^{\infty
}\rightarrow \mathcal{C}_{X_{1}}^{\infty }$, the sheaf $h^{\ast }(\mathcal{E}%
)\otimes _{\mathcal{C}_{X_{2}}^{\infty }}\mathcal{C}_{X_{1}}^{\infty }$ is a
well defined locally free sheaf of $\mathcal{C}_{X_{1}}^{\infty }$-modules
over $X_{1}$. Furthermore, $h^{\ast }(\mathcal{E})\otimes _{\mathcal{C}%
_{X_{2}}^{\infty }}\mathcal{C}_{X_{1}}^{\infty }$ is the sheaf of germs of
differentiable sections of the bundle $h^{\ast }(E)\rightarrow X_{1}$. So,
the section $(h|_{X_{1}^{\prime }})^{\ast }(s):X_{1}^{^{\prime }}\rightarrow
h^{\ast }(E)|_{X_{1}^{\prime }}$ is extendable, because the section $t\in
(h^{\ast }(\mathcal{E})\otimes _{\mathcal{C}_{X_{2}}^{\infty }}\mathcal{C}%
_{X_{1}}^{\infty })(X_{1})$ defined by $t=h^{\ast }(s)\otimes 1$ is an
extension of it.
\end{proof}

The hypotheses of Proposition \ref{Fibr. e Mappe} can be weakened. Let $%
E^{\prime }\rightarrow X_{2}^{\prime }$ be an extendable bundle over $%
X_{2}^{\prime }$ and consider its pull back $(h|_{h^{-1}(X_{2}^{\prime
})})^{\ast }(E^{\prime })\rightarrow h^{-1}(X_{2}^{\prime })$ via $%
h|_{h^{-1}(X_{2}^{\prime })}$. A priori, such a bundle is not extendable,
because, by the very definition of extendable bundle, a necessary condition
for $(h|_{h^{-1}(X_{2}^{\prime })})^{\ast }(E^{\prime })$\ to be extendable
is to be defined at least on the whole of $X_{1}^{\prime }$, the regular
part of $X_{1}$ (cp. Definition\ \ref{Def. Fibrato Estendibile}). So, in
order to generalize Proposition \ref{Fibr. e Mappe}, we have to assume that $%
X_{1}^{\prime }\subseteq h^{-1}(X_{2}^{\prime })$, that is%
\begin{equation}
h^{-1}(Sing(X_{2}))\subseteq Sing(X_{1}).
\label{Singolarita' per Tirare Indietro}
\end{equation}

\begin{proposition}
\label{Forme Est. e Mappe}Let $X_{1}$ and $X_{2}$ be finite dimensional
complex analytic varieties and $h:X_{1}\rightarrow X_{2}$ an analytic map
such that $h^{-1}(Sing(X_{2}))\subseteq Sing(X_{1})$. Let $E^{\prime
}\rightarrow X_{2}^{\prime }$ be an $\mathcal{S}_{E^{\prime }}$-extendable
differentiable real (complex) vector bundle and $s^{\prime }:X_{2}^{\prime
}\rightarrow E^{\prime }$ be an $\mathcal{S}_{E^{\prime }}$-extendable
section of $E^{\prime }$. Then the restriction $(h|_{h^{-1}(X_{2}^{\prime
})})^{\ast }(E^{\prime })|_{X_{1}^{\prime }}\rightarrow X_{1}^{\prime }$ of $%
(h|_{h^{-1}(X_{2}^{\prime })})^{\ast }(E^{\prime })\rightarrow
h^{-1}(X_{2}^{\prime })$ at $X_{1}^{\prime }$ is a $(h^{\ast }(\mathcal{S}%
_{E^{\prime }})\otimes _{\mathcal{C}_{X_{2}}^{\infty }}\mathcal{C}%
_{X_{1}}^{\infty })$-extendable vector bundle and $(h|_{X_{1}^{\prime
}})^{\ast }(s^{\prime }):X_{1}^{^{\prime }}\rightarrow h^{\ast }(E^{\prime
})|_{X_{1}^{\prime }}$ is a $(h^{\ast }(\mathcal{S}_{E^{\prime }})\otimes _{%
\mathcal{C}_{X_{2}}^{\infty }}\mathcal{C}_{X_{1}}^{\infty })$-extendable
differentiable section of $h^{\ast }(E^{\prime })|_{X_{1}^{\prime
}}\rightarrow X_{1}^{\prime }$.
\end{proposition}

\begin{proof}
Since $\mathcal{S}_{E^{\prime }}$ is a coherent sheaf of $\mathcal{C}%
_{X_{2}}^{\infty }$-modules over $X_{2}$ whose restriction at $X_{2}^{\prime
}$ is a locally free sheaf of $\mathcal{C}_{X_{2}}^{\infty }$-modules, its
pull back $h^{\ast }(\mathcal{S}_{E^{\prime }})$ via $h$ is a coherent sheaf
of $\mathcal{C}_{X_{2}}^{\infty }$-modules over $X_{1}$ whose restrictions
at $h^{-1}(X_{2})$ and at $X_{1}^{\prime }$ are a locally free sheaves of $%
\mathcal{C}_{X_{2}}^{\infty }$-modules. Moreover, the map $h$ induces a
morphism $h^{\ast }:\mathcal{C}_{X_{2}}^{\infty }\rightarrow \mathcal{C}%
_{X_{1}}^{\infty }$. So, $h^{\ast }(\mathcal{S}_{E^{\prime }})\otimes _{%
\mathcal{C}_{X_{2}}^{\infty }}\mathcal{C}_{X_{1}}^{\infty }$ is a well
defined coherent sheaf of $\mathcal{C}_{X_{1}}^{\infty }$-modules over $%
X_{1} $ whose restrictions at $h^{-1}(X_{2})$ and at $X_{1}^{\prime }$ are a
locally free sheaves of $\mathcal{C}_{X_{1}}^{\infty }$-modules. Then the
vector bundle $(h|_{h^{-1}(X_{2})})^{\ast }(E^{\prime })|_{X_{1}^{\prime
}}\rightarrow X_{1}^{\prime }$ is extendable is $(h^{\ast }(\mathcal{S}%
_{E^{\prime }})\otimes _{\mathcal{C}_{X_{2}}^{\infty }}\mathcal{C}%
_{X_{1}}^{\infty })$-extendable, because the restriction of $(h^{\ast }(%
\mathcal{S}_{E^{\prime }})\otimes _{\mathcal{C}_{X_{2}}^{\infty }}\mathcal{C}%
_{X_{1}}^{\infty })$ at $X_{1}^{\prime }$ coincides with the locally free
sheaf of $\mathcal{C}_{X_{1}}^{\infty }$-modules of germs of differentiable
sections of $(h|_{h^{-1}(X_{2})})^{\ast }(E^{\prime })|_{X_{1}^{\prime }}$.
Finally, the section $(h|_{X_{1}^{\prime }})^{\ast }(s):X_{1}^{^{\prime
}}\rightarrow h^{\ast }(E^{\prime })|_{X_{1}^{\prime }}$ is extendable,
because it admits the extension $h^{\ast }(s)\otimes 1\in (h^{\ast }(%
\mathcal{S}_{E^{\prime }})\otimes _{\mathcal{C}_{X_{2}}^{\infty }}\mathcal{C}%
_{X_{1}}^{\infty })(X_{1})$, with $s\in \mathcal{S}_{E^{\prime }}(X_{2})$
any extension of $s^{\prime }:X_{2}^{\prime }\rightarrow E^{\prime }$.
\end{proof}

\subsection{Extendable differential forms\label{Extendable forms
(Sottosezione)}}

In this subsection we study the extendable vector bundles we are mainly
interested in. Actually, the definition of extendable vector bundles given
in Subsection \ref{Extendable vector bundles (Sottosezione)} (cp. Definition %
\ref{Def. Fibrato Estendibile})\ has been based on such examples.

Let $X$ be an abstract finite dimensional complex analytic variety. We need
the following observations.

\begin{enumerate}
\item Denote by $TX^{\prime }$ and $T^{\ast }X^{\prime }$ the holomorphic
tangent and, respectively, cotangent bundles of the manifold $X^{\prime }$.
Then, in case $Sing(X)\neq \emptyset $, the bundles $TX^{\prime }\ $and $%
T^{\ast }X^{\prime }$ are not the restriction at $X^{\prime }$\ of any
vector bundle defined on the whole of $X$. The same holds for every their
tensor power, for every their non trivial algebraic quotient\ and for every
their vector subbundle. In particular, for any $N,$ $N^{\ast },$ $p\in 
\mathbb{N}
$ the bundles $TX^{\prime \otimes N}\otimes T^{\ast }X^{\prime \otimes
N^{\ast }}$ and $\Lambda ^{p}T^{\ast }X^{\prime }$ are not the restriction
at $X^{\prime }$ of any bundle defined on $X$.

\item Let $\mathbf{T}X$ and $\mathbf{T}^{\ast }X$ be the holomorphic tangent
and, respectively, holomorphic cotangent varieties of $X$. The following
diffeomorphisms (denoted by $\approx $) of real vector bundles hold: $%
TX^{\prime }\approx \mathbf{T}X^{\prime }=\mathbf{T}X|_{X^{\prime }}$ and $%
T^{\ast }X^{\prime }\approx \mathbf{T}^{\ast }X^{\prime }=\mathbf{T}^{\ast
}X|_{X^{\prime }}$ (cp. Proposition 6.2 of Chapter I of \cite{SuwaLibro}).

\item Let $\mathcal{O}_{X}(\mathbf{T}X)$ be the sheaf of germs of
holomorphic vector fields on $X$. Then $\mathcal{O}_{X}(\mathbf{T}X)$ is a
sheaf of $\mathcal{O}_{X}$-modules over $X$ that, even if it is not
necessarily locally free, is always coherent.

\item Denote by $\Omega _{X}$ the sheaf of germs of holomorphic
differentials on $X$. Then $\Omega _{X}$ is a coherent sheaf of $\mathcal{O}%
_{X}$-modules over $X$. Furthermore, in case $X$ is reduced and irreducible, 
$\Omega _{X}$ is locally free if and only if $X$ is regular (see Theorem
8.15 of \cite{Hartshorne}).

\item Fix $p\in 
\mathbb{N}
$. We wish to consider the sheaf over $X$ of germs of holomorphic $p$-form
on $X$. Actually, several different definitions of such a sheaf can be
given. So, for the moment, just choose one of them, call it the sheaf over $%
X $ of germs of holomorphic $p$-form on $X$ and denote it by $\Omega
_{X}^{p} $. Anyway, as it will be clear from the following discussion, to
our purposes, it is not important which of the several candidates has been
chosen. Indeed, what it is really important is the following property common
to any sheaf candidate to be the sheaf of germs of holomorphic $p$-form on $%
X $: it is a sheaf of $\mathcal{O}_{X}$-modules that, even if it is not
necessarily locally free, it is always coherent (cp. \cite{Ferrari 2}).
\end{enumerate}

Let us recall that $\mathcal{C}_{X}^{\infty }$ is the sheaf of germs of
differentiable functions on $X$ (see Section \ref{Notazioni (Sezione)}).

\begin{remark}
\label{Fibr. Tang. Non Est. (Remark)}Let $X$ be an abstract finite
dimensional complex analytic variety. Then $TX^{\prime }\rightarrow
X^{\prime }$ is $(\mathcal{O}_{X}(\mathbf{T}X)\otimes _{\mathcal{O}_{X}}%
\mathcal{C}_{X}^{\infty })$-extendable and $T^{\ast }X^{\prime }\rightarrow
X^{\prime }$ is $(\Omega _{X}\otimes _{\mathcal{O}_{X}}\mathcal{C}%
_{X}^{\infty })$-extendable. Analogously, every tensor power, every non
trivial algebraic quotient\ and every vector subbundle of $TX^{\prime }$ and 
$T^{\ast }X^{\prime }$ is an extendable vector bundle. In particular, for
any $N,$ $N^{\ast },$ $p\in 
\mathbb{N}
$ the bundle $TX^{\prime \otimes N}\otimes T^{\ast }X^{\prime \otimes
N^{\ast }}$ is $((\Omega _{X}^{\otimes N}\otimes _{\mathcal{O}_{X}}\mathcal{O%
}_{X}(\mathbf{T}X)^{\otimes N^{\ast }})\otimes _{\mathcal{O}_{X}}\mathcal{C}%
_{X}^{\infty })$-extendable and the bundle $\Lambda ^{p}T^{\ast }X^{\prime }$
is $(\Omega _{X}^{p}\otimes _{\mathcal{O}_{X}}\mathcal{C}_{X}^{\infty })$%
-extendable.
\end{remark}

We introduce more manageable notations.

\begin{notation}
Let $X$ be an abstract finite dimensional complex analytic variety. For each 
$N,$ $N^{\ast },$ $p\in 
\mathbb{N}
$ the sheaf $(\Omega _{X}^{\otimes N}\otimes _{\mathcal{O}_{X}}\mathcal{O}%
_{X}(\mathbf{T}X)^{\otimes N^{\ast }})\otimes _{\mathcal{O}_{X}}\mathcal{C}%
_{X}^{\infty }$ will be denoted by $\mathcal{S}_{N,N^{\ast }}$ and the sheaf 
$\Omega _{X}^{p}\otimes _{\mathcal{O}_{X}}\mathcal{C}_{X}^{\infty }$ will be
denoted by $\mathcal{S}_{p}$.
\end{notation}

Thus the bundles $TX^{\prime \otimes N}\otimes T^{\ast }X^{\prime \otimes
N^{\ast }}$ and $\Lambda ^{p}T^{\ast }X^{\prime }$ are $\mathcal{S}%
_{N,N^{\ast }}$-extendable and, respectively, $\mathcal{S}_{p}$-extendable.
As a note, both $\mathcal{S}_{0,0}|_{X^{\prime }}$ and $\mathcal{S}%
_{0}|_{X^{\prime }}$ coincide with $\mathcal{C}_{X^{\prime }}^{\infty }=%
\mathcal{C}_{X}^{\infty }|_{X^{\prime }}$.

\begin{example}
\label{Oss. Fibrati Estendibili}Let $X$ be an abstract finite dimensional
complex analytic variety and $E\rightarrow X$ a differentiable real
(complex) vector bundle. Then for each $N,$ $N^{\ast },$ $p\in 
\mathbb{N}
$ the bundle $TX^{\prime \otimes N}\otimes T^{\ast }X^{\prime \otimes
N^{\ast }}\otimes \Lambda ^{p}T^{\ast }X^{\prime }\otimes E|_{X^{\prime }}$
is $(\mathcal{S}_{N,N^{\ast }}\otimes \mathcal{S}_{p}\otimes \mathcal{E})$%
-extendable.
\end{example}

The following definition is justified by the very important role played by
the $\mathcal{S}_{p}$-extendable sections of $\Lambda ^{p}T^{\ast }X^{\prime
}$.

\begin{definition}
\label{Def. Forme Estendibili}Let $X$ be an abstract finite dimensional
complex analytic variety and $p\in 
\mathbb{N}
$. The space%
\begin{equation}
\Gamma _{ext}(\Lambda ^{p}T^{\ast }X^{\prime })=\{\omega \in \Gamma
(X^{\prime },\Lambda ^{p}T^{\ast }X^{\prime }):\omega \text{ }\emph{is\ }%
\mathcal{S}_{p}\emph{-extendable}\}  \label{Spazio delle Forme Est.}
\end{equation}%
of $\mathcal{S}_{p}$-extendable differentiable sections of $\Lambda
^{p}T^{\ast }X^{\prime }$ is called \emph{space of extendable differentiable 
}$p$\emph{-forms on }$X$.
\end{definition}

In what follows, we heavily use Notations of Section \ref{Notazioni
(Sezione)} (see (\ref{Notazione Atlante})) and Subsection \ref{Extendable
vector bundles (Sottosezione)}.

Let $X$\ be an abstract finite dimensional complex analytic variety, $p\in 
\mathbb{N}
$ and $\mathcal{C}_{p}=\{(A_{l},n_{l},U_{l},W_{l},F_{l})\}_{l\in L}$ an
atlas of trivializing extension for the $\mathcal{S}_{p}$-extendable bundle $%
\Lambda ^{p}T^{\ast }X^{\prime }$ (cp. Remark \ref{Atlante Associato}). Up
to shrink the open sets of $\mathcal{C}_{p}$, if necessary, we can improve (%
\ref{Trivializing Extension}) and get the following commutative diagram of
surjective maps (for the notations see Section \ref{Notazioni (Sezione)}).%
\begin{equation}
\begin{array}{ccc}
(\mathcal{C}_{%
\mathbb{C}
^{n_{l}}}^{\infty })^{\nu _{l}}|_{U_{l}} & \rightarrow & \mathcal{E}%
_{U_{l}}^{p} \\ 
\downarrow &  & \downarrow \\ 
((\mathcal{C}_{X}^{\infty })^{\nu _{l}}|_{A_{l}})_{l} & \rightarrow & (%
\mathcal{S}_{p}|_{A_{l}})_{l}%
\end{array}
\label{Miglioramento di ''Trivializing Extension''}
\end{equation}%
Indeed, it suffices to observe that for any $x\in X$ the index $l\in I(x)$
can be chosen in such a way that $n_{l}$ verifies the following inequality $%
\tbinom{n_{l}}{p}\leq \nu _{l}$.

Thus, thanks to (\ref{Miglioramento di ''Trivializing Extension''}), (\ref%
{Estensione Locale di una Sezione})\ can be also improved. Namely, $\omega
\in \Gamma (\Lambda ^{p}T^{\ast }X^{\prime })$ is extendable if and only if
for each $x\in X$ there are $A_{l}\in \mathcal{C}_{p}$ and%
\begin{equation}
\tilde{\omega}_{l}\in \Gamma (\Lambda ^{p}T^{\ast }U_{l})
\label{Estensione Locale di una Forma (Candidata)}
\end{equation}%
such that%
\begin{equation}
\omega |_{A_{l}^{\prime }}=[F_{l}|_{A_{l}^{\prime }}]^{\ast }(\tilde{\omega}%
_{l}).  \label{Estensione Locale di una Forma}
\end{equation}

We have the following remark (cp. \cite{Perrone}).

\begin{remark}
\label{Estensione su Ogni Atlante}Let $X$\ be a finite dimensional complex
analytic variety and $\mathcal{C}=\{A_{l}\}_{l\in L}$ an atlas of $X$. Then $%
\omega \in \Gamma (\Lambda ^{p}T^{\ast }X^{\prime })$ is extendable on $X$
if and only if it is extendable on $A_{l}\in \mathcal{C}$ for any $l\in L$.
\end{remark}

We explicitly note some facts concerning extendable forms.

\begin{remark}
\label{Fatti sulle Forme Est. (Limitatezza)}Let $X$ be a finite dimensional
complex analytic variety and $\omega \in \Gamma _{ext}(\oplus _{p\in 
\mathbb{N}
}\Lambda ^{p}T^{\ast }X^{\prime })$ an extendable form on $X$. Then for any
point $x\in X$, even singular, there exists a neighborhood $A_{x}$ of $x$
such that $\omega |_{A_{x}\cap X^{Reg}}$ is bounded. This follows from the
very definition of extendable sections (cp. Definition \ref{Def. Forme
Estendibili}) and from (\ref{Miglioramento di ''Trivializing Extension''}), (%
\ref{Estensione Locale di una Forma (Candidata)}), (\ref{Estensione Locale
di una Forma}).
\end{remark}

\begin{remark}
\label{Fatti sulle Forme Est. (Non si Tira Indietro)}Let $X_{1},$ $X_{2}$ be
finite dimensional complex analytic varieties and $h:X_{1}\rightarrow X_{2}$
an analytic map. If $\omega \in \Gamma _{ext}(\oplus _{p\in 
\mathbb{N}
}\Lambda ^{p}T^{\ast }X_{2}^{\prime })$ is an extendable form on $X_{2}$,
then its pull back $(h|_{h^{-1}(X_{2}^{\prime })\cap X_{1}^{\prime }})^{\ast
}(\omega )$\ is not, in general, an extendable form on $X_{1}$. Indeed, in
order to have the extensibility of $(h|_{h^{-1}(X_{2}^{\prime })\cap
X_{1}^{\prime }})^{\ast }(\omega )$, Condition (\ref{Singolarita' per Tirare
Indietro})\ must be verified (cp. Proposition \ref{Forme Est. e Mappe}).
\end{remark}

Let $X$ be an abstract finite dimensional complex analytic variety and fix $%
p\in 
\mathbb{N}
$. Let $Z$ be the closure of a non empty open set of $X$ which is also a
polyhedron of $X$ and denote by $\mathfrak{i}:Z\hookrightarrow X$\ the
inclusion. As a matter of notations, set $Sing(Z)=Z\cap Sing(X)$ and $%
Z^{\prime }=Z\cap X^{\prime }$. Then%
\begin{equation}
\mathfrak{i}^{-1}(Sing(X))\subseteq Sing(Z).
\label{Singolarita' per Tirare Indietro a Z}
\end{equation}%
The set%
\begin{equation}
\Gamma _{ext}(\Lambda ^{p}T^{\ast }X^{\prime })_{Z}=\{\omega \in \Gamma
_{ext}(\Lambda ^{p}T^{\ast }X^{\prime }):(\mathfrak{i}|_{Z^{\prime }})^{\ast
}(\omega )=0\}  \label{Spazio delle Forme Est. Nulle su Z}
\end{equation}%
is called \emph{space of extendable (differentiable) differential }$p$\emph{%
-forms (on }$X$\emph{) vanishing on }$Z$. It is easy to prove that $\Gamma
_{ext}(\Lambda ^{p}T^{\ast }X^{\prime })_{Z}$ is a subspace of $\Gamma
_{ext}(\Lambda ^{p}T^{\ast }X^{\prime })$.

Next, consider the vector bundle $(\mathfrak{i}|_{Z^{\prime }})^{\ast
}(\Lambda ^{p-1}T^{\ast }X^{\prime })\rightarrow Z^{\prime }$ and look at (%
\ref{Singolarita' per Tirare Indietro a Z}). Then, arguing as in the proof
of Proposition \ref{Forme Est. e Mappe}, it can be proved that such a bundle
is $(\mathfrak{i}^{\ast }(\mathcal{S}_{p-1})\otimes _{\mathcal{C}%
_{X}^{\infty }}\mathcal{C}_{Z}^{\infty })$-extendable. The set%
\begin{equation}
\Gamma _{ext}(\Lambda ^{p}\mathfrak{i})=\Gamma _{ext}(\Lambda ^{p}T^{\ast
}X^{\prime })\oplus \Gamma _{ext}((\mathfrak{i}|_{Z^{\prime }})^{\ast
}(\Lambda ^{p-1}T^{\ast }X^{\prime }))
\label{Spazio delle Forme Est. Relative a (X,Z)}
\end{equation}%
is called \emph{space of extendable (differentiable) differential }$p$\emph{%
-forms relative to the pair }$(X,Z)$.

\begin{remark}
\label{Somma di forme estendibili e' estendibile}Let $X$, $Z$ and $\mathfrak{%
i}$ be as above. It is easy to check that $\Gamma _{ext}(\oplus _{p\in 
\mathbb{N}
}\Lambda ^{p}T^{\ast }X^{\prime })$, $\Gamma _{ext}(\oplus _{p\in 
\mathbb{N}
}\Lambda ^{p}T^{\ast }X^{\prime })_{Z}$\ and $\Gamma _{ext}(\oplus _{p\in 
\mathbb{N}
}\Lambda ^{p}\mathfrak{i})$\ are complex vector spaces endowed with a
structure of a graded algebra.
\end{remark}

Let $X$ be a finite dimensional complex analytic variety. The two following
short observations are in order.

\begin{enumerate}
\item $T^{%
\mathbb{C}
}X^{\prime \otimes N}$ is an extendable vector bundle. Indeed, $T^{%
\mathbb{C}
}X^{\prime }$ splits as $T^{%
\mathbb{C}
}X^{\prime }=\mathbf{T}X^{\prime }\oplus \mathbf{\bar{T}}X^{\prime }$ and
both $\mathbf{T}X^{\prime }$ and $\mathbf{\bar{T}}X^{\prime }$ are
extendable bundles.

\item Up to write \emph{holomorphic} and $\mathcal{O}_{X}$ instead of \emph{%
differentiable} and, respectively, $\mathcal{C}_{X}^{\infty }$, the above
discussion can be repeated word by word in the holomorphic category. The
space and the sheaf of germs of $\mathcal{S}_{E^{\prime }}$-extendable
holomorphic sections of an $\mathcal{S}_{E^{\prime }}$-extendable
holomorphic complex vector bundle $E^{\prime }\rightarrow X^{\prime }$ will
be denoted by $\digamma _{ext}(X^{\prime },E^{\prime })$ and, respectively,
by $_{ext}\mathcal{O}_{X}(E^{\prime })$.
\end{enumerate}

\section{Extendable cohomologies\label{Extendable cohomologies (Sezione)}}

\subsection{Extendable cohomology groups\label{Extendable cohomology groups
(Sottosezione)}}

In this subsection we introduce the extendable cohomology groups of a
complex analytic variety.

\begin{enumerate}
\item Let $X$ be an abstract finite dimensional complex analytic variety.
Then%
\begin{equation*}
d^{p}:\Gamma _{ext}(\Lambda ^{p}T^{\ast }X^{\prime })\rightarrow \Gamma
_{ext}(\Lambda ^{p+1}T^{\ast }X^{\prime }),
\end{equation*}%
the restriction at $\Gamma _{ext}(\Lambda ^{p}T^{\ast }X^{\prime })$\ of the 
$p^{th}$\ exterior differential $d^{p}:\Gamma (\Lambda ^{p}T^{\ast
}X^{\prime })\rightarrow \Gamma (\Lambda ^{p+1}T^{\ast }X^{\prime })$, is
well defined for any $p\in 
\mathbb{N}
$. So, $\Gamma _{ext}(\oplus _{p\in 
\mathbb{N}
}\Lambda ^{p}T^{\ast }X^{\prime })$ endowed with $d=\oplus _{p\in 
\mathbb{N}
}d^{p}$\ is a cochains complex, because $d\circ d=0$.

\begin{proof}
Let $\mathcal{C}=\{A_{l}\}_{l\in L}$ be an atlas of trivializing extension
for $\oplus _{p\in 
\mathbb{N}
}\Lambda ^{p}T^{\ast }X^{\prime }$. Given $\omega \in \Gamma _{ext}(\Lambda
^{p}T^{\ast }X^{\prime })$, let $\tilde{\omega}$ be an extension of $\omega $
on $A_{l}$. Then $d^{p}\left( \tilde{\omega}\right) $ is an extension of $%
d^{p}\left( \omega \right) $ on $A_{l}$, because $d^{p}$ commutes with the
pull back operators (see (\ref{Miglioramento di ''Trivializing Extension''}%
), (\ref{Estensione Locale di una Forma (Candidata)}), (\ref{Estensione
Locale di una Forma})).
\end{proof}

\item Let $X$ be an abstract finite dimensional complex analytic variety, $Z$
the closure of a non empty open set of $X$ which is also a polyhedron of $X$
and $\mathfrak{i}:Z\hookrightarrow X$\ the inclusion. Set $Sing(Z)=Z\cap
Sing(X)$ and $Z^{\prime }=Z\cap X^{\prime }$ and let $\omega |_{Z}=(%
\mathfrak{i}|_{Z^{\prime }})^{\ast }(\omega )$ denote the pull back of $%
\omega \in \Gamma (\oplus _{p\in 
\mathbb{N}
}\Lambda ^{p}T^{\ast }X^{\prime })$ to $Z^{\prime }$. Then

\begin{enumerate}
\item The restriction of $d^{p}$ at $\Gamma _{ext}(\Lambda ^{p}T^{\ast
}X^{\prime })_{Z}$ is well defined, because $d$ is a local operator.

\item For each $p\in 
\mathbb{N}
$ let $d_{\mathfrak{i}}^{p}:\Gamma _{ext}(\Lambda ^{p}\mathfrak{i}%
)\rightarrow \Gamma _{ext}(\Lambda ^{p+1}\mathfrak{i})$ be the operator
given by%
\begin{equation}
d_{\mathfrak{i}}^{p}(\varkappa ,\upsilon )=(d^{p}\varkappa ,\varkappa
|_{Z}-d^{p-1}\upsilon ).  \label{Differenziale del Cilindro Mappante}
\end{equation}%
Set $d_{\mathfrak{i}}=\oplus _{p\in 
\mathbb{N}
}d_{\mathfrak{i}}^{p}$. A straightforward computation shows that $d_{%
\mathfrak{i}}\circ d_{\mathfrak{i}}=0$. So, $\Gamma _{ext}(\oplus _{p\in 
\mathbb{N}
}\Lambda ^{p}\mathfrak{i})$ endowed with $d_{\mathfrak{i}}$\ is a cochains
complex. As a note, a $d_{\mathfrak{i}}$-closed element of $\Gamma
_{ext}(\oplus _{p\in 
\mathbb{N}
}\Lambda ^{p}\mathfrak{i})$\ corresponds to a $d$-closed element of $\Gamma
_{ext}(\oplus _{p\in 
\mathbb{N}
}\Lambda ^{p}T^{\ast }X)$ whose restriction at $Z$ is exact.
\end{enumerate}
\end{enumerate}

\begin{definition}
Let $X$ be an abstract finite dimensional complex analytic variety, $Z$ the
closure of a non empty open set which is also a polyhedron of $X$ and $%
\mathfrak{i}:Z\hookrightarrow X$ the inclusion. Fix $p\in 
\mathbb{N}
$.

\begin{enumerate}
\item Set $Z_{ext}^{p}\left( X\right) =\ker (d^{p}|_{\Gamma _{ext}(\Lambda
^{p}T^{\ast }X^{\prime })})$ and $B_{ext}^{p}\left( X\right) =\func{Im}%
(d^{p-1}|_{\Gamma _{ext}(\Lambda ^{p-1}T^{\ast }X^{\prime })})$. The group%
\begin{equation*}
H_{ext}^{p}\left( X\right) =\frac{Z_{ext}^{p}\left( X\right) }{%
B_{ext}^{p}\left( X\right) }
\end{equation*}%
is \emph{the }$p^{th}$\emph{\ extendable cohomology group of }$X$.

\item Set $Z_{ext}^{p}(X)_{Z}=\ker (d^{p}|_{\Gamma _{ext}(\Lambda
^{p}T^{\ast }X^{\prime })_{Z}})$ and $B_{ext}^{p}(X)_{Z}=\func{Im}%
(d^{p-1}|_{\Gamma _{ext}(\Lambda ^{p-1}T^{\ast }X^{\prime })_{Z}})$. The
group%
\begin{equation*}
H_{ext}^{p}(X)_{Z}=\frac{Z_{ext}^{p}(X)_{Z}}{B_{ext}^{p}(X)_{Z}}
\end{equation*}%
is \emph{the }$p^{th}$\emph{\ extendable cohomology group of }$X$\emph{\
vanishing on }$Z$.

\item Set $Z_{ext}^{p}\left( X,Z\right) =\ker (d_{\mathfrak{i}}^{p})$ and $%
B_{ext}^{p}\left( X,Z\right) =\func{Im}(d_{\mathfrak{i}}^{p-1})$. The group%
\begin{equation*}
H_{ext}^{p}\left( X,Z\right) =\frac{Z_{ext}^{p}\left( X,Z\right) }{%
B_{ext}^{p}\left( X,Z\right) }
\end{equation*}%
is \emph{the }$p^{th}$\emph{\ extendable cohomology group relative to the
pair }$(X,Z)$.
\end{enumerate}
\end{definition}

The following example is an adjustment of a real analytic example given by
Bloom and Herrera (see\ pages 287-288 of \cite{Bloom-Herrera}).

\begin{example}
\label{Es. Bloom Herrera Complesso}Let $z$ denote the coordinate on $%
\mathbb{C}
$\ and $(z_{1},$ $z_{2})$ the coordinates on $%
\mathbb{C}
^{2}$. Consider the map%
\begin{equation*}
\begin{array}{rccl}
f: & 
\mathbb{C}
& \rightarrow & 
\mathbb{C}
^{2} \\ 
& z & \mapsto & (z^{5},\text{ }z^{6}+z^{7})%
\end{array}%
\end{equation*}%
and let $B$ be a neighborhood of $0$ in $%
\mathbb{C}
$ such that $X=f(B)$ is an irreducible complex analytic variety. Recall
that, since the complex dimension of $X$ is $1$, any holomorphic
differential $2$-form of type $(2,0)$ defined on $X^{Reg}$ is identically
zero. Let $\omega $ be a holomorphic differential $1$-form of type $(1,0)$
defined on a neighborhood of $0$ in $%
\mathbb{C}
^{2}$ and not identically zero on $X^{Reg}$. Since $d(\omega )=\partial
(\omega )+\bar{\partial}(\omega )$, we have $d(\omega )=0$ on $X^{Reg}$.
Indeed, on one hand, $\bar{\partial}(\omega )=0$, because $\omega $ is
holomorphic, and, on the other hand, $\partial (\omega )=0$, because $%
\partial (\omega )$ is a holomorphic differential $2$-form of type $(2,0)$
on $X^{Reg}$. Let $\mathbf{\omega }\in \mathcal{E}_{%
\mathbb{C}
^{2},0}^{1}$ be the germ at $0$ of $\omega $. Taken $f^{\ast }(\mathbf{%
\omega })\in \mathcal{E}_{%
\mathbb{C}
,0}^{1}$, it results $d(f^{\ast }(\mathbf{\omega }))=0$. Indeed, $f^{\ast }(%
\mathbf{\omega })$ is closed, because it is a form of type $(1,0)$ on $%
\mathbb{C}
$. Then, by Poincar\'{e} lemma in $%
\mathbb{C}
$, there exists $\mathbf{h}\in \mathcal{E}_{%
\mathbb{C}
,0}^{0}$ such that $f^{\ast }(\mathbf{\omega })=d(\mathbf{h})$. In fact,
since $\mathbf{\omega }$ is holomorphic of type $(1,0)$, such a germ $%
\mathbf{h}$ is in $\mathcal{O}_{%
\mathbb{C}
,0}$. If there existed an element $\mathbf{g}\in \mathcal{E}_{%
\mathbb{C}
^{2},0}^{0}$ not identically vanishing on $X^{Reg}$ and such that $\mathbf{%
\omega }=d(\mathbf{g})$, then we would have $d(f^{\ast }(\mathbf{g}%
))=f^{\ast }(d(\mathbf{g}))=f^{\ast }(\mathbf{\omega })$ and so $\mathbf{h}%
=f^{\ast }(\mathbf{g})+Const.$ . A\ necessary condition for $\mathbf{h}$\ to
be of such a form is that the formal power series of $\mathbf{h}$\ at $0$
can be expressed as a power series in $z^{5}$ and $(z^{6}+z^{7})$. Now,
consider the holomorphic differential $1$-form $\omega =z_{1}dz_{2}$ of type 
$(1,0)$ defined on a neighborhood of $0$ in $%
\mathbb{C}
^{2}$ and not identically zero on $X^{Reg}$. Then, since the resulting $%
\mathbf{h}$ does not have that form, the Poincar\'{e} lemma does not hold.
So, $[\omega ]_{H_{ext}^{1}\left( X\right) }\neq 0$ and $H_{ext}^{1}\left(
X\right) \neq 0$.
\end{example}

The following important remark is in order.

\begin{remark}
Example \ref{Es. Bloom Herrera Complesso} shows that

\begin{enumerate}
\item There exist complex analytic varieties whose extendable cohomology
groups are not trivial

\item A Poincar\'{e} Lemma for extendable differential forms does not hold.
\end{enumerate}
\end{remark}

\subsection{A technical lemma}

The following technical lemma shows the existence of extendable partition of
unity. Notations of Section \ref{Notazioni (Sezione)} will be heavily used.

\begin{lemma}
\label{Esiste Part. d'Unita' "Liscia" copy(1)}Let $X$ be an abstract finite
dimensional complex analytic variety and $\mathcal{B}=\{B_{\beta }\}_{\beta
\in \mathfrak{B}}$ an atlas of $X$ enjoying (2) of Lemma \ref{Lemma su
Topologia e Ricoprimenti}. Then there exists a partition of unity $\{\rho
_{\beta }:X\rightarrow 
\mathbb{R}
\}_{\beta \in \mathfrak{B}}$ subordinated to $\mathcal{B}$ such that for any 
$\beta \in \mathfrak{B}$ the function $r_{\beta }:W_{\beta }\rightarrow 
\mathbb{R}
$ given by%
\begin{equation}
r_{\beta }=\rho _{\beta }|_{A_{\beta }}\circ F_{\beta }^{-1}
\label{Parte dell'Unita' su un Aperto Coordinato}
\end{equation}%
is extendable to a differentiable function $R_{\beta }:U_{\beta }\rightarrow 
\mathbb{R}
$. $\left\{ \rho _{\beta }:X\rightarrow 
\mathbb{R}
\right\} _{\beta \in \mathfrak{B}}$ is an \emph{extendable partition of
unity subordinated to }$\mathcal{B}$. As a note, $\rho _{\beta }|_{X^{\prime
}}:X^{\prime }\rightarrow 
\mathbb{R}
$ is an extendable $0$-form for any $\beta \in \mathfrak{B}$.
\end{lemma}

\begin{proof}
Let $\mathcal{C}=\{C_{\beta }\}_{\beta \in \mathfrak{B}}$, $\mathcal{D}%
=\{D_{\beta }\}_{\beta \in \mathfrak{B}}$ and $\mathcal{E}=\{E_{\beta
}\}_{\beta \in \mathfrak{B}}$\ be open coverings of $X$ all refining $%
\mathcal{B}=$$\left\{ B_{\beta }\right\} _{\beta \in \mathfrak{B}}$. We can
assume without loss of generality that the sets of indices of $\mathcal{C}$, 
$\mathcal{D}$ and $\mathcal{E}$\ coincide with $\mathfrak{B}$, the set of
indices of $\mathcal{B}$. Moreover, we can also assume that for each $\beta
\in \mathfrak{B}$ it holds%
\begin{equation*}
E_{\beta }\subset \overline{E_{\beta }}\subset D_{\beta }\subset \overline{%
D_{\beta }}\subset C_{\beta }\subset \overline{C_{\beta }}\subset B_{\beta }.
\end{equation*}%
These facts follow from Lemma \ref{Lemma su Topologia e Ricoprimenti} (1).
Then, $\mathcal{C}$, $\mathcal{D}$ and $\mathcal{E}$ are coordinate open
coverings of $X$ enjoying (2)\ of Lemma \ref{Lemma su Topologia e
Ricoprimenti} as well as $\mathcal{B}$.

For each $\beta \in \mathfrak{B}$ consider the homeomorphism $F_{\beta
}:B_{\beta }\rightarrow W_{\beta }\subseteq U_{\beta }$ and the images in $%
U_{\beta }$ of $E_{\beta },$ $\overline{E_{\beta }},$ $D_{\beta },$ $%
\overline{D_{\beta }},$ $C_{\beta },$ $\overline{C_{\beta }},$ $B_{\beta }$
via $F_{\beta }$ (cp. (\ref{Notazione Biolomorfismi 1})). Then, there exist
open sets $\boldsymbol{C}_{\beta },$ $\boldsymbol{D}_{\beta },$ $\boldsymbol{%
E}_{\beta }$ of $U_{\beta }$ such that $F_{\beta }\left( C_{\beta }\right) =%
\boldsymbol{C}_{\beta }\cap W_{\beta }$, $F_{\beta }\left( D_{\beta }\right)
=\boldsymbol{D}_{\beta }\cap W_{\beta }$ and $F_{\beta }\left( E_{\beta
}\right) =\boldsymbol{E}_{\beta }\cap W_{\beta }$, because $F_{\beta }\left(
C_{\beta }\right) ,$ $F_{\beta }\left( D_{\beta }\right) ,$ $F_{\beta
}\left( E_{\beta }\right) $ are open sets of $W_{\beta }$ and the topology
of $W_{\beta }=F_{\beta }\left( B_{\beta }\right) $ is induced by Euclidean
topology of $U_{\beta }\subseteq 
\mathbb{C}
^{n_{\beta }}$.

For any $\beta \in \mathfrak{B}$ let $G_{\beta }:U_{\beta }\rightarrow 
\mathbb{R}
$ be a positive, bounded, differentiabile real valued function such that $%
\overline{\boldsymbol{E}_{\beta }}\subseteq supp\left( G_{\beta }\right)
\subseteq \boldsymbol{D}_{\beta }$, denote by $g_{\beta }=G_{\beta
}|_{W_{\beta }}:W_{\beta }\rightarrow 
\mathbb{R}
$ the restriction of $G_{\beta }$ at $W_{\beta }$ and consider the function $%
F_{\beta }^{\ast }\left( g_{\beta }\right) :B_{\beta }\rightarrow 
\mathbb{R}
$ given by $F_{\beta }^{\ast }(g_{\beta })(x)=g_{\beta }\circ F_{\beta }(x)$%
. Since $F_{\beta }^{\ast }(g_{\beta })$ is identically zero away from $%
\overline{D_{\beta }}$, let us agree that on $X\setminus B_{\beta }$ the
symbol $F_{\beta }^{\ast }(g_{\beta })$ denotes the function identically
equal to zero. So, $F_{\beta }^{\ast }\left( g_{\beta }\right) :X\rightarrow 
\mathbb{R}
$ is the function defined on the whole of $X$ given by%
\begin{equation}
F_{\beta }^{\ast }\left( g_{\beta }\right) (x)=\left\{ 
\begin{array}{lll}
g_{\beta }\circ F_{\beta }(x) & \text{ } & for\text{ }x\in B_{\beta } \\ 
0 & \text{ } & for\text{ }x\in X\setminus B_{\beta }%
\end{array}%
\right.  \label{g F-(tirata indietro)}
\end{equation}%
Moreover, let us agree that for any $(\gamma ,\beta )\in \mathfrak{B}\times 
\mathfrak{B}$ such that $B_{\gamma }\cap B_{\beta }=\emptyset $ the symbol $%
F_{\left( \gamma ,\beta \right) }^{\ast }\left( g_{\gamma }\right) $ denotes
the function identically equal to zero (cp. (\ref{Notazione Biolomorfismi 2}%
)).

For each $\beta \in \mathfrak{B}$ let $\rho _{\beta }$ be the function
defined by%
\begin{equation}
\begin{array}{cccl}
\rho _{\beta }: & X & \rightarrow & 
\mathbb{R}
\\ 
& x & \mapsto & \rho _{\beta }\left( x\right) =\tfrac{F_{\beta }^{\ast
}\left( g_{\beta }\right) \left( x\right) }{\tsum\nolimits_{\gamma \in 
\mathfrak{B}}F_{\beta }^{\ast }\circ F_{(\gamma ,\beta )}^{\ast }\left(
g_{\gamma }\right) \left( x\right) }%
\end{array}
\label{Def. della Partizione dell'Unita' Ro (1)}
\end{equation}%
Again we have $\rho _{\beta }=0$ away from $\overline{D_{\beta }}$. So,%
\begin{equation}
\rho _{\beta }\left( x\right) =\left\{ 
\begin{array}{lll}
\tfrac{g_{\beta }\circ F_{\beta }(x)}{\tsum\nolimits_{\gamma \in \mathfrak{B}%
}F_{\beta }^{\ast }\circ F_{(\gamma ,\beta )}^{\ast }\left( g_{\gamma
}\right) \left( x\right) } & \text{ } & for\text{ }x\in B_{\beta } \\ 
0 & \text{ } & for\text{ }x\in X\setminus B_{\beta }%
\end{array}%
\right.  \label{Def. della Partizione dell'Unita' Ro (2)}
\end{equation}

It is easy to check that $\{\rho _{\beta }:X\rightarrow 
\mathbb{R}
\}_{\beta \in \mathfrak{B}}$ is a partition of unity subordinated to $%
\mathcal{B}$. Thus, it only remains to prove that for each $\beta \in 
\mathfrak{B}$ the function $r_{\beta }:W_{\beta }\rightarrow 
\mathbb{R}
$ admits a differentiable extension $R_{\beta }$ defined on the whole of $%
U_{\beta }$.

For this, fix $\beta \in \mathfrak{B}$. It follows from Lemma \ref{Lemma su
Topologia e Ricoprimenti} (2) that there exists a finite subset of indices $%
\mathfrak{B}\left( \beta \right) \subseteq \mathfrak{B}$ such that $%
\overline{B_{\gamma }}\cap B_{\beta }\neq \emptyset $ if and only if $\gamma
\in \mathfrak{B}\left( \beta \right) $. Note that, by (\ref{Parte
dell'Unita' su un Aperto Coordinato})\ and (\ref{Def. della Partizione
dell'Unita' Ro (1)}), in order to have a differentiable extension $R_{\beta
} $ of $r_{\beta }$, it suffices to find a differentiable extension $\mathbb{%
L}_{(\gamma ,\beta )}\left( G_{\gamma }\right) :U_{\beta }\rightarrow 
\mathbb{R}
$ of $F_{(\gamma ,\beta )}^{\ast }\left( g_{\gamma }\right) $ for any $%
\gamma \in \mathfrak{B}\left( \beta \right) $.

Now, for each $\gamma \in \mathfrak{B}\left( \beta \right) $\ consider the
biholomorphism $F_{(\gamma ,\beta )}:F_{\beta }(B_{\beta }\cap B_{\gamma
})\rightarrow F_{\gamma }(B_{\beta }\cap B_{\gamma })$ (see (\ref{Notazione
Biolomorfismi 2})) and note that only two cases are possible. Namely, $%
n_{\beta }\leq n_{\gamma }$ or $n_{\beta }\gneq n_{\gamma }$, where $%
n_{\alpha }=\dim _{%
\mathbb{C}
}U_{\alpha }$ for any $\alpha \in \mathfrak{B}$ (cp. Section \ref{Notazioni
(Sezione)}).

Suppose $n_{\beta }\leq n_{\gamma }$. Then there exist an open subset $%
Q_{\beta \gamma }$ of $U_{\beta }$\ such that $Q_{\beta \gamma }\cap
W_{\beta }=F_{\beta }(B_{\beta }\cap B_{\gamma })$ and a holomorphic
injective map $\mathbb{F}_{(\gamma ,\beta )}:Q_{\beta \gamma }\rightarrow 
\mathbb{F}_{(\gamma ,\beta )}(Q_{\beta \gamma })$ such that $F_{(\gamma
,\beta )}$ is the restriction of $\mathbb{F}_{(\gamma ,\beta )}$ at $%
F_{\beta }(B_{\beta }\cap B_{\gamma })$. Let $O_{\beta \gamma }$ be an open
subset of $U_{\beta }$\ such that $O_{\beta \gamma }\subset Q_{\beta \gamma
} $\ and $O_{\beta \gamma }\cap W_{\beta }=F_{\beta }(B_{\beta }\cap
C_{\gamma })$ and consider the restriction $\mathbb{F}_{(\gamma ,\beta
)}|_{O_{\beta \gamma }}:O_{\beta \gamma }\rightarrow \mathbb{F}_{(\gamma
,\beta )}(O_{\beta \gamma })$ of $\mathbb{F}_{(\gamma ,\beta )}$\ at $%
O_{\beta \gamma }$. Let $\mathbb{L}_{(\gamma ,\beta )}\left( G_{\gamma
}\right) :U_{\beta }\rightarrow 
\mathbb{R}
$ be the map defined by%
\begin{equation}
\mathbb{L}_{(\gamma ,\beta )}\left( G_{\gamma }\right) \left( u\right)
=\left\{ 
\begin{array}{lll}
(\mathbb{F}_{(\gamma ,\beta )}|_{O_{\beta \gamma }})^{\ast }(G_{\gamma }|_{%
\mathbb{F}_{(\gamma ,\beta )}(O_{\beta \gamma })})(u) & \text{ } & for\text{ 
}u\in O_{\beta \gamma } \\ 
0 & \text{ } & for\text{ }u\in U_{\beta }\setminus O_{\beta \gamma }%
\end{array}%
\right.  \label{L (I caso)}
\end{equation}%
It follows from $\overline{\boldsymbol{E}_{\gamma }}\subseteq supp\left(
G_{\gamma }\right) \subseteq \boldsymbol{D}_{\gamma }$ that $\mathbb{L}%
_{(\gamma ,\beta )}\left( G_{\gamma }\right) $ is a well defined, continuous
map that is also is differentiable, because of the holomorphy of $\mathbb{F}%
_{(\gamma ,\beta )}|_{O_{\beta \gamma }}$ and the differentiability of $%
G_{\gamma }$. Furthermore, by its very definition, $\mathbb{L}_{(\gamma
,\beta )}\left( G_{\gamma }\right) $ is such that%
\begin{equation}
\mathbb{L}_{(\gamma ,\beta )}\left( G_{\gamma }\right) |_{F_{\beta }(%
\overline{D_{\beta }})}=F_{(\gamma ,\beta )}^{\ast }\left( g_{\gamma
}\right) |_{F_{\beta }(\overline{D_{\beta }})}.
\label{L (I caso) proprieta'}
\end{equation}%
Then, since $\rho _{\beta }=0$ away from $\overline{D_{\beta }}$, $\mathbb{L}%
_{(\gamma ,\beta )}\left( G_{\gamma }\right) $ is the wanted extension of $%
F_{(\gamma ,\beta )}^{\ast }\left( g_{\gamma }\right) $ and we are done in
the case $n_{\beta }\leq n_{\gamma }$.

Suppose $n_{\beta }\gneq n_{\gamma }$. Then there exist an open subset $%
Q_{\beta \gamma }$ of $U_{\beta }$\ such that $Q_{\beta \gamma }\cap
W_{\beta }=F_{\beta }(B_{\beta }\cap B_{\gamma })$ and a holomorphic
submersion $\mathbb{F}_{(\gamma ,\beta )}:Q_{\beta \gamma }\rightarrow 
\mathbb{F}_{(\gamma ,\beta )}(Q_{\beta \gamma })$ such that $F_{(\gamma
,\beta )}$ is the restriction of $\mathbb{F}_{(\gamma ,\beta )}$ at $%
F_{\beta }(B_{\beta }\cap B_{\gamma })$. Let $O_{\beta \gamma }$ be an open
subset of $U_{\beta }$\ such that $O_{\beta \gamma }\subset Q_{\beta \gamma
} $\ and $O_{\beta \gamma }\cap W_{\beta }=F_{\beta }(B_{\beta }\cap
C_{\gamma })$ and consider the restriction $\mathbb{F}_{(\gamma ,\beta
)}|_{O_{\beta \gamma }}:O_{\beta \gamma }\rightarrow \mathbb{F}_{(\gamma
,\beta )}(O_{\beta \gamma })$ of $\mathbb{F}_{(\gamma ,\beta )}$\ at $%
O_{\beta \gamma }$. By using the local form of submersion, the function $(%
\mathbb{F}_{(\gamma ,\beta )}|_{O_{\beta \gamma }})^{\ast }(G_{\gamma }|_{%
\mathbb{F}_{(\gamma ,\beta )}(O_{\beta \gamma })})$ is an extension of $%
F_{(\gamma ,\beta )}^{\ast }\left( g_{\gamma }\right) |_{F_{\beta }(B_{\beta
}\cap C_{\gamma })}$\ to the whole of $O_{\beta \gamma }$.

Now, let $N_{\beta \gamma }$ be an open subset of $U_{\beta }$ such that $%
N_{\beta \gamma }\subset \overline{N_{\beta \gamma }}\subset O_{\beta \gamma
}$ and $N_{\beta \gamma }\cap W_{\beta }\supset F_{\beta }(B_{\beta }\cap
D_{\gamma })$. Then the restriction $(\mathbb{F}_{(\gamma ,\beta )}|_{%
\overline{N_{\beta \gamma }}})^{\ast }(G_{\gamma }|_{\mathbb{F}_{(\gamma
,\beta )}(\overline{N_{\beta \gamma }})})$ of $G_{\gamma }|_{\mathbb{F}%
_{(\gamma ,\beta )}(O_{\beta \gamma })}\circ \mathbb{F}_{(\gamma ,\beta
)}|_{O_{\beta \gamma }}$ at $\overline{N_{\beta \gamma }}$ is also a
differentiable extension of $F_{(\gamma ,\beta )}^{\ast }\left( g_{\gamma
}\right) |_{\overline{N_{\beta \gamma }}\cap W_{\beta }}$. Let us agree that
the function $(\mathbb{F}_{(\gamma ,\beta )}|_{\overline{N_{\beta \gamma }}%
})^{\ast }(G_{\gamma }|_{\mathbb{F}_{(\gamma ,\beta )}(\overline{N_{\beta
\gamma }})})$ is defined also on $W_{\beta }\setminus N_{\beta \gamma }$ as
the identically vanishing function.

Let $V_{\beta \gamma }$ be an open subset of $U_{\beta }$ that is a
neighborhood of $W_{\beta }\setminus N_{\beta \gamma }$ such that $V_{\beta
\gamma }\subset \overline{V_{\beta \gamma }}\subset U_{\beta }$, $V_{\beta
\gamma }\cap N_{\beta \gamma }\neq \emptyset $ and $\overline{V_{\beta
\gamma }}\cap \boldsymbol{D}_{\beta }=\emptyset $. Then, by Tietze's
extension theorem, the function%
\begin{equation}
(\mathbb{F}_{(\gamma ,\beta )}|_{\overline{V_{\beta \gamma }}\cap \overline{%
N_{\beta \gamma }}})^{\ast }(G_{\gamma }|_{\mathbb{F}_{(\gamma ,\beta )}(%
\overline{V_{\beta \gamma }}\cap \overline{N_{\beta \gamma }})}):(\overline{%
V_{\beta \gamma }}\cap \overline{N_{\beta \gamma }})\cup (W_{\beta
}\setminus N_{\beta \gamma })\rightarrow 
\mathbb{R}
\label{G F-(tirata indietro)}
\end{equation}%
can be extended to a continuous, non identically vanishing function defined
on $(\overline{V_{\beta \gamma }}\cap \overline{N_{\beta \gamma }})\cup 
\overline{V_{\beta \gamma }}$. Moreover, by the approximation theorem (see 
\cite{Steenrod}, 6.7), we can assume without loss of generality that such a
function is also differentiable. Thus, we get a differentiable function
defined on the closed subset $\overline{N_{\beta \gamma }}\cup \overline{%
V_{\beta \gamma }}$ of $U_{\beta }$.

By using the same technique (extension and approximation theorems), we get a
differentiable function $\mathbb{L}_{(\gamma ,\beta )}\left( G_{\gamma
}\right) :U_{\beta }\rightarrow 
\mathbb{R}
$ such that%
\begin{equation}
\mathbb{L}_{(\gamma ,\beta )}\left( G_{\gamma }\right) |_{F_{\beta }(%
\overline{D_{\beta }})}=F_{(\gamma ,\beta )}^{\ast }\left( g_{\gamma
}\right) |_{\overline{D_{\beta }}}.  \label{L (II caso) proprieta'}
\end{equation}%
Then $\mathbb{L}_{(\gamma ,\beta )}\left( G_{\gamma }\right) $ is the wanted
extension of $F_{(\gamma ,\beta )}^{\ast }\left( g_{\gamma }\right) $,
because $\rho _{\beta }=0$ away from $\overline{D_{\beta }}$. We are done
also in the case $n_{\beta }\gneq n_{\gamma }$.

The proof is thereby concluded. Indeed, the function $R_{\beta }:U_{\beta
}\rightarrow 
\mathbb{R}
$ defined by%
\begin{equation}
R_{\beta }(u)=\left\{ 
\begin{array}{lll}
\tfrac{G_{\beta }(u)}{\tsum\nolimits_{\gamma \in \mathfrak{B}}\mathbb{L}%
_{(\gamma ,\beta )}\left( G_{\gamma }\right) (u)} & \text{ } & for\text{ }%
u\in \tbigcup\nolimits_{\gamma \in \mathfrak{B}}O_{\beta \gamma } \\ 
0 & \text{ } & for\text{ }u\in U_{\beta }\setminus
(\tbigcup\nolimits_{\gamma \in \mathfrak{B}}O_{\beta \gamma })%
\end{array}%
\right.  \label{Estensione della Part. dell'Unita'}
\end{equation}%
is a differentiable extension of the function $r_{\beta }:W_{\beta
}\rightarrow 
\mathbb{R}
$ (see (\ref{Parte dell'Unita' su un Aperto Coordinato})).
\end{proof}

Lemma \ref{Esiste Part. d'Unita' "Liscia" copy(1)} has several important
consequences. We begin with the following corollary.

\begin{corollary}
\label{Fasci Fini II}Let $X$ be a finite dimensional complex analytic
variety. Then the sheaves $\mathcal{S}_{N,N^{\ast }}$ and $\mathcal{S}_{p}$\
are fine and soft for any $N,$ $N^{\ast },$ $p\in 
\mathbb{N}
$.
\end{corollary}

As a further consequence of Lemma \ref{Esiste Part. d'Unita' "Liscia"
copy(1)}, we can prove that the spaces of extendable sections of extendable
bundles of Example \ref{Oss. Fibrati Estendibili} are not trivial.

\begin{remark}
\label{Sul Fibr. Tang.}Let $E\rightarrow X$ be a differentiable real
(complex) vector bundle defined over an abstract finite dimensional complex
analytic variety. Then for each $N,$ $N^{\ast },$ $p\in 
\mathbb{N}
$ the bundle $\Gamma _{ext}(TX^{\prime \otimes N}\otimes T^{\ast }X^{\prime
\otimes N^{\ast }}\otimes \Lambda ^{p}T^{\ast }X^{\prime }\otimes
E|_{X^{\prime }})\neq \{0\}$. For this, consider an extendable partition of
unity subordinated to a suitable open covering of $X$ and argue locally.
\end{remark}

Next, let $X$ be an abstract finite dimensional complex analytic variety, $Z$
the closure of a non empty open set that is also a polyhedron of $X$ and $%
\mathfrak{i}:Z\hookrightarrow X$ the inclusion. Let $\alpha ^{p}:\Gamma
_{ext}(\Lambda ^{p}T^{\ast }X^{\prime })_{Z}\rightarrow \Gamma
_{ext}(\Lambda ^{p}\mathfrak{i})$ be the map that to any $\omega \in \Gamma
_{ext}(\Lambda ^{p}T^{\ast }X^{\prime })_{Z}$\ associates $\alpha
^{p}(\omega )=(\omega ,$ $0)$. A straightforward computation shows that $%
\alpha ^{p}$ induces a homomorphism $A^{p}:H_{ext}^{p}(X)_{Z}\rightarrow
H_{ext}^{p}(X,Z)$ in cohomology. Furthermore, the following proposition,
whose proof is a direct consequence of Lemma \ref{Esiste Part. d'Unita'
"Liscia" copy(1)}, holds.

\begin{proposition}
\label{Forme Relative come Forme Globali}Let $X$ be a finite dimensional
complex analytic variety, $Z$ the closure of a non empty open set which is
also a polyhedron of $X$ and $\mathfrak{i}:Z\hookrightarrow X$ the
inclusion. Then $A^{p}:H_{ext}^{p}(X)_{Z}\rightarrow H_{ext}^{p}(X,Z)$ is an
isomorphism for\ any $p\in 
\mathbb{N}
$.
\end{proposition}

\begin{proof}[Proof (Sketch)]
Injectivity of $A^{p}$.\newline
Let $\omega \in \Gamma _{ext}(\Lambda ^{p}T^{\ast }X^{\prime })_{Z}$ be such
that $d^{p}(\omega )=0$ and $A^{p}([\omega
]_{H_{ext}^{p}(X)_{Z}})=0_{H_{ext}^{p}(X,Z)}$. This means that $\alpha
^{p}(\omega )=(\omega ,$ $0)$ is $d_{\mathfrak{i}}$-exact. Then there exists 
$(\vartheta ,$ $\varsigma )\in \Gamma _{ext}(\Lambda ^{p-1}\mathfrak{i})$
such that $(\omega ,$ $0)=d_{\mathfrak{i}}^{p-1}(\vartheta ,$ $\varsigma
)=(d^{p-1}\vartheta ,$ $\vartheta |_{Z}-d^{p-2}\varsigma )$. Thus, $\omega $
is an exact global form and $\vartheta |_{Z}-d^{p-2}\varsigma =0$. Actually,
this is not enough to conclude, because we want a primitive of $\omega $
vanishing on $Z$. Let $\tilde{\varsigma}$\ be\ an extension of $\varsigma $\
to the whole of $X$. Note that such an extension there exists, because the
sheaf $\mathcal{S}_{p-2}$\ is soft (cp. Corollary \ref{Fasci Fini II} of
Lemma \ref{Esiste Part. d'Unita' "Liscia" copy(1)}). Now, set $(\vartheta
^{\prime },$ $\varsigma ^{\prime })=(\vartheta -d^{p-2}\tilde{\varsigma},$ $%
0)\in \Gamma _{ext}(\Lambda ^{p-1}\mathfrak{i})$. Then $d_{\mathfrak{i}%
}^{p-1}(\vartheta ^{\prime },$ $\varsigma ^{\prime })=(\omega ,$ $0)$ and $%
\vartheta ^{\prime }|_{Z}=\vartheta |_{Z}-d^{p-2}\varsigma =0$. We are done,
because $\vartheta ^{\prime }$ is a primitive of $\omega $ vanishing on $Z$,
because.

Surjectivity of $A^{p}$.\newline
Let $(\varkappa ,$ $\upsilon )\in \Gamma _{ext}(\Lambda ^{p}\mathfrak{i})$
be $d_{\mathfrak{i}}^{p}$-closed. We want to find a $d^{p}$-closed $\omega
\in \Gamma _{ext}(\Lambda ^{p}T^{\ast }X^{\prime })_{Z}$ such that $%
[(\varkappa ,$ $\upsilon )]_{H_{ext}^{p}(X,Z)}=[(\omega ,$ $%
0)]_{H_{ext}^{p}(X,Z)}$. That is, we want to find a $d^{p}$-closed $\omega
\in \Gamma _{ext}(\Lambda ^{p}T^{\ast }X^{\prime })_{Z}$ and $(\vartheta ,$ $%
\varsigma )\in \Gamma _{ext}(\Lambda ^{p-1}\mathfrak{i})$ such that $%
(\varkappa ,$ $\upsilon )=(\omega ,$ $0)+d_{\mathfrak{i}}^{p-1}(\vartheta ,$ 
$\varsigma )=(\omega +d^{p-1}\vartheta ,$ $\vartheta |_{Z}-d^{p-2}\varsigma
) $. Since the sheaf $\mathcal{S}_{p-1}$\ is soft, there exists a global
extension $\tilde{\upsilon}$\ of $\upsilon $. Set $(\vartheta ,$ $\varsigma
)=(\tilde{\upsilon},$ $0)\in \Gamma _{ext}(\Lambda ^{p-1}\mathfrak{i})$ and
define $\omega =\varkappa -d^{p-1}\tilde{\upsilon}$. Then $(\varkappa ,$ $%
\upsilon )=(\omega ,$ $0)+d_{\mathfrak{i}}^{p-1}(\vartheta ,$ $\varsigma )$.
Moreover, $d^{p}\varkappa =0$ and $\varkappa |_{Z}-d^{p-1}\upsilon =0$,
because $(\varkappa ,$ $\upsilon )\in \Gamma _{ext}(\Lambda ^{p}\mathfrak{i}%
) $ is $d_{\mathfrak{i}}^{p}$-closed. Then, by the very definition of $%
\omega $, it results $\omega |_{Z}=0$. We are done, because $\omega $ is an
extendable differentiable form vanishing on $Z$ such that $[(\varkappa ,$ $%
\upsilon )]_{H_{ext}^{p}(X,Z)}=A^{p}([\omega ]_{H_{ext}^{p}(X)_{Z}})$.
\end{proof}

\subsection{Extendable \v{C}ech cohomology groups\label{Ext. Cech cohomology
groups (Sottosezione)}}

Let $X$ be a finite dimensional complex analytic variety and $Z$ the closure
of a non empty open set that is also a polyhedron of $X$. An open covering $%
\mathcal{V}=\{V_{j}\}_{j\in J}$ of $X$\ is \emph{adapted to }$Z$ if there
exists a unique $j(Z)\in J$ such that $Z\subset V_{j(Z)}$ and if for any $%
j\in J\setminus \{j(Z)\}$ it holds $Z\cap V_{j}=\emptyset $. Sometimes the
set $V_{j(Z)}$ is simply denoted by $V_{Z}$. Let us agree that every open
covering of $X$ is adapted to the empty set.

\begin{lemma}
\label{Raffinamento adattato a Z}Let $X$ and $Z$ be as above. If $\mathcal{V}
$ is an open covering of $X$ adapted to $Z$, then there exists an open
covering $\mathcal{B}$\ of $X$ adapted to $Z$ which is a locally finite
refinement of $\mathcal{V}$.
\end{lemma}

\begin{proof}
Let $\mathcal{B}^{\ast }=\{B_{\beta ^{\ast }}^{\ast }\}_{\beta ^{\ast }\in 
\mathfrak{B}^{\ast }}$ be a locally finite refinement of $\mathcal{V}$ and $%
\lambda :\mathfrak{B}^{\ast }\rightarrow J$ a refinement map associated with 
$\mathcal{B}^{\ast }$\ and $\mathcal{V}$. Write $\mathfrak{B}^{\ast }$, the
set of indices $\mathcal{B}^{\ast }$, as the disjoint union of $\{\beta
^{\ast }\in \mathfrak{B}^{\ast }:Z\cap B_{\beta ^{\ast }}^{\ast }\neq
\emptyset \}$ and $\{\beta ^{\ast }\in \mathfrak{B}^{\ast }:Z\cap B_{\beta
^{\ast }}^{\ast }=\emptyset \}$ and note that, if $B_{\bar{\beta}^{\ast
}}^{\ast }\in \mathcal{B}^{\ast }$ is such that $Z\cap B_{\bar{\beta}^{\ast
}}^{\ast }\neq \emptyset $, then $\lambda (\bar{\beta}^{\ast })=j(Z)$,
because $\mathcal{V}$\ is adapted to $Z$. Now, set $B_{Z}=\tbigcup%
\nolimits_{\beta ^{\ast }\in \mathfrak{B}^{\ast }:Z\cap B_{\beta ^{\ast
}}^{\ast }\neq \emptyset }B_{\beta ^{\ast }}^{\ast }$ and for any $\beta
^{\ast }\in \mathfrak{B}^{\ast }$ such that $Z\cap B_{\beta ^{\ast }}^{\ast
}=\emptyset $ set $\beta =\beta ^{\ast }\ $and $B_{\beta }=B_{\beta ^{\ast
}}^{\ast }$. Write $\mathfrak{B}=\{Z\}\cup (\mathfrak{B}^{\ast }\setminus
\{\beta ^{\ast }\in \mathfrak{B}^{\ast }:Z\cap B_{\beta ^{\ast }}^{\ast
}\neq \emptyset \})$. Then $\mathcal{B}=\{B_{\beta }\}_{\beta \in \mathfrak{B%
}}$ is the wanted covering.
\end{proof}

Let $X$ be a finite dimensional complex analytic variety, $Z$ either the
empty set or the closure of a non empty open set that is a polyhedron of $X$
and $\mathcal{V}=\left\{ V_{j}\right\} _{j\in J}$ an open covering of $X$
adapted to $Z$.

As a matter of notations, for any $q\in 
\mathbb{N}
$ and $\boldsymbol{j}=(j_{0},...,j_{q})\in J^{q+1}$ set $V_{\boldsymbol{j}%
}=V_{j_{0}}\cap ...\cap V_{j_{q}}$ and $V_{\boldsymbol{j}}^{\prime }=V_{%
\boldsymbol{j}}\setminus Sing(X)$. Moreover, for any $q\in 
\mathbb{N}
$, $m\in \left\{ 0,..,q\right\} $ and $\boldsymbol{j}=(j_{0},...,j_{q})\in
J^{q+1}$ set $\boldsymbol{j}_{m}=(j_{0},...,\widehat{j_{m}},...,j_{q})\in
J^{q}$ and denote by $\partial _{q-1,m}^{\boldsymbol{j}}:V_{\boldsymbol{j}%
}\hookrightarrow V_{\boldsymbol{j}_{m}}$ the inclusion. Note that, if $q\in 
\mathbb{N}
\setminus \{0\}$ and $\boldsymbol{j}\in J^{q+1}$, then $V_{\boldsymbol{j}%
}\cap Z=\emptyset $.

Let $p\in 
\mathbb{N}
$. The set%
\begin{equation}
C_{q}^{p}(X,\mathcal{V},V_{Z})=\tprod\nolimits_{\boldsymbol{j}\in
J^{q+1}}\Gamma _{ext}(\Lambda ^{p}T^{\ast }V_{\boldsymbol{j}}^{\prime })_{Z}
\label{Forme Est. associate al ricoprimento V}
\end{equation}%
is \emph{the space of extendable }$p$\emph{-forms associated with }$\mathcal{%
V}$\emph{\ and vanishing on }$Z$. Note that, since $\mathcal{V}$\ is adapted
to $Z$, the unique space really containing forms vanishing on $Z$ is $\Gamma
_{ext}(\Lambda ^{p}T^{\ast }V_{j(Z)}^{\prime })_{Z}\subset C_{0}^{p}(X,%
\mathcal{V},V_{Z})$.

The elements of $\left\{ \partial _{q-1,m}\right\} _{q\in 
\mathbb{N}
\setminus \{0\},\text{ }m\in \left\{ 0,...,q\right\} }$ give rise to a
sequence of families of maps%
\begin{equation*}
\left\{ \delta _{q-1,m}^{p}:\tprod\nolimits_{\boldsymbol{k}\in J^{q}}\Gamma
_{ext}(\Lambda ^{p}T^{\ast }V_{\boldsymbol{k}}^{\prime })_{Z}\text{ }%
\rightarrow \text{ }\tprod\nolimits_{\boldsymbol{j}\in J^{q+1}}\Gamma
_{ext}(\Lambda ^{p}T^{\ast }V_{\boldsymbol{j}}^{\prime })_{Z}\right\} _{q\in 
\mathbb{N}
,\text{ }m\in \left\{ 0,...,q\right\} ,\text{ }p\in 
\mathbb{N}
}
\end{equation*}%
depending on $p\in 
\mathbb{N}
$. For each $p\in 
\mathbb{N}
$ the map $\delta _{q-1,m}^{p}$\ is described as follows. The image $%
((\delta _{q-1,m}^{p}(\omega _{\boldsymbol{k}})_{\boldsymbol{k}\in J^{q}})_{%
\boldsymbol{j}})_{\boldsymbol{j}\in J^{q+1}}$ of $(\omega _{\boldsymbol{k}%
})_{\boldsymbol{k}\in J^{q}}\in \tprod\nolimits_{\boldsymbol{k}\in
J^{q}}\Gamma _{ext}(\Lambda ^{p}T^{\ast }V_{\boldsymbol{k}}^{\prime })_{Z}$
is the element of $\tprod\nolimits_{\boldsymbol{j}\in J^{q+1}}\Gamma
_{ext}(\Lambda ^{p}T^{\ast }V_{\boldsymbol{j}}^{\prime })_{Z}$ whose $%
\boldsymbol{j}^{th}$ term $(\delta _{q-1,m}^{p}(\omega _{\boldsymbol{k}})_{%
\boldsymbol{k}\in J^{q}})_{\boldsymbol{j}}$ is given by $(\partial
_{q-1,m}|_{V_{\boldsymbol{j}}^{\prime }})^{\ast }(\omega _{\boldsymbol{j}%
_{m}})$. Sometimes we write $\omega _{\boldsymbol{j}_{m}}$ instead of $%
((\partial _{q-1,m}|_{V_{\boldsymbol{j}}^{\prime }})^{\ast }(\omega _{%
\boldsymbol{k}})_{\boldsymbol{k}\in J^{q}})_{\boldsymbol{j}}$. So, $(\delta
_{q-1,m}^{p}(\omega _{\boldsymbol{k}})_{\boldsymbol{k}\in J^{q}})_{%
\boldsymbol{j}}=\omega _{\boldsymbol{j}_{m}}$.

For any $j_{0}\in J$ denote by $\partial ^{j_{0}}:V_{j_{0}}\hookrightarrow X$
the inclusion and define $\varrho :\tcoprod\nolimits_{j_{0}\in
J}V_{j_{0}}\rightarrow X$ by setting $\varrho |_{V_{j_{0}}}=\partial
^{j_{0}} $. Then $\varrho $\ gives rise to a sequence of maps $%
\{P^{p}:\Gamma _{ext}(\Lambda ^{p}T^{\ast }X^{\prime })_{Z}\rightarrow
\tprod\nolimits_{j_{0}\in J}\Gamma _{ext}(\Lambda ^{p}T^{\ast
}V_{j_{0}}^{\prime })_{Z}\}_{p\in 
\mathbb{N}
}$ depending on $p\in 
\mathbb{N}
$. For each $p\in 
\mathbb{N}
$ the map $P^{p}$ is described as follows. The image of $\omega \in \Gamma
_{ext}(\Lambda ^{p}T^{\ast }X^{\prime })$ is the element $((P^{p}(\omega
)_{j_{0}})_{j_{0}\in J}$ of $\tprod\nolimits_{j_{0}\in J}\Gamma
_{ext}(\Lambda ^{p}T^{\ast }V_{j_{0}}^{\prime })_{Z}$ whose $j_{0}^{th}$
term is $(\varrho |_{V_{j_{0}}^{\prime }})^{\ast }\left( \omega \right) $.

For each $p\in 
\mathbb{N}
$ and $q\in 
\mathbb{N}
\setminus \{0\}$ the map%
\begin{equation}
\begin{array}{rccl}
\delta _{q-1}^{p}: & C_{q}^{p}(X,\mathcal{V},V_{Z}) & \rightarrow & 
C_{q+1}^{p}(X,\mathcal{V},V_{Z}) \\ 
& (\omega _{\boldsymbol{k}})_{\boldsymbol{k}\in J^{q}} & \mapsto & ((\delta
_{q-1}^{p}(\omega _{\boldsymbol{k}})_{\boldsymbol{k}\in J^{q}})_{\boldsymbol{%
j}})_{\boldsymbol{j}\in J^{q+1}}%
\end{array}
\label{Operatore Differenza (I)}
\end{equation}%
defined by%
\begin{equation}
(\delta _{q-1}^{p}(\omega _{\boldsymbol{k}})_{\boldsymbol{k}\in J^{q}})_{%
\boldsymbol{j}}=\tsum\nolimits_{m=0}^{q}\left( -1\right) ^{m}(\delta
_{q-1,m}^{p}(\omega _{\boldsymbol{k}})_{\boldsymbol{k}\in J^{q}})_{%
\boldsymbol{j}}  \label{Operatore Differenza (II)}
\end{equation}%
is called\emph{\ }$p$\emph{-difference operator associated with }$\mathcal{V}
$. A straightforward computation shows that for any $p\in 
\mathbb{N}
$ and $q\in 
\mathbb{N}
\setminus \{0\}$ it results $\delta _{q}^{p}\circ \delta _{q-1}^{p}=0$ and $%
\delta _{0}^{p}\circ P^{p}=0$.

\begin{proposition}
\label{Delta - Aciclicita' per la Coom.}Let $X$ be a finite dimensional
complex analytic variety, $Z$ either the empty set or the closure of a non
empty open set that is a polyhedron of $X$ and $\mathcal{V}=\left\{
V_{j}\right\} _{j\in J}$ an open covering of $X$ adapted to $Z$. Fix $p\in 
\mathbb{N}
$. Then the $\delta ^{p}$-cohomology of $0\rightarrow \Gamma _{ext}(\Lambda
^{p}T^{\ast }X^{\prime })_{Z}\rightarrow C_{0}^{p}(X,\mathcal{V}%
,V_{Z})\rightarrow ...$ is identically zero.
\end{proposition}

\begin{proof}
By Lemma \ref{Raffinamento adattato a Z}, there exists\ an open covering $%
\mathcal{B}$ of $X$ adapted to $Z$ which is a locally finite refinement of $%
\mathcal{V}$. Moreover, by Lemma \ref{Lemma su Topologia e Ricoprimenti}, we
can assumewithout loss of generality that $\mathcal{B}$ enjoys (2) of Lemma %
\ref{Lemma su Topologia e Ricoprimenti}. It follows from Lemma \ref{Esiste
Part. d'Unita' "Liscia" copy(1)} that there exists an extendable partition
of unity subordinated to $\mathcal{B}$. Finally, proceed as in \cite{Gunning}%
, Vol. III, Ch. E.
\end{proof}

Let $X$ be a finite dimensional complex analytic variety, $Z$ either the
empty set or the closure of a non empty open set that is a polyhedron of $X$
and $\mathcal{V}=\left\{ V_{j}\right\} _{j\in J}$ an open covering of $X$
adapted to $Z$. For each $r\in 
\mathbb{N}
$ set $K^{r}(X,\mathcal{V},V_{Z})=\oplus _{p,q\in 
\mathbb{N}
:p+q=r}C_{q}^{p}(X,\mathcal{V},V_{Z})$. The space $K(X,\mathcal{V}%
,V_{Z})=\oplus _{r\in 
\mathbb{N}
}K_{ext}^{r}(X,\mathcal{V},V_{Z})$ is \emph{the space of extendable form
associated with }$X$\emph{, }$Z$\emph{\ and }$\mathcal{V}$.

For any $r\in 
\mathbb{N}
$ let $D^{r}:K^{r}(X,\mathcal{V},V_{Z})\rightarrow K^{r+1}(X,\mathcal{V}%
,V_{Z})$ be the operator that on each $C_{q}^{p}(X,\mathcal{V},V_{Z},)$ with 
$p+q=r$ is defined by%
\begin{equation}
D^{r}|_{C_{q}^{p}(X,\mathcal{V},V_{Z})}=\delta _{q}^{p}+\left( -1\right)
^{q}d^{p}  \label{Differenziale di Cech}
\end{equation}%
and that is identically zero otherwise. It easy to prove that $D^{r+1}\circ
D^{r}=0$ for any $r\in 
\mathbb{N}
$. Denote by $D:K(X,\mathcal{V},V_{Z})\rightarrow K(X,\mathcal{V},V_{Z})$
the operator that on each $K^{r}(X,\mathcal{V},V_{Z})$ is given by $%
D|_{K^{r}(X,\mathcal{V},V_{Z})}=D^{r}$.

\begin{definition}
Let $X$ be a finite dimensional complex analytic variety, $Z$ either the
empty set or the closure of a non empty open set which is a polyhedron of $X$
and $\mathcal{V}$ an open covering of $X$ adapted to $Z$. Fix $r\in 
\mathbb{N}
$. Set $\check{Z}_{ext}^{r}(X,\mathcal{V},V_{Z})=\ker (D^{r})$ and $\check{B}%
_{ext}^{r}(X,\mathcal{V},V_{Z})=\func{Im}(D^{r-1})$. The group%
\begin{equation*}
\check{H}_{ext}^{r}(X,\mathcal{V},V_{Z})=\frac{\check{Z}_{ext}^{r}(X,%
\mathcal{V},V_{Z})}{\check{B}_{ext}^{r}(X,\mathcal{V},V_{Z})}
\end{equation*}%
is \emph{the }$r^{th}$\emph{\ extendable \v{C}ech cohomology group
associated with }$X$\emph{, }$Z$\emph{\ and }$\mathcal{V}$.
\end{definition}

We have the following theorem.

\begin{theorem}
\label{Isomorfismo tra Coomologie}Let $X$ be a complex analytic variety of
finite complex dimension, $Z$ either the empty set or the closure of a non
empty open set which is a polyhedron of $X$ and $\mathcal{V}$ an open
covering of $X$ adapted to $Z$. Fix $r\in 
\mathbb{N}
$. Then $P^{r}:\Gamma _{ext}\left( \Lambda ^{r}T^{\ast }X^{\prime }\right)
_{Z}\rightarrow K^{r}(X,\mathcal{V},V_{Z})$ induces an isomorphism $P^{r\ast
}:H_{ext}^{r}\left( X,Z\right) \rightarrow \check{H}_{ext}^{r}(X,\mathcal{V}%
,V_{Z})$.
\end{theorem}

\begin{proof}
Proceed as for the proof of Proposition II.8.8 of \cite{Bott-Tu}.
\end{proof}

We have the following remark.

\begin{remark}
\label{Omotopia di Catene (Relativa)}Let $X$ be a complex analytic variety
of finite complex dimension, $Z$ either the empty set or the closure of a
non empty open set which is a polyhedron of $X$ and $\mathcal{V}$ an open
covering of $X$ adapted to $Z$. Then

\begin{enumerate}
\item The possibility of adapting the proofs of Proposition \ref{Delta -
Aciclicita' per la Coom.}\ and Theorem \ref{Isomorfismo tra Coomologie}\ to
the case of extendable forms is a direct consequence of Lemma \ref{Esiste
Part. d'Unita' "Liscia" copy(1)}.

\item The complexes of cochains $\Gamma _{ext}(\oplus _{r\in 
\mathbb{N}
}\Lambda ^{r}T^{\ast }X^{\prime })_{Z}$\ and $K(X,\mathcal{V},V_{Z})$\ are
chain homotopic. More precisely, let $\mathcal{B}=\{B_{\beta }\}_{\beta \in 
\mathfrak{B}}$\ be an open covering of $X$ refining $\mathcal{V}$ and
enjoying (2) of Lemma \ref{Lemma su Topologia e Ricoprimenti}. Let $\{\rho
\}=\{\rho _{\beta }:X\rightarrow 
\mathbb{R}
\}_{\beta \in \mathfrak{B}}$ be an extendable partition of unity
subordinated to $\mathcal{B}$. Then there exist a homotopy operator $L:K(X,%
\mathcal{V},V_{Z})\rightarrow K(X,\mathcal{V},V_{Z})$ and a chain map $%
\mathbf{\phi }=\mathbf{\phi }_{X,Z,\mathcal{V},\mathcal{B},\{\rho \}}:K(X,%
\mathcal{V},V_{Z})\rightarrow \Gamma _{ext}(\oplus _{r\in 
\mathbb{N}
}\Lambda ^{r}T^{\ast }X^{\prime })_{Z}$ such that $\mathbf{\phi }\circ
P=id_{\Gamma _{ext}(\oplus _{r\in 
\mathbb{N}
}\Lambda ^{r}T^{\ast }X^{\prime })_{Z}}$ and $P\circ \mathbf{\phi }=D\circ
L-L\circ D+id_{K(X,\mathcal{V},V_{Z})}$. For a proof, see \cite{Bott-Tu},
Ch. II, Sec. 9.

Note that, if $\mathcal{V}=\left\{ V_{0},V_{1}\right\} $\ only contains two
open sets, then for each $r\in 
\mathbb{N}
$ it results%
\begin{equation*}
K_{ext}^{r}(X,\mathcal{V})=\Gamma _{ext}(\Lambda ^{r}T^{\ast }V_{0}^{\prime
})\oplus \Gamma _{ext}(\Lambda ^{r}T^{\ast }V_{1}^{\prime })\oplus \Gamma
_{ext}(\Lambda ^{r-1}T^{\ast }V_{(0,1)}^{\prime }).
\end{equation*}%
Moreover, if $\{\rho _{0},\rho _{1}:X\rightarrow 
\mathbb{R}
\}$ is an extendable partition of unity subordinated to $\mathcal{B}=%
\mathcal{V}$, then $\mathbf{\phi }$ is the map that to $\check{\omega}=(%
\check{\omega}_{0},$ $\check{\omega}_{1},$ $\check{\omega}_{01})\in
K_{ext}^{r}(X,\mathcal{V})$ associates the $r$-form%
\begin{equation}
\omega =\rho _{0}\check{\omega}_{0}+\rho _{1}\check{\omega}_{1}-d\rho
_{0}\wedge \check{\omega}_{01}  \label{Forma reincollata}
\end{equation}
\end{enumerate}
\end{remark}

As a note, the theory that in this section has been developed for the real
bundle $TX^{\prime }$ holds also for the vector bundles $\mathbf{T}X^{\prime
}$, $\mathbf{\bar{T}}X^{\prime }$ and $T^{%
\mathbb{C}
}X^{\prime }$.

\section{Integration\label{Integration (Sezione)}}

\subsection{Integration of extendable forms}

For the proof of the following theorem, see \cite{Hardt} and \cite%
{Lojasiewicz}.

\begin{theorem}
\label{Teorema di Triangolabilita'}Let $X$ be an abstract finite dimensional
complex analytic variety. Then $X$ is a triangulable topological space.
\end{theorem}

As a matter of notations, let $X$ be an abstract finite dimensional complex
analytic variety. For each triangulation $\mathbb{T}$ of $X$ the simplicial
complex associated with $X$ and $\mathbb{T}$ is denoted by $(X,\mathbb{T})$.
Let $h\in \{0,...,\dim _{%
\mathbb{R}
}(X)\}$. The set of the $h$-simplices of $\mathbb{T}$ and the $h$-skeleton
of $(X,\mathbb{T})$\ are denoted by $\mathbb{T}_{h}$ and, respectively, $%
Skel^{h}(X,\mathbb{T})$. The groups of $h$-chains\ and $h$-cochains with
coefficients in $%
\mathbb{C}
$ associated with $\mathbb{T}$ are denoted by $C_{h}^{\mathbb{T}}\left(
X\right) $ and, respectively, $C_{\mathbb{T}}^{h}\left( X\right) $. The
groups $C_{h}^{\mathbb{T}}\left( X\right) $\ and $C_{\mathbb{T}}^{h}\left(
X\right) $\ are endowed with a structure of complex vector spaces. Moreover, 
$C_{\mathbb{T}}^{h}\left( X\right) $ is the dual space of $C_{h}^{\mathbb{T}%
}\left( X\right) $. The $h^{th}$ simplicial homology and cohomology spaces
of $X$ with coefficients in $%
\mathbb{C}
$ are denoted by $H_{h}\left( X\right) $ and, respectively, by $H^{h}\left(
X\right) $. Recall that $H^{h}(X)$ is the dual space of $H_{h}(X)$.

Let $\mathcal{C}=\{C_{l}\}_{l\in L}$ be an open covering of $X$. A
triangulation $\mathbb{T}$ of $X$ is $\mathcal{C}$\emph{-small} if for any
simplex $\Delta \in \mathbb{T}$ there is $C_{l(\Delta )}\in \mathcal{C}$
such that $\Delta \subseteq C_{l(\Delta )}$.

\begin{lemma}
\label{Triangolazioni e Ricoprimenti}Let $X$ be a triangulable topological
space, $\mathcal{C}$ an open covering of $X$ and $\mathbb{T}$ a
triangulation of $X$. Then there exists a natural number $b(\mathcal{C)}\in 
\mathbb{N}
$ such that $\mathbb{T}^{b(\mathcal{C)}}$, the $b(\mathcal{C)}^{th}$
barycentric subdivision of $\mathbb{T}$, is $\mathcal{C}$-small.
\end{lemma}

\begin{proof}
See \cite{Munkres}, Ch. 4, Sec. 31.
\end{proof}

Recall that, if $\Delta ^{h}$ denote the standard simplex of real dimension $%
h$, then%
\begin{equation}
H_{k}(\Delta ^{h})=H_{k}(\overline{\Delta ^{h}})=\left\{ 
\begin{array}{lll}
\mathbb{C}
& \text{ } & for\text{ }k\in 
\mathbb{N}
:k=0 \\ 
0 & \text{ } & for\text{ }k\in 
\mathbb{N}
:k>0%
\end{array}%
\right.  \label{Omologia dei Simplessi}
\end{equation}

Let $X$ be a finite dimensional complex analytic variety. For any $\omega \ $%
in $\Gamma _{ext}(\Lambda ^{p}T^{\ast }X^{\prime })$ there is an atlas $%
\mathcal{A}\left( \omega \right) =\{(A_{l},n_{l},U_{l},W_{l},F_{l})\}_{l\in
L}$ of trivializing extension for $\Lambda ^{p}T^{\ast }X^{\prime }$ such
that $\omega $ is completely extendable on $A_{l}$ for any $l\in L$ (cp. (%
\ref{Notazione Biolomorfismi 1}), (\ref{Notazione Atlante}), Remark \ref%
{Atlante Associato}, Terminology \ref{Terminologia}). An atlas as $\mathcal{A%
}\left( \omega \right) $ is an\emph{\ atlas of extensibility of }$\omega $.
Suppose that an atlas of extensibility $\mathcal{A}\left( \omega \right) $
is given for any $\omega \in \Gamma _{ext}(\Lambda ^{p}T^{\ast }X^{\prime })$%
.

As a matter of terminology, let $X_{\blacklozenge }\subseteq X$ be a complex
analytic subvariety of $X$. A triangulation $\mathbb{T}$ of $X$ is \emph{%
compatible with }$X_{\blacklozenge }$\emph{\ }if $X_{\blacklozenge }$ is
homeomorphic to a simplicial subcomplex of $(X,\mathbb{T})$ via the
homeomorphism between $X$ and $(X,\mathbb{T})$.

\begin{proposition}
\label{Integrale per Triangolazione(Omega)}Let $X$ be an abstract complex
analytic variety of complex dimension $n$. Let $\omega \in \Gamma
_{ext}(\Lambda ^{p}T^{\ast }X^{\prime })$ be an extendable $p$-form and $%
\mathcal{A}\left( \omega \right) $\ an atlas of extensibility of $\omega $.
Let $\mathbb{T}$ be an $\mathcal{A}\left( \omega \right) $-small
triangulation of $X$ compatible with the analytic subvariety $Sing\left(
X\right) $ of $X$. Then the map%
\begin{equation}
\begin{array}{rccl}
\tint \omega : & \mathbb{T}_{p} & \mathbb{\rightarrow } & \mathbb{%
\mathbb{C}
} \\ 
& \Delta & \mapsto & \tint\nolimits_{\Delta }\omega
=\tint\nolimits_{F_{A}\left( \Delta \right) }\tilde{\omega},%
\end{array}
\label{Integr. di Forme Estend.}
\end{equation}%
where $A\in \mathcal{A}\left( \omega \right) $ is such that $A\supseteq
\Delta $ and $\tilde{\omega}\in \Gamma (\Lambda ^{p}T^{\ast }U_{A})$ is any
extension of $\omega $, is well defined.
\end{proposition}

\begin{proof}
We have to show that $\tint \omega $ is well defined.

Let $\Delta \in \mathbb{T}_{p}$ be such that $\mathring{\Delta}\cap
Sing\left( X\right) =\emptyset $. Then the proof goes as in the non singular
case. Indeed, on one hand the topological boundary of $\partial \Delta $ is
Lebesgue trascurable. On the other hand, $\omega $ is bounded, because it is
extendable (see the proof of (1) at the beginning of Subsection \ref%
{Extendable cohomology groups (Sottosezione)}). Then it is possible to
define $\tint\nolimits_{\Delta }\omega =\tint\nolimits_{\mathring{\Delta}%
}\omega $ as in the non singular case. As a note, if $p=2n$, then $\mathring{%
\Delta}\cap Sing\left( X\right) =\emptyset $, because $\dim _{%
\mathbb{R}
}(Sing(X))\leq 2(n-1)$.

Next, suppose that $\Delta \in \mathbb{T}_{p}$ is such that $\mathring{\Delta%
}\cap Sing\left( X\right) \neq \emptyset $. We have to prove that the
definition of $\tint\nolimits_{\Delta }\omega $ does not depend neither on $%
A\in \mathcal{A}\left( \omega \right) $ such that $A\supseteq \Delta $ nor
on the extension $\tilde{\omega}\in \Gamma (\Lambda ^{p}T^{\ast }U_{A})$ of $%
\omega $ (cp. (\ref{Integr. di Forme Estend.})).

(1) Let $(A,n_{A},U_{A},W_{A},F_{A})\in \mathcal{A}\left( \omega \right) $
be such that $A\supseteq \Delta $ and write $A^{\prime }=A\setminus
Sing\left( X\right) $. If $\tilde{\omega}_{1},$ $\tilde{\omega}_{2}\in
\Gamma (\Lambda ^{p}T^{\ast }U_{A})$ are of two extensions $\omega $, then $%
\tilde{\omega}_{1}|_{F_{A}\left( A^{\prime }\right) }=\tilde{\omega}%
_{2}|_{F_{A}\left( A^{\prime }\right) }$, because $(F_{A}|_{A^{\prime
}})^{\ast }(\tilde{\omega}_{1})=\omega |_{A^{\prime }}=(F_{A}|_{A^{\prime
}})^{\ast }(\tilde{\omega}_{2})$. Furthermore, it results $\tilde{\omega}%
_{1}|_{F_{A}\left( A\right) }=\tilde{\omega}_{2}|_{F_{A}\left( A\right) }$,
because $\tilde{\omega}_{1}$ and $\tilde{\omega}_{2}$ are continuous and
coincide on the dense subset $F_{A}\left( A^{\prime }\right) $ of $%
F_{A}\left( A\right) $. Then $\tint\nolimits_{F_{A}\left( \Delta \right) }%
\tilde{\omega}_{1}=\tint\nolimits_{F_{A}\left( \Delta \right) }\tilde{\omega}%
_{2}$ and, in this case, the definition of $\tint\nolimits_{\Delta }\omega $%
\ is independent of the extension of $\omega $.

(2) Let $(A_{1},n_{1},U_{1},W_{1},F_{1}),$ $(A_{2},n_{2},U_{2},W_{2},F_{2})%
\in \mathcal{A}\left( \omega \right) $ be such that both $A_{1}$ and $A_{2}$
contain $\Delta $. Set $A_{(1,2)}=A_{1}\cap A_{2}$ and write $A_{1}^{\prime
}=A_{1}\setminus Sing\left( X\right) $, $A_{2}^{\prime }=A_{2}\setminus
Sing\left( X\right) $, $A_{(1,2)}^{\prime }=A_{(1,2)}\setminus Sing\left(
X\right) $ (see Section \ref{Notazioni (Sezione)}). If $\tilde{\omega}%
_{1}\in \Gamma (\Lambda ^{p}T^{\ast }U_{1})$ and $\tilde{\omega}_{2}\in
\Gamma (\Lambda ^{p}T^{\ast }U_{2})$ are of two extensions of $\omega $,
then $\omega |_{A_{1}^{\prime }}=(F_{1}|_{A_{1}^{\prime }})^{\ast }(\tilde{%
\omega}_{1})$ and $\omega |_{A_{2}^{\prime }}=(F_{2}|_{A_{2}^{\prime
}})^{\ast }(\tilde{\omega}_{2})$. So,%
\begin{eqnarray}
\tilde{\omega}_{1}|_{F_{1}(A_{(1,2)}^{\prime })} &=&(F_{\left( 2,1\right)
}|_{F_{2}(A_{(1,2)}^{\prime })})^{\ast }(\tilde{\omega}_{2})
\label{Integrale: Indipendenza dalla Rappresentazione} \\
\tilde{\omega}_{2}|_{F_{2}(A_{(1,2)}^{\prime })} &=&(F_{\left( 1,2\right)
}|_{F_{1}(A_{(1,2)}^{\prime })})^{\ast }(\tilde{\omega}_{1})  \notag
\end{eqnarray}%
and%
\begin{eqnarray}
d\tilde{\omega}_{1}|_{F_{1}(A_{(1,2)}^{\prime })} &=&(F_{\left( 2,1\right)
}|_{F_{2}(A_{(1,2)}^{\prime })})^{\ast }(d\tilde{\omega}_{2})
\label{Integrale: Indipendenza dalla Rappresentazione (3/2)} \\
d\tilde{\omega}_{2}|_{F_{2}(A_{(1,2)}^{\prime })} &=&(F_{\left( 1,2\right)
}|_{F_{1}(A_{(1,2)}^{\prime })})^{\ast }(d\tilde{\omega}_{1})  \notag
\end{eqnarray}%
We claim that, actually, it holds%
\begin{eqnarray}
\tilde{\omega}_{1}|_{F_{1}(A_{(1,2)})} &=&\tilde{\omega}_{2}\circ F_{\left(
2,1\right) }  \label{Integrale: Indipendenza dalla Rappresentazione (2)} \\
\tilde{\omega}_{2}|_{F_{2}(A_{(1,2)})} &=&\tilde{\omega}_{1}\circ F_{\left(
1,2\right) }  \notag
\end{eqnarray}%
In order to prove this, let $O_{1}$ be an open subset of $U_{1}$ such that $%
O_{1}\cap F_{1}(A_{1})=F_{1}(A_{(1,2)})$\ and $\mathbb{F}_{\left( 2,1\right)
}:O_{1}\rightarrow \mathbb{F}_{\left( 2,1\right) }\left( O_{1}\right) $ a
holomorphic extension of the biholomorphism $F_{\left( 2,1\right)
}:F_{1}(A_{(1,2)})\rightarrow F_{2}(A_{(1,2)})$ (cp. (\ref{Notazione
Biolomorfismi 2})). Let $O_{2}$ be an open subset of $U_{2}$ such that $%
O_{2}\cap F_{2}(A_{2})=F_{2}(A_{(1,2)})$. Shrinking $O_{1}$, if necessary,
we can assume without loss of generality that $\mathbb{F}_{\left( 2,1\right)
}(O_{1})\subseteq O_{2}$. Let $\iota :\mathbb{F}_{\left( 2,1\right)
}(O_{1})\hookrightarrow O_{2}$ denote the inclusion. Then%
\begin{eqnarray}
\tilde{\omega}_{1}|_{F_{1}(A_{(1,2)})} &=&(\mathbb{F}_{\left( 2,1\right)
}|_{F_{1}(A_{(1,2)})})^{\ast }\circ \iota ^{\ast }(\tilde{\omega}%
_{2}|_{O_{2}})  \label{omega(1) e' la tirata indietro di omega(2)} \\
&=&\tilde{\omega}_{2}|_{O_{2}}\circ \iota \circ \mathbb{F}_{\left(
2,1\right) }|_{F_{1}(A_{(1,2)})}  \notag
\end{eqnarray}%
Indeed, on one hand%
\begin{equation*}
F_{2}(A_{(1,2)})=F_{2}\circ F_{1}^{-1}\circ F_{1}(A_{(1,2)})\subseteq 
\mathbb{F}_{\left( 2,1\right) }\left( O_{1}\right) \subseteq O_{2}
\end{equation*}%
On the other hand, $\tilde{\omega}_{1}|_{F_{1}(A_{(1,2)})}$\ and $\tilde{%
\omega}_{2}|_{O_{2}}\circ \iota \circ \mathbb{F}_{\left( 2,1\right)
}|_{F_{1}(A_{(1,2)})}$ coincide on the dense subset $F_{1}(A_{(1,2)}^{\prime
})$\ of $F_{1}(A_{(1,2)})$, because of (\ref{Integrale: Indipendenza dalla
Rappresentazione}). Then, it results $\tilde{\omega}_{1}|_{F_{1}(A_{(1,2)})}=%
\tilde{\omega}_{2}|_{O_{2}}\circ \iota \circ \mathbb{F}_{\left( 2,1\right)
}|_{F_{1}(A_{(1,2)})}$. Then we get $\tilde{\omega}_{1}|_{F_{1}(A_{(1,2)})}=%
\tilde{\omega}_{2}\circ F_{2}\circ F_{1}^{-1}|_{F_{1}(A_{(1,2)})}$.
Analogously, $\tilde{\omega}_{2}|_{F_{2}(A_{(1,2)})}=\tilde{\omega}_{1}\circ
F_{1}\circ F_{2}^{-1}|_{F_{2}(A_{(1,2)})}$. So, (\ref{Integrale:
Indipendenza dalla Rappresentazione (2)}) are proved.

Furthermore, a similar argument shows that%
\begin{eqnarray}
d(\tilde{\omega}_{1}|_{F_{1}(A_{(1,2)})}) &=&d\tilde{\omega}_{2}\circ
F_{\left( 2,1\right) }
\label{Integrale: Indipendenza dalla Rappresentazione (3)} \\
d(\tilde{\omega}_{2}|_{F_{2}(A_{(1,2)})}) &=&d\tilde{\omega}_{1}\circ
F_{\left( 1,2\right) }  \notag
\end{eqnarray}

Now, since $\Delta \in \mathbb{T}_{p}$ is such that $\mathring{\Delta}\cap
Sing\left( X\right) \neq \emptyset $, it results $p=\dim _{%
\mathbb{R}
}\left( \Delta \right) \lneq 2n$, because the triangulation $\mathbb{T}$ is
compatible with $Sing(X)$. Then there exists a simplex $E\in \mathbb{T}$
such that $\dim _{%
\mathbb{R}
}\left( E\right) \gneq \dim _{%
\mathbb{R}
}\left( \Delta \right) $, $\partial \left( E\right) \supseteq \Delta $ and $%
\mathring{E}\cap Sing\left( X\right) =\emptyset $. Thus, by the geometry of
simplices and (\ref{Omologia dei Simplessi}), there are a $p$-chain $\Delta
^{\prime }\subset \mathring{E}\subset E$ and a $(p+1)$-chain $\Gamma \subset 
\mathring{E}\subset E$ such that $\Delta ^{\prime },$ $\Gamma \subset
A_{(1,2)}$ and $\Delta +\Delta ^{\prime }=\partial \left( \Gamma \right) $.
Note that $\Delta ^{\prime }\cap Sing\left( X\right) =\emptyset $ and $%
\mathring{\Gamma}\cap Sing\left( X\right) =\emptyset $, because $\Delta
^{\prime },$ $\Gamma \subset \mathring{E}$.

Next, since $F_{\left( 2,1\right) }|_{F_{1}(A_{12}^{\prime
})}:F_{1}(A_{12}^{\prime })\rightarrow F_{2}(A_{12}^{\prime })$ and $%
F_{\left( 1,2\right) }|_{F_{2}(A_{12}^{\prime })}:F_{2}(A_{12}^{\prime
})\rightarrow F_{1}(A_{12}^{\prime })$ are biholomorphisms between complex
manifolds, they preserve orientations. Then%
\begin{equation}
\tint\nolimits_{F_{1}(\Delta ^{\prime })}\tilde{\omega}_{1}=\tint%
\nolimits_{F_{(2,1)}\circ F_{1}(\Delta ^{\prime })}(F_{\left( 1,2\right)
}|_{F_{2}\left( A_{12}^{\prime }\right) })^{\ast }(\tilde{\omega}%
_{1})=\tint\nolimits_{F_{2}(\Delta ^{\prime })}\tilde{\omega}_{2}
\label{Integrale: Indipendenza dalla Rappresentazione passo 2}
\end{equation}%
and%
\begin{equation}
\tint\nolimits_{F_{1}(\mathring{\Gamma})}d\tilde{\omega}_{1}=\tint%
\nolimits_{F_{\left( 2,1\right) }\circ F_{1}(\mathring{\Gamma})}(F_{\left(
1,2\right) }|_{F_{2}\left( A_{12}^{\prime }\right) })^{\ast }(d\tilde{\omega}%
_{1})=\tint\nolimits_{F_{2}(\mathring{\Gamma})}d\tilde{\omega}_{2}
\label{Integrale: Indipendenza dalla Rappresentazione (4)}
\end{equation}

So, by (\ref{Integrale: Indipendenza dalla Rappresentazione (3)}) and (\ref%
{Integrale: Indipendenza dalla Rappresentazione (4)}), it results%
\begin{eqnarray}
\tint\nolimits_{F_{1}(\Delta )+F_{1}(\Delta ^{\prime })}\tilde{\omega}_{1}
&=&\tint\nolimits_{F_{1}(\Gamma )}d\tilde{\omega}_{1}=\tint\nolimits_{F_{1}(%
\Gamma )}d\tilde{\omega}_{2}\circ F_{\left( 2,1\right) }
\label{Integrale: Indipendenza dalla Rappresentazione (5)} \\
&=&\tint\nolimits_{F_{\left( 2,1\right) }\circ F_{1}(\Gamma )}d\tilde{\omega}%
_{2}=\tint\nolimits_{F_{2}(\Delta )+F_{2}(\Delta ^{\prime })}\tilde{\omega}%
_{2},  \notag
\end{eqnarray}%
because $F_{\left( 2,1\right) }\circ F_{1}(\Gamma )=F_{2}(\Gamma )$. Then,
by (\ref{Integrale: Indipendenza dalla Rappresentazione (2)}), it results $%
\tint\nolimits_{F_{1}(\Delta )}\tilde{\omega}_{1}=\tint\nolimits_{F_{2}(%
\Delta )}\tilde{\omega}_{2}$ and we are done.
\end{proof}

Let $X$ an abstract finite dimensional complex analytic variety and take $%
\omega \in \Gamma _{ext}(\Lambda ^{p}T^{\ast }X^{\prime })$. By Lemma \ref%
{Triangolazioni e Ricoprimenti}, the operator $\tint \omega :\mathbb{T}%
_{p}\rightarrow \mathbb{%
\mathbb{C}
}$ is well defined for any triangulation $\mathbb{T}$ of $X$ compatible with 
$Sing\left( X\right) $.

\begin{theorem}
\label{Teorema di Stokes}\textbf{(Stokes)} Let $X$ be an abstract finite
dimensional complex analytic variety and $\mathbb{T}$ a triangulation of $X$
compatible with $Sing\left( X\right) $. Let $\omega \in \Gamma
_{ext}(\Lambda ^{p}T^{\ast }X^{\prime })$, $C\in C_{p+1}^{\mathbb{T}}\left(
X\right) $ and $\iota :\partial C\hookrightarrow C$ be the inclusion. Then%
\begin{equation*}
\tint\nolimits_{C}d\omega =\tint\nolimits_{\partial C}\iota ^{\ast }\left(
\omega \right) .
\end{equation*}
\end{theorem}

\begin{proof}[Proof (Sketch)]
Let $\mathcal{A}(\omega )$ be an atlas of extensibility of $\omega $ and $%
\mathbb{T}^{b}$ an $\mathcal{A}(\omega )$-small barycentric subdivision of $%
\mathbb{T}$ (cp. Lemma \ref{Triangolazioni e Ricoprimenti}). Then $C\ $%
belongs to $C_{p+1}^{\mathbb{T}^{b}}\left( X\right) $. Namely, $%
C=\tsum\nolimits_{\lambda \in \Lambda }\Delta _{\lambda }$, with $\Delta
_{\lambda }\in \mathbb{T}_{p+1}^{b}$ for any $\lambda \in \Lambda $. Thus,
it suffices to show that the thesis holds when $C$ is a simplex $\Delta \in 
\mathbb{T}_{p+1}^{b}$.

If $\Delta \cap Sing\left( X\right) =\emptyset $, then there is nothing to
prove, being the classical result. If $\Delta \cap Sing\left( X\right) \neq
\emptyset $, then there exist $A\in \mathcal{A}(\omega )$ and $\tilde{\omega}%
\in \Gamma _{ext}(\Lambda ^{p}T^{\ast }U_{A})$ such that $A\supseteq \Delta $
and $\omega |_{A}=(F_{A}|_{A\setminus Sing\left( X\right) })^{\ast }(\tilde{%
\omega})$ (cp. Proposition \ref{Integrale per Triangolazione(Omega)}). Set $%
\Delta _{A}=F_{A}\left( \Delta \right) $. Then $\Delta _{A}\subseteq U_{A}$.
Moreover, $\partial \Delta _{A}=\partial F_{A}\left( \Delta \right)
=F_{A}\left( \partial \Delta \right) $, because $F_{A}$\ is a homeomorphism.
So, letting $\iota _{A}:\partial \Delta _{A}\hookrightarrow \Delta _{A}$
denote the inclusion, the wanted result follows from the classical Stokes'
theorem%
\begin{equation}
\tint\nolimits_{\Delta }d\omega =\tint\nolimits_{\Delta _{A}}d\tilde{\omega}%
=\tint\nolimits_{\partial \Delta _{A}}(\iota _{A})^{\ast }\left( \tilde{%
\omega}\right) =\tint\nolimits_{\partial \Delta }\iota ^{\ast }\left( \omega
\right) .  \label{Teorema di Stokes (dimostrazione)}
\end{equation}
\end{proof}

We have the following remark.

\begin{remark}
\label{Form. Integr. per Var. Comp. (Osservazione)}Let $X$ be a compact
irreducible complex analytic variety of complex dimension $n$ and $\mathbb{T}
$ a triangulation of $X$ compatible with $Sing\left( X\right) $. Take $%
\omega \in \Gamma _{ext}(\Lambda ^{2n}T^{\ast }X)$ and let $\mathcal{A}%
\left( \omega \right) $ be an atlas of extensibility of $\omega $. Let $b\in 
\mathbb{N}
$ be such that $\mathbb{T}^{b}$ is $\mathcal{A}\left( \omega \right) $%
-small. Then the number%
\begin{equation}
\tint\nolimits_{X}\omega =\tsum\nolimits_{\Delta \in \mathbb{T}%
^{b}}\tint\nolimits_{\Delta }\omega  \label{Form. Integr. per Var. Comp.}
\end{equation}%
is well defined. Indeed, even if in general two triangulations of $X$ do not
have a common refinement, actually, (\ref{Form. Integr. per Var. Comp.})
does not depend on the chosen triangulation $\mathbb{T}$ of $X$ compatible
with $Sing\left( X\right) $. In order to prove this, recall that, if $%
\mathbb{T}_{1}$and $\mathbb{T}_{2}$ are triangulations of $X$ compatible
with $Sing\left( X\right) $, then there exist triangulations $\mathbb{T}_{3}$%
, $\mathbb{T}_{I}$, $\mathbb{T}_{II}$ of $X$ compatible with $Sing\left(
X\right) $ such that $\mathbb{T}_{I}$ is a refinement of both $\mathbb{T}%
_{1} $ and $\mathbb{T}_{3}$ and $\mathbb{T}_{II}$ is a refinement of both $%
\mathbb{T}_{2}$\ and $\mathbb{T}_{3}$. The claimed independence follows from%
\begin{equation*}
\tsum\nolimits_{\Delta _{1}\in \mathbb{T}_{1}}\Delta
_{1}=\tsum\nolimits_{\Delta _{I}\in \mathbb{T}_{I}}\Delta
_{I}=\tsum\nolimits_{\Delta _{3}\in \mathbb{T}_{3}}\Delta
_{3}=\tsum\nolimits_{\Delta _{II}\in \mathbb{T}_{II}}\Delta
_{II}=\tsum\nolimits_{\Delta _{2}\in \mathbb{T}_{2}}\Delta _{2}.
\end{equation*}%
As a matter of terminology, $\tint\nolimits_{X}\omega $ is \emph{the
integral of }$\omega $\emph{\ on }$X$.
\end{remark}

\subsection{Integration of extendable cohomology classes}

Let $X$ be a finite dimensional complex analytic variety and $Z$ either the
empty set or the closure of an open set which is a polyhedron of $X$. For
any $p\in 
\mathbb{N}
$ the $p^{th}$ simplicial cohomology group relative to the pair $(X,Z)$ is
denoted by $H^{p}(X,Z)$ (cp. Theorem \ref{Teorema di Triangolabilita'}). Let 
$\omega \in \Gamma _{ext}(\Lambda ^{p}T^{\ast }X^{\prime })_{Z}$ be an
extendable $p$-form and $\mathbb{T}$\ a triangulation of $X$ compatible with 
$Sing\left( X\right) $ and $Z$. Extending by linearity $\tint \omega :%
\mathbb{T}_{p}\rightarrow \mathbb{%
\mathbb{C}
}$ to the whole of $C_{p}^{\mathbb{T}}\left( X\right) $, we get a well
defined map $\tint \omega :C_{p}^{\mathbb{T}}\left( X\right) \rightarrow 
\mathbb{%
\mathbb{C}
}$ that, in fact, lies in $C_{\mathbb{T}}^{p}\left( X\right) $.

So, we get a homomorphism between complex vector spaces%
\begin{equation}
\begin{array}{rccl}
\eta _{Z}^{p}: & \Gamma _{ext}(\Lambda ^{p}T^{\ast }X^{\prime })_{Z} & 
\rightarrow & C_{\mathbb{T}}^{p}\left( X\right) \\ 
& \omega & \mapsto & \eta _{Z}^{p}(\omega )=\tint \omega%
\end{array}
\label{Operatore di Integrazione}
\end{equation}%
$\eta _{Z}^{p}$ is called \emph{operator of integration of degree }$p$.

Furthermore, for any $\omega \in \Gamma _{ext}(\Lambda ^{p}T^{\ast
}X^{\prime })_{Z}$ the cochain $\eta _{Z}^{p}(\omega )$ lies in $C_{\mathbb{T%
}}^{p}\left( X,Z\right) $, because $\omega |_{Z}=0$. So,%
\begin{equation*}
\eta _{Z}^{p}:\Gamma _{ext}(\Lambda ^{p}T^{\ast }X^{\prime })_{Z}\rightarrow
C_{\mathbb{T}}^{p}\left( X,Z\right) .
\end{equation*}%
With slight abuses of notations, given any $\mathbf{\Delta }\in C_{p}^{%
\mathbb{T}}\left( X,Z\right) $, the cochain $\eta _{Z}^{p}\left( \omega
\right) $ acts as follows $\Delta +C_{p}^{\mathbb{T}}\left( Z\right) \ni 
\mathbf{\Delta }\mapsto \eta _{Z}^{p}\left( \omega \right) \left( \mathbf{%
\Delta }\right) =\tint\nolimits_{\Delta }\omega \in 
\mathbb{C}
$.

Let%
\begin{equation}
\begin{array}{rccl}
H_{Z}^{p}: & H_{ext}^{p}(X)_{Z} & \rightarrow & H^{p}(X,Z) \\ 
& [\omega ] & \mapsto & H_{Z}^{p}([\omega ])=\tint [\omega ]%
\end{array}
\label{Omomorfismo di Integrazione}
\end{equation}%
be the operator induced by $\eta _{Z}^{p}$. Then $H_{Z}^{p}([\omega ])\in
H^{p}(X,Z)$ is the map that associates to each $[\mathbf{C}]\in H_{p}(X,Z),$ 
$\mathbf{C}\in C+C_{p}^{\mathbb{T}}\left( Z\right) $ the number%
\begin{equation}
\tint\nolimits_{\lbrack \mathbf{C}]}[\omega ]=\tint\nolimits_{C}\omega .
\label{Integrale Omologico (Definizione)}
\end{equation}%
In order to prove that $H_{Z}^{p}$\ is well defined, we use Stokes' theorem.
Let $\mathbf{C}\in C+C_{p}^{\mathbb{T}}\left( Z\right) $ be a $p$-cycle of $%
(X,\mathbb{T})$ relative to $(X,Z)$ and $\omega \in \Gamma _{ext}(\Lambda
^{p}T^{\ast }X^{\prime })_{Z}$ a closed extendable $p$-form vanishing on $Z$%
. Since $\partial ^{p}C\subseteq Z\ $and $d^{p}\omega =0,$ $\omega |_{Z}=0$,
for each $E\in C_{p+1}^{\mathbb{T}}\left( X,Z\right) $\ and for each $\sigma
\in \Gamma _{ext}(\Lambda ^{p-1}T^{\ast }X)_{Z}$ we have%
\begin{eqnarray}
\tint\nolimits_{C+\partial ^{p+1}E+C_{p}^{T}\left( Z\right) }(\omega
+d^{p-1}\sigma ) &=&\tint\nolimits_{C}\omega +\tint\nolimits_{\partial
^{p}C}\sigma +\tint\nolimits_{C_{p}^{\mathbb{T}}\left( Z\right) }(\omega
+d^{p-1}\sigma )|_{Z}  \label{Omomorfismo di Integrazione (dimostrazione)} \\
&=&\tint\nolimits_{C}\omega +\tint\nolimits_{\partial ^{p}C}\sigma
=\tint\nolimits_{C}\omega ,  \notag
\end{eqnarray}%
because $(\omega +d^{p-1}\sigma )|_{Z}=0,$ $\partial ^{p}C\subseteq Z$ and $%
\sigma |_{Z}=0$. So, $H_{Z}^{p}:H_{ext}^{p}(X,Z)\rightarrow H^{p}(X,Z)$ is
well defined, because $\tint\nolimits_{[\mathbf{C}]}[\omega ]\ $does not
depend on the representatives chosen in $[C]_{H_{p}(X,Z)}$\ and $[\omega
]_{H_{ext}^{p}(X)_{Z}}$. Moreover, $H_{Z}^{p}$ is a homomorphism of vector
spaces, because such is $\eta _{Z}^{p}$. $H_{Z}^{p}\ $is called \emph{%
homomorphism of integration of degree }$p$\emph{\ relative to }$(X,Z)$.

\begin{lemma}
\label{Sottovarieta' Compatte}Let $X$ be a finite dimensional complex
analytic variety and $S$ a closed compact complex analytic subvariety of $X$%
. Such an $S$ is a polyhedron of $X$. Then there are arbitrarily small open
neighborhoods $U_{S}$ of $S$ in $X$ such that

\begin{enumerate}
\item $\overline{U_{S}}$ and $X\setminus U_{S}$ are polyhedra in $X$

\item The inclusions $S\rightarrow U_{S}\rightarrow \overline{U_{S}}$ and $%
(X\setminus \overline{U_{S}})\rightarrow (X\setminus U_{S})\rightarrow
(X\setminus S)$\ are homotopy equivalences.
\end{enumerate}
\end{lemma}

\begin{proof}
If $X$ is compact, then proceed as in \cite{Munkres}, Ch. 8, Sec. 72. If $X$
is not compact, we can agree as follows. Let $X^{\ast }$ be an open
neighborhood of $S$ in $X$ such that $\overline{X^{\ast }}$ is a compact
polyhedron in $X$ with respect to the same triangulation for which $S$ is a
polyhedron. Then the statement holds for $S$ considered as a polyhedron in $%
\overline{X^{\ast }}$. So, the following%
\begin{equation*}
\begin{array}[b]{ccccc}
S & \rightarrow & U_{S} & \rightarrow & \overline{U_{S}} \\ 
\overline{X^{\ast }}\setminus \overline{U_{S}} & \rightarrow & \overline{%
X^{\ast }}\setminus U_{S} & \rightarrow & \overline{X^{\ast }}\setminus S%
\end{array}%
\end{equation*}%
are homotopy equivalences, with $U_{S}$\ an open neighborhoods of $S$ in $%
\overline{X^{\ast }}$ and then in $X$. Furthermore, by their very
constructions (cp. \cite{Munkres}, page 429), the homotopy equivalences $%
\overline{X^{\ast }}\setminus \overline{U_{S}}\rightarrow \overline{X^{\ast }%
}\setminus U_{S}\rightarrow \overline{X^{\ast }}\setminus S$ do not involve
the boundary of $\overline{X^{\prime }}$\ in $X$. Then, in order to define
the wanted homotopy equivalences%
\begin{equation*}
X\setminus \overline{U_{S}}\rightarrow X\setminus U_{S}\rightarrow
X\setminus S,
\end{equation*}%
just glue the maps $\overline{X^{\ast }}\setminus \overline{U_{S}}%
\rightarrow \overline{X^{\ast }}\setminus U_{S}\rightarrow \overline{X^{\ast
}}\setminus S$ and $X\setminus \overline{X^{\ast }}\rightarrow X\setminus 
\overline{X^{\ast }}\rightarrow X\setminus \overline{X^{\ast }}$.
\end{proof}

\begin{remark}
Let $X$ be a finite dimensional complex analytic variety, $S$ a closed
compact complex analytic subvariety of $X$ and $U_{S}$\ an open neighborhood
of $S$ enjoying (1) and (2) of Lemma \ref{Sottovarieta' Compatte}. Set $%
Z=X\setminus U_{S}$ and for any $p\in 
\mathbb{N}
$ consider the homomorphism $H_{X\setminus U_{S}}^{p}:H_{ext}^{p}\left(
X,X\setminus U_{S}\right) \rightarrow H^{p}\left( X,X\setminus U_{S}\right) $%
. Then, since $X\setminus S$\ and $X\setminus U_{S}$\ are homotopically
equivalent, we get a homomorphism $H_{Z}^{p}:H_{ext}^{p}\left( X,Z\right)
\rightarrow H^{p}\left( X,X\setminus S\right) $, because the groups $%
H^{p}\left( X,X\setminus S\right) $ and $H^{p}\left( X,X\setminus
U_{S}\right) $ are isomorphic.
\end{remark}

An easy verification shows that the following properties of complex analytic
varieties hold. For the necessary background, refer to \cite{Munkres}.

\begin{remark}
\label{Topologia Algebrica (Osservazione)}Let $X$ be a complex analytic
variety of complex dimension $n$ and $\mathbb{T}$\ a triangulation of $X$
compatible with $Sing(X)$.

\begin{enumerate}
\item The pair $(X,Sing(X))$ is a relative homology $2n$-manifold (cp. \cite%
{Munkres}).

\item If $X$ is also irreducible, then $(X,Sing(X))$ is a relative pseudo $%
2n $-manifold, $X$ is the closure of the union of the $2n$-simplices of $%
\mathbb{T}$, any $(2n-1)$ simplex is a face of exactly two $2n$-simplices of 
$\mathbb{T}$ and the $2n$-simplices of $\mathbb{T}$\ can be coherently
oriented (cp. \cite{Munkres}).

\item If $X$ is a compact and irreducible, then $(X,Sing(X))$ is an
orientable relative pseudo $2n$-manifold, the groups $H_{2n}(X)$\ and $%
H_{2n}(X,Sing(X))$\ are isomorphic and the $2n$-dimensional simplices of $(X,%
\mathbb{T)}$ can be oriented in such a way that their sum is a non vanishing
cycle. Moreover, such a sum is indipendent of the chosen triangulation (cp.
Remark \ref{Form. Integr. per Var. Comp. (Osservazione)}). The class $[X]$
represented by such a cycle is \emph{the fundamental class of }$X$ (cp. \cite%
{Munkres}).
\end{enumerate}
\end{remark}

We have the following theorem.

\begin{theorem}
\label{Integrale Omologico per Var. Compatte}Let $X$ be a compact
irreducible complex analytic variety of complex dimension $n$. If $\omega
\in \Gamma _{ext}(\Lambda ^{2n}T^{\ast }X)$ be such that $d\omega =0$, then
the number%
\begin{equation}
\tint\nolimits_{\lbrack X]}[\omega ]=\tint\nolimits_{X}\omega
\label{Form. Integr. Omolog. per Var. Comp.}
\end{equation}%
is well defined. $\tint\nolimits_{[X]}[\omega ]$ is \emph{the integral of }$%
[\omega ]$\emph{\ on }$[X]$.
\end{theorem}

\begin{proof}[Proof (Sketch)]
The wanted result follows from (\ref{Integrale Omologico (Definizione)}).
Indeed, $(X,Sing\left( X\right) )$ is an orientable relative pseudomanifold
of real dimension $2n$ (cp. Remark \ref{Topologia Algebrica (Osservazione)}).
\end{proof}

As a matter of notations, if $Z=\emptyset $, we omit to write the subscript $%
Z$ everywhere it should appear. The next remark concerns the injectivity of
homomorphisms of integration.

\begin{remark}
Let $X$ be an abstract finite dimensional complex analytic variety and $p\in 
\mathbb{N}
$. Then $H^{p}:H_{ext}^{p}\left( X\right) \rightarrow H^{p}\left( X\right) $
is generally not injective.

Indeed, let $X$ be as in Example \ref{Es. Bloom Herrera Complesso}.
Shrinking the neighborhood $B$ of $0$ in $%
\mathbb{C}
$, if necessary, we can assume that $X$ is topologically contractible (cp. 
\cite{Goresky-MacPherson}). So, $H^{1}\left( X\right) =0$. Then $%
H^{1}:H_{ext}^{1}\left( X\right) \rightarrow H^{1}\left( X\right) $ is not
injective, because $H_{ext}^{1}\left( X\right) \neq 0$ (cp. Example \ref{Es.
Bloom Herrera Complesso}).
\end{remark}

The following proposition concerns the surjectivity of homomorphisms of
integration. For the necessary background in algebraic topology, see \cite%
{Munkres}.

\begin{proposition}
Let $X$ be an abstract finite dimensional complex analytic variety. Fix $%
p\in 
\mathbb{N}
$. Then $H^{p}:H_{ext}^{p}\left( X\right) \rightarrow H^{p}\left( X\right) $
is surjective.
\end{proposition}

\begin{proof}
Let $\mathbb{T}$ be a triangulation of $X$ and $c\in C_{\mathbb{T}%
}^{p}\left( X\right) $ be a $p$-cocycle. We look for a closed extendable $p$%
-form $\omega _{c}\in \Gamma _{ext}(\Lambda ^{p}T^{\ast }X^{\prime })$ such
that $[c]_{H^{p}(X)}=H^{p}([\omega _{c}]_{H_{ext}^{p}(X)})$. Let $\mathcal{V}%
(\mathbb{T)}$ be the open covering of $X$ by the open stars of vertices of $%
\mathbb{T}$. The nerve $N(\mathcal{V}(\mathbb{T)})$ of such a covering\ is
an abstract simplicial complex. Moreover, the vertex correspondence $C_{0}^{%
\mathbb{T}}\left( X\right) \ni \Delta ^{0}\mapsto St(\Delta ^{0})\in N(%
\mathcal{V}(\mathbb{T)})$ is an isomorphism between $(X,\mathbb{T})$\ and $N(%
\mathcal{V}(\mathbb{T)})$ (cp. \cite{Munkres}, Theorem 73.2). Then to any $p$%
-cocycle $c\in C_{\mathbb{T}}^{p}\left( X\right) $ corresponds a unique
cocycle $\mathbf{c}_{c}\in C^{p}(N(\mathcal{V}(\mathbb{T)}))$. Furthermore,
by the very definition of $C^{p}(N(\mathcal{V}(\mathbb{T)}))$, we have $%
C^{p}(N(\mathcal{V}(\mathbb{T)}))=C_{0}^{p}(X,\mathcal{V})$ (see (\ref{Forme
Est. associate al ricoprimento V})).

Let $\check{H}^{p}\left( X\right) $\ be the $p^{th}$ \v{C}ech cohomology
group of $X$ (see \cite{Munkres}, Ch. 8, Sec. 73). Let $\mathbb{T}$\ be so
fine that $\check{H}^{p}(X,N(\mathcal{V}(\mathbb{T)}))$ is isomorphic to $%
\check{H}^{p}\left( X\right) $ and let $\mathcal{B}$\ be an open covering of 
$X$ enjoying (2) of Lemma \ref{Lemma su Topologia e Ricoprimenti} and
refining $\mathcal{V}(\mathbb{T)}$. Let $\{\rho _{\beta }\}_{\beta \in 
\mathfrak{B}}$ be an extendable partition of unity subordinated to $\mathcal{%
B}$ and $\mathbf{\phi }:K(X,\mathcal{V}(\mathbb{T)})\rightarrow \Gamma
_{ext}(\oplus _{r\in 
\mathbb{N}
}\Lambda ^{r}T^{\ast }X^{\prime })$\ the chain map associated with $\mathcal{%
V}(\mathbb{T)}$, $\mathcal{B}$ and $\{\rho \}$ (see Remark (2)\ of \ref%
{Omotopia di Catene (Relativa)}). Then the closed extendable $p$-form $%
\omega _{c}\in \Gamma _{ext}(\Lambda ^{p}T^{\ast }X^{\prime })$ defined via
the collating formula $\omega _{c}=\mathbf{\phi }(\mathbf{c}_{c})$ of Bott
and Tu does the desired job (see \cite{Bott-Tu}, Ch. II, Sec. 9, Proposition
9.8).
\end{proof}

\subsection{Integration of extendable \v{C}ech cohomology classes\label%
{Integrazione delle Classi di Coom. di Cech Est. (Sottosezione)}}

Let $X$ be an abstract finite dimensional complex analytic variety, $Z$
either the empty set or the closure of an open set that is also a polyhedron
of $X$ and $\mathcal{V}$\ an open covering of $X$ adapted to $Z$. Let $r\in 
\mathbb{N}
$ and consider the group homomorphism $\check{H}_{Z}^{r}=H_{Z}^{r}\circ
(P^{r\ast })^{-1}:\check{H}_{ext}^{r}\left( X,\mathcal{V},V_{Z}\right)
\rightarrow H^{r}\left( X,Z\right) $ (cp. Theorem \ref{Isomorfismo tra
Coomologie}).

Let $\mathbb{T}$ be a triangulation of $X\ $compatible with $Sing\left(
X\right) $ and $Z$ and consider the operator $\eta _{Z}^{r}:\Gamma
_{ext}(\Lambda ^{r}T^{\ast }X)_{Z}\rightarrow C_{\mathbb{T}}^{r}(X,Z)$. Let $%
\mathcal{B}=\{B_{\beta }\}_{\beta \in \mathfrak{B}}$ be a refinement of $%
\mathcal{V}$ enjoying (2) of Lemma \ref{Lemma su Topologia e Ricoprimenti}, $%
\{\rho \}=\{\rho _{\beta }:X\rightarrow 
\mathbb{R}
\}_{\beta \in \mathfrak{B}}$ an extendable partition of unity subordinated
to $\mathcal{B}$ and $\mathbf{\phi }:K\left( X,\mathcal{V},V_{Z}\right)
\rightarrow \oplus _{r\in 
\mathbb{N}
}\Gamma _{ext}\left( \Lambda ^{r}T^{\ast }X\right) _{Z}$ the chain map
associated with $\mathcal{V}$, $\mathcal{B}$ and $\{\rho \}$ (cp. Remark \ref%
{Omotopia di Catene (Relativa)} (2)).

Let%
\begin{equation}
\check{\eta}_{Z}^{r}:\check{Z}_{ext}^{r}\left( X,\mathcal{V},V_{0}\right)
\rightarrow Z_{\mathbb{T}}^{r}\left( X,Z\right)
\label{Cech Integrale Relativo}
\end{equation}%
denote the homomorphism defined by%
\begin{equation}
\check{\eta}_{Z}^{r}=\eta _{Z}^{r}\circ \mathbf{\phi }
\label{Relaz. tra Oper. d'Integr. Relat. Piu' Esplicita}
\end{equation}%
Then $\check{H}_{Z}^{r}$ is induced by $\check{\eta}_{Z}^{r}$. Indeed, for
any $\check{\omega}\in \check{Z}_{ext}^{r}\left( X,\mathcal{V},V_{Z}\right) $
it results%
\begin{eqnarray}
\lbrack \check{\eta}_{Z}^{r}\left( \check{\omega}\right) ]_{H^{r}\left(
X,Z\right) } &=&[\eta _{Z}^{r}\circ \mathbf{\phi }(\check{\omega}%
)]_{H^{r}\left( X,Z\right) }=H_{Z}^{r}[\mathbf{\phi }\left( \check{\omega}%
\right) ]_{H_{ext}^{r}\left( X,Z\right) }
\label{Equiv. tra Oper. d'Integr. Rel.} \\
&=&H_{Z}^{r}\circ (P^{r\ast })^{-1}[\check{\omega}]_{\check{H}%
_{ext}^{r}\left( X,\mathcal{V},V_{Z}\right) }=\check{H}_{Z}^{r}[\check{\omega%
}]_{\check{H}_{ext}^{r}\left( X,\mathcal{V},V_{Z}\right) }  \notag
\end{eqnarray}

\begin{definition}
Let $X$ be a complex analytic variety of complex dimension $n$ and $\mathcal{%
V}=\left\{ V_{0},V_{1}\right\} $\ an open covering of $X$. A \emph{honeycomb
cell system associated with }$\mathcal{V}$\ is a family $\mathcal{R}=\left\{
R_{0},R_{1}\right\} $ of subsets of $X$ such that

\begin{enumerate}
\item For $j\in \{0,1\}$ $R_{j}$ is the closure of an $n$-dimensional
complex analytic subvariety of $X$ with piecewise differentiable boundary.

\item It holds $R_{0}\subsetneq V_{0}$ and $R_{1}\subsetneq V_{1}$

\item It holds $R_{0}\cup R_{1}=X$ and $\mathring{R}_{0}\cap \mathring{R}%
_{1}=\emptyset $

\item Both $\partial R_{0}$ and $\partial R_{1}$ are a complete intersection
and $R_{0}\cap R_{1}$ is a variety of real dimension $2n-1$ with piecewise
differentiable boundary.

\item $R_{(0,1)}$ denotes the hypersuface $R_{0}\cap R_{1}$ oriented in the
following way: the orientation on $R_{(0,1)}$ is that one defined by the
interior normal of $R_{0}$ and by the exterior normal of $R_{1}$.

\item $R_{(1,0)}$ denotes the hypersuface $R_{0}\cap R_{1}$ oriented in the
following way: the orientation on $R_{(1,0)}$ is that one defined by the
interior normal of $R_{1}$ and by the exterior normal of $R_{0}$.
\end{enumerate}
\end{definition}

We have the following remark.

\begin{remark}
Let $X$ be a complex analytic variety of finite dimension and $\mathcal{V}%
=\left\{ V_{0},V_{1}\right\} $ an open covering of $X$. Then there exists a
honeycomb cell system $\mathcal{R}=\left\{ R_{0},R_{1}\right\} $ associated
with $\mathcal{V}$. Indeed, every analytic variety $X$ admits a
triangulation\ such that each simplex of it corresponds to an analytic
subvariety of $X$ (see \cite{Lojasiewicz} and \cite{Hardt}). So, the
existence of a honeycomb cell system $\mathcal{R}$ associated to $\mathcal{V}
$ can be deduced from Lemma \ref{Triangolazioni e Ricoprimenti}.
\end{remark}

Next, let $X$ be an abstract finite dimensional complex analytic variety and 
$Z$ either the empty set or the closure of an open set that is also a
polyhedron of $X$. Let $\mathbb{T}$ be a triangulation of $X\ $compatible
with $Sing\left( X\right) $ and $Z$, take $r\in 
\mathbb{N}
$ and consider the $(r-1)$-skeleton $Skel^{r-1}(X,\mathbb{T)}$ of $\left( X,%
\mathbb{T}\right) $. In what follows we denote by $\sim $\ the homotopic
equivalence. Let $Z_{r-1}$ be the closure of an open set of $X$ that is a
polyhedron of $X$ such that $Z_{r-1}\supseteq Skel^{r-1}\left( X,\mathbb{T}%
\right) $ and $Z_{r-1}\sim Skel^{r-1}\left( X,\mathbb{T}\right) $. Set $%
Z_{\bullet }=Z\cup Z_{r-1}$.

Let $V_{0}$ be an open neighborhood of $Z_{\bullet }$ such that $V_{0}\sim
Z_{\bullet }$ and $V_{1}$ an open set of $X$ such that $V_{1}\cap Z_{\bullet
}=\emptyset $ and $V_{0}\cup V_{1}=X$. Set $\mathcal{V}=\left\{
V_{0},V_{1}\right\} $ and note that $\mathcal{V}$ is adapted to $Z_{\bullet
} $ . By Proposition \ref{Esiste Part. d'Unita' "Liscia" copy(1)}, there
exists an extendable partition of unity $\{\rho \}=\{\rho _{0}:X\rightarrow 
\mathbb{R}
,$ $\rho _{1}:X\rightarrow 
\mathbb{R}
\}$ subordinated to $\mathcal{B}=\mathcal{V}$. Moreover, by construction,
for any $\Delta \in \mathbb{T}_{r}$ it results $\rho _{1}|_{\partial \Delta
}\equiv 0$, because $\partial \Delta \subseteq Skel^{r-1}(X,\mathbb{T)}%
\subseteq Z_{\bullet }\subseteq V_{0}$.

Let $\mathcal{R}=\{R_{0},R_{1}\}$ be a honeycomb cell system associated with 
$\mathcal{V}$. Building $R_{0}$ by means of very fine triangulations of $X$,
we can assume without loss of generality that $R_{0}\supseteq Z_{\bullet }$,
that $R_{0}\sim Z_{\bullet }$ and that for any $\Delta \in \mathbb{T}_{r}$
it holds%
\begin{equation}
\Delta \cap R_{(1,0)}\sim \partial \Delta .
\label{R_0,1 e Scheletro Relativo}
\end{equation}%
So, the inclusions $Z_{\bullet }\subseteq R_{0}\subseteq V_{0}$ are\
homotopic equivalences and $R_{1}\cap Z_{\bullet }=\emptyset $, because $%
R_{1}\subseteq V_{1}$.

Let $\mathbf{\phi }$ be the chain map associated with $\mathcal{V}$, $%
\mathcal{B}=\mathcal{V}$ and $\{\rho \}$ (cp. (2) of Remark \ref{Omotopia di
Catene (Relativa)}). Take $\omega \in Z_{ext}^{r}\left( X,Z\right) $ and let 
$\check{\omega}=(\check{\omega}_{0},$ $\check{\omega}_{1},$ $\check{\omega}%
_{01})\in \check{Z}_{ext}^{r}\left( X,\mathcal{V},V_{0}\right) $ be such
that $\omega =\mathbf{\phi }(\check{\omega})$. Namely, by (\ref{Forma
reincollata}),%
\begin{equation}
\omega =\rho _{0}\check{\omega}_{0}+\rho _{1}\check{\omega}_{1}-d\rho
_{0}\wedge \check{\omega}_{01}.  \label{omega e omega-Cech Relativo}
\end{equation}%
Since $\omega \in Z_{ext}^{r}\left( X,Z\right) $ and $\check{\omega}\in 
\check{Z}_{ext}^{r}\left( X,\mathcal{V},V_{0}\right) $, we have $\omega
|_{Z}=0$ and $\check{\omega}_{0}|_{Z}=0$. Moreover, $d\omega =0$ and $D%
\check{\omega}=0$. The latter is%
\begin{equation}
(d\check{\omega}_{0},d\check{\omega}_{1},d\check{\omega}_{01}-\check{\omega}%
_{0}|_{V_{(0,1)}}+\check{\omega}_{1}|_{V_{(0,1)}})=\left( 0,\text{ }0,\text{ 
}0\right) .  \label{d(omega) e D(omega-Cech Relativo)}
\end{equation}

In what follows, with slight abuses of notations, for any $\mathbf{\Delta }%
\in C_{r}^{\mathbb{T}}\left( X,Z\right) $ we write $\mathbf{\Delta }\in
\Delta +C_{p}^{\mathbb{T}}\left( Z\right) $. In the above situation, the
integral of $\omega \in Z_{ext}^{r}\left( X,Z\right) $ over an element $%
\mathbf{\Delta }\in C_{r}^{\mathbb{T}}\left( X,Z\right) $ can be written in
terms of the integrals of the components of $\check{\omega}\in \check{Z}%
_{ext}^{r}\left( X,\mathcal{V},V_{0}\right) $ such that $\omega =\mathbf{%
\phi }(\check{\omega})$. Indeed, the image of $\check{\omega}$ via $\check{%
\eta}_{Z}^{r}$ (cp. (\ref{Cech Integrale Relativo}) and (\ref{Relaz. tra
Oper. d'Integr. Relat. Piu' Esplicita})) is the map $\check{\eta}_{Z}^{r}(%
\check{\omega})\in C_{\mathbb{T}}^{r}\left( X,Z\right) $ that to any $%
\mathbf{\Delta }\in \Delta +C_{r}^{\mathbb{T}}\left( Z\right) $ associates
the number%
\begin{equation}
\check{\eta}_{Z}^{r}(\check{\omega}_{0},\text{ }\check{\omega}_{1},\text{ }%
\check{\omega}_{01})(\mathbf{\Delta })=\tint\nolimits_{\Delta \cap R_{0}}%
\check{\omega}_{0}+\tint\nolimits_{\Delta \cap R_{1}}\check{\omega}%
_{1}-\tint\nolimits_{\Delta \cap R_{(1,0)}}\check{\omega}_{01}
\label{Cech Integrazione Relativo}
\end{equation}

In order to prove (\ref{Cech Integrale Relativo}), recall (\ref{Cech
Integrazione Relativo}). Then%
\begin{equation}
\check{\eta}_{Z}^{r}(\check{\omega})(\mathbf{\Delta })=\eta _{Z}^{r}\circ 
\mathbf{\phi }(\check{\omega})(\mathbf{\Delta })=\tint\nolimits_{\Delta
+C_{r}^{\mathbb{T}}\left( Z\right) }\mathbf{\phi }(\check{\omega})
\label{Cech Integrale Relativo (2)}
\end{equation}%
So, by (\ref{omega e omega-Cech Relativo}),%
\begin{eqnarray}
\check{\eta}_{Z}^{r}(\check{\omega})(\mathbf{\Delta })
&=&\tint\nolimits_{(\Delta +C_{r}^{\mathbb{T}}\left( Z\right) )\cap
R_{0}}(\rho _{0}\check{\omega}_{0}+\rho _{1}\check{\omega}%
_{1})+\tint\nolimits_{(\Delta +C_{r}^{\mathbb{T}}\left( Z\right) )\cap
R_{1}}(\rho _{0}\check{\omega}_{0}+\rho _{1}\check{\omega}_{1})+  \notag \\
&&+\tint\nolimits_{(\Delta +C_{r}^{\mathbb{T}}\left( Z\right) )\cap
R_{0}}d\rho _{1}\wedge \check{\omega}_{01}-\tint\nolimits_{(\Delta +C_{r}^{%
\mathbb{T}}\left( Z\right) )\cap R_{1}}d\rho _{0}\wedge \check{\omega}_{01} 
\notag \\
&=&\tint\nolimits_{(\Delta +C_{r}^{\mathbb{T}}\left( Z\right) )\cap
R_{0}}(\rho _{0}\check{\omega}_{0}+\rho _{1}\check{\omega}%
_{1})+\tint\nolimits_{(\Delta +C_{r}^{\mathbb{T}}\left( Z\right) )\cap
R_{1}}(\rho _{0}\check{\omega}_{0}+\rho _{1}\check{\omega}_{1})+  \notag \\
&&+\tint\nolimits_{(\Delta +C_{r}^{\mathbb{T}}\left( Z\right) )\cap
R_{0}}(d(\rho _{1}\check{\omega}_{01})-\rho _{1}d\check{\omega}_{01})+
\label{Buona Posizione Cech Integrale (0)} \\
&&-\tint\nolimits_{(\Delta +C_{r}^{\mathbb{T}}\left( Z\right) )\cap
R_{1}}(d(\rho _{0}\check{\omega}_{01})-\rho _{0}d\check{\omega}_{01})  \notag
\\
_{\text{Thm. \ref{Teorema di Stokes}, (\ref{d(omega) e D(omega-Cech
Relativo)})}} &=&\tint\nolimits_{(\Delta +C_{r}^{\mathbb{T}}\left( Z\right)
)\cap R_{0}}(\rho _{0}\check{\omega}_{0}+\rho _{1}\check{\omega}%
_{0})+\tint\nolimits_{(\Delta +C_{r}^{\mathbb{T}}\left( Z\right) )\cap
R_{1}}(\rho _{0}\check{\omega}_{1}+\rho _{1}\check{\omega}_{1})+  \notag \\
&&+\tint\nolimits_{\partial ((\Delta +C_{r}^{\mathbb{T}}\left( Z\right)
)\cap R_{0})}\rho _{1}\check{\omega}_{01}-\tint\nolimits_{\partial ((\Delta
+C_{r}^{\mathbb{T}}\left( Z\right) )\cap R_{1})}\rho _{0}\check{\omega}_{01}
\notag
\end{eqnarray}%
Then, by (\ref{Buona Posizione Cech Integrale (0)}),%
\begin{eqnarray}
\check{\eta}_{Z}^{r}(\check{\omega})(\mathbf{\Delta })
&=&\tint\nolimits_{(\Delta +C_{r}^{\mathbb{T}}\left( Z\right) )\cap R_{0}}%
\check{\omega}_{0}+\tint\nolimits_{(\Delta +C_{r}^{\mathbb{T}}\left(
Z\right) )\cap R_{1}}\check{\omega}_{1}+  \notag \\
&&+\tint\nolimits_{\partial ((\Delta +C_{r}^{\mathbb{T}}\left( Z\right)
)\cap R_{0})}\rho _{1}\check{\omega}_{01}-\tint\nolimits_{\partial ((\Delta
+C_{r}^{\mathbb{T}}\left( Z\right) )\cap R_{1})}\rho _{0}\check{\omega}_{01}
\notag \\
_{(1)} &=&\tint\nolimits_{(\Delta +C_{r}^{\mathbb{T}}\left( Z\right) )\cap
R_{0}}\check{\omega}_{0}+\tint\nolimits_{(\Delta +C_{r}^{\mathbb{T}}\left(
Z\right) )\cap R_{1}}\check{\omega}_{1}+  \notag \\
&&+\tint\nolimits_{\partial \Delta +\Delta \cap R_{(0,1)}+B_{r-1}^{\mathbb{T}%
}\left( Z\right) \cap R_{0}}\rho _{1}\check{\omega}_{01}+
\label{Buona Posizione Cech Integrale} \\
&&-\tint\nolimits_{\Delta \cap R_{(1,0)}}\rho _{0}\check{\omega}_{01}  \notag
\\
_{(2)} &=&\tint\nolimits_{(\Delta +C_{r}^{\mathbb{T}}\left( Z\right) )\cap
R_{0}}\check{\omega}_{0}+\tint\nolimits_{(\Delta +C_{r}^{\mathbb{T}}\left(
Z\right) )\cap R_{1}}\check{\omega}_{1}+  \notag \\
&&+\tint\nolimits_{\Delta \cap R_{(0,1)}}\rho _{1}\check{\omega}%
_{01}-\tint\nolimits_{\Delta \cap R_{(1,0)}}\rho _{0}\check{\omega}_{01} 
\notag \\
_{(3)} &=&\tint\nolimits_{(\Delta +C_{r}^{\mathbb{T}}\left( Z\right) )\cap
R_{0}}\check{\omega}_{0}+\tint\nolimits_{(\Delta +C_{r}^{\mathbb{T}}\left(
Z\right) )\cap R_{1}}\check{\omega}_{1}  \notag \\
&&-\tint\nolimits_{\Delta \cap R_{(1,0)}}\check{\omega}_{01},  \notag
\end{eqnarray}%
with $B_{r-1}^{\mathbb{T}}\left( Z\right) =\partial C_{r}^{\mathbb{T}}\left(
Z\right) $. Indeed, $(1)$ follows from $\partial (\Delta \cap
R_{0})=\partial \Delta +\Delta \cap R_{(0,1)}$, $\partial (\Delta \cap
R_{1})=\Delta \cap R_{(1,0)}$ and $R_{1}\cap Z=\emptyset $, $(2)$ follows
from $\rho _{1}|_{Skel^{r-1}(X,\mathbb{T)}}\equiv 0$ and $B_{r-1}^{\mathbb{T}%
}\left( Z\right) \subset Skel^{r-1}(X,\mathbb{T)}$, $(3)$ follows from $%
R_{(0,1)}=-R_{(1,0)}$, \thinspace $\rho _{0}+\rho _{1}\equiv 1$.

Actually, the last term of (\ref{Buona Posizione Cech Integrale}) is
independent of the choice of the representative of $\mathbf{\Delta }\in
\Delta +C_{r}^{\mathbb{T}}\left( Z\right) $. Indeed, if $\Upsilon \in C_{r}^{%
\mathbb{T}}\left( Z\right) $, then $\Upsilon \cap R_{0}\subseteq Z$ and $%
\Upsilon \cap R_{1}=\emptyset $.So, $\tint\nolimits_{\Upsilon \cap R_{0}}%
\check{\omega}_{0}+\tint\nolimits_{\Upsilon \cap R_{1}}\check{\omega}_{1}$,
because $\check{\omega}_{0}|_{Z}\equiv 0$. Then (\ref{Cech Integrale
Relativo}) follows from (\ref{Cech Integrale Relativo (2)}), (\ref{Buona
Posizione Cech Integrale}) and we are done.

Next, $\check{\eta}_{Z}^{r}$ is independent of the chosen honeycomb cell
system $\mathcal{R}$ associated with $\mathcal{V}$. In fact, by (\ref{Relaz.
tra Oper. d'Integr. Relat. Piu' Esplicita}), $\check{\eta}_{Z}^{r}=\eta
_{Z}^{r}\circ \mathbf{\phi }$ and the right hand side does not depend on the
choice of $\mathcal{R}$. Furhtermore, $\check{\eta}_{Z}^{r}$ is independent
of the partition of unity $\{\rho _{0},\rho _{1}\}$ subordinated to $%
\mathcal{V}$. In fact, take into account (\ref{Cech Integrale Relativo}) and
note that its right hand side does not depend on $\{\rho _{0},\rho _{1}\}$.

As a note, for any $\omega \in Z_{ext}^{r}\left( X,Z\right) $ and $\mathbf{%
\Delta }\in \Delta +C_{r}^{\mathbb{T}}\left( Z\right) $\ it results%
\begin{equation}
\tint\nolimits_{\Delta }\omega =\tint\nolimits_{\Delta \cap R_{0}}\check{%
\omega}_{0}+\tint\nolimits_{\Delta \cap R_{1}}\check{\omega}%
_{1}-\tint\nolimits_{\Delta \cap R_{(1,0)}}\check{\omega}_{01},
\label{Varie Scritture dell'Integrale}
\end{equation}%
because of $\eta _{Z}^{r}(\omega )(\mathbf{\Delta })=\eta _{Z}^{r}\circ 
\mathbf{\phi }(\check{\omega})(\mathbf{\Delta })=\check{\eta}_{Z}^{r}(\check{%
\omega})(\mathbf{\Delta })$ (see (\ref{Relaz. tra Oper. d'Integr. Relat.
Piu' Esplicita})).

\begin{remark}
The explicit expression (\ref{Cech Integrazione Relativo}) of $\check{\eta}%
_{Z}^{r}$ is strongly related to properties of $\mathcal{R}$. Indeed, the
hypotheses $R_{0}\tilde{\supseteq}Z_{\bullet }$ and (\ref{R_0,1 e Scheletro
Relativo}) have been given\ in order to achieve Formula (\ref{Cech
Integrazione Relativo}). If $\mathcal{R}^{\prime }$\ is a honeycomb cell
system associated with $\mathcal{V}$ that does not enjoy any particular
property, then the explicit formula for $\check{\eta}_{Z}^{r}$ associated
with $\mathcal{R}^{\prime }$\ is, in general, more complicated than Formula (%
\ref{Cech Integrazione Relativo}). For example, it can happen that the
explicit expression of $\check{\eta}_{Z}^{r}$ depends on the chosen
partition of unity subordinated to $\mathcal{V}$, even if $\check{\eta}%
_{Z}^{r}$ is independent of it (cp. \cite{Suwa Sao Carlos Versione II}).
\end{remark}

\section{Vector bundles\label{Vector bundles (Sezione)}}

\subsection{Extendable Chern classes\label{Ext. Chern classes (Sottosezione)}%
}

In the following, let $\mathbb{K}$ denote either $%
\mathbb{C}
$ or $%
\mathbb{R}
$. Let $X$ be an abstract finite dimensional complex analytic variety. Let $%
E\rightarrow X$ be a differentiable (holomorphic) $\mathbb{K}$-vector bundle
of rank $e$ over $X$ and $\mathcal{C}=\{(A_{l},n_{l},U_{l},W_{l},F_{l})\}_{l%
\in L}$ an atlas of $X$ trivializing $E$. Then for each $l\in L$ there is a
differentiable (holomorphic) real (complex) vector bundle $E_{l}\rightarrow
U_{l}$ of rank $e$ such that $E|_{A_{l}}=(F_{l})^{\ast }(E_{l}|_{W_{l}})$.
For this, just to take the trivial bundle $E_{l}=U_{l}\times \mathbb{K}%
^{e}\rightarrow U_{l}$.

In what follows, we will consider $E|_{X^{\prime }}\rightarrow X^{\prime }$
as an $\mathcal{E}$-extendable bundle, with $\mathcal{E}$ the sheaf of germs
of differentiable sections of $E$ (see Example \ref{Esempio E-->X}). If this
is the case, then $\mathcal{C}=\{A_{l}\}_{l\in L}$, the atlas trivializing $%
E $, is associated with $E|_{X^{\prime }}$ (cp. Remark \ref{Atlante
Associato}). Shrinking the open sets of $\mathcal{C}$, if necessary, we get
the following commutative diagram%
\begin{equation}
\begin{array}{ccccc}
(\mathcal{C}_{%
\mathbb{C}
^{n_{l}}}^{\infty })^{\nu _{l}}|_{U_{l}} & \rightarrow & \mathcal{C}_{%
\mathbb{C}
^{n_{l}}}^{\infty }(%
\mathbb{C}
^{n_{l}},\mathbb{K}^{e})|_{U_{l}} & \rightarrow & 0 \\ 
\downarrow &  & \downarrow &  &  \\ 
((\mathcal{C}_{X}^{\infty })^{\nu _{l}}|_{A_{l}})_{l} & \rightarrow & (%
\mathcal{E}|_{A_{l}})_{l} & \rightarrow & 0 \\ 
\downarrow &  & \downarrow &  &  \\ 
0 &  & 0 &  & 
\end{array}
\label{Miglioramento di ''Trivializing Extension'' (per E)}
\end{equation}%
that improves (\ref{Trivializing Extension}) (cp. Notation \ref{Fasci
(Notazioni)}).

Recall that for any $N,$ $N^{\ast },$ $p\in 
\mathbb{N}
$ the bundle $TX^{\prime \otimes N}\otimes T^{\ast }X^{\prime \otimes
N^{\ast }}\otimes \Lambda ^{p}T^{\ast }X^{\prime }\otimes E|_{X^{\prime }}$
is $(\mathcal{S}_{N,N^{\ast }}\otimes \mathcal{S}_{p}\otimes \mathcal{E})$%
-extendable (cp. Example \ref{Oss. Fibrati Estendibili}).

\begin{definition}
\label{Def. Conness. Est.}Let $X$ be an abstract finite dimensional complex
analytic variety and $E\rightarrow X$ a differentiable (holomorphic) $%
\mathbb{K}$-vector bundle on $X$ of rank $e$. An \emph{extendable linear
connection for }$E$ is a $\mathbb{K}$-linear map $\nabla :\Gamma
_{ext}(E|_{X^{\prime }})\rightarrow \Gamma _{ext}(T^{\mathbb{K}}{}^{\ast
}X^{\prime }\otimes E|_{X^{\prime }})$ such that $\forall f\in \Gamma
_{ext}((X\times \mathbb{K)}|_{X^{\prime }})$ and $\forall s\in \Gamma
_{ext}(E|_{X^{\prime }})$ it holds $\nabla (fs)=df\otimes s+f\nabla s$.
\end{definition}

Let $X$ be an abstract complex analytic variety of complex dimension $n$, $%
E\rightarrow X$ a differentiable (holomorphic) real (complex) vector bundle
of rank $e$ over $X$ and $\nabla $ an extendable connection for $E$.\ If $%
\mathcal{C}=\{A_{l}\}_{l\in L}$ is an atlas of $X$ trivializing $E$ that is
also associated with $T^{\mathbb{K}}X^{\prime }$, then for any $l\in L$
there exists a linear connection $^{l}\tilde{\nabla}:\Gamma
(E_{l})\rightarrow \Gamma (T^{\mathbb{K\ast }}U_{l}\otimes E_{l})$ for $%
E_{l} $ such that $\nabla |_{A_{l}^{\prime }}=(F_{l}|_{A_{l}^{\prime
}})^{\ast }(^{l}\tilde{\nabla})$. Note that the connection forms of $\nabla $
with respect to any given extendable frames of $E|_{A_{l}^{\prime }}$ and $%
T^{\mathbb{K}}X^{\prime }|_{A_{l}^{\prime }}$ are extendable differential
forms.

For this, let $(\partial )=(\partial _{1},$ $...,$ $\partial _{n})$ and $(%
\mathbf{e})=(\mathbf{e}_{1,}$ $...,$ $\mathbf{e}_{e})$ be extendable frames
of $T^{\mathbb{K}}X^{\prime }|_{A_{l}^{\prime }}$ and, respectively, $%
E|_{A_{l}^{\prime }}$. The connection forms $(\theta _{\alpha \beta
})_{\alpha ,\beta }$ of $\nabla |_{A_{l}^{\prime }}$ with respect to $%
(\partial )$ and $(\mathbf{e})$ are defined by $\nabla \mathbf{e}_{\alpha
}=\tsum\nolimits_{\beta }\theta _{\alpha \beta }\otimes \mathbf{e}_{\beta }$%
. Now, let $(\tilde{\partial})=(\tilde{\partial}_{1},$ $...,$ $\tilde{%
\partial}_{n})$ and $(\mathbf{\tilde{e}})=(\mathbf{\tilde{e}}_{1,}$ $...,$ $%
\mathbf{\tilde{e}}_{e})$ be frames of $T^{\mathbb{K}}U_{l}$ and,
respectively, $E_{l}$ extending $(\partial )$ and, respectively, $(\mathbf{e}%
)$ (cp. (\ref{Miglioramento di ''Trivializing Extension''}) and (\ref%
{Miglioramento di ''Trivializing Extension'' (per E)})). Then for any $\iota
\in \{1,..,n\},$ $\upsilon \in \{1,...,e\}$ we have $\partial _{\iota
}=(F_{l}|_{A_{l}^{\prime }})^{\ast }(\tilde{\partial}_{\iota })$ and $%
\mathbf{e}_{\upsilon }=(F_{l}|_{A_{l}^{\prime }})^{\ast }(\mathbf{\tilde{e}}%
_{\upsilon })$. Denote by $(\tilde{\theta}_{\alpha \beta })_{\alpha ,\beta }$
the connection forms of $^{l}\tilde{\nabla}$ with respect to $(\tilde{%
\partial})$ and $(\mathbf{\tilde{e}})$ and recall that these forms are
defined by $^{l}\tilde{\nabla}\mathbf{\tilde{e}}_{\alpha
}=\tsum\nolimits_{\beta }\tilde{\theta}_{\alpha \beta }\otimes \mathbf{%
\tilde{e}}_{\beta }$. Then $\theta _{\alpha \beta }=(F_{l}|_{A_{l}^{\prime
}})^{\ast }(\tilde{\theta}_{\alpha \beta })$ for any $\alpha ,\beta \in
\{1,...,n\}$.

\begin{proposition}
Let $X$ be an abstract finite dimensional complex analytic variety and $%
E\rightarrow X$ a differentiable (holomorphic) real (complex) vector bundle.
Then there exists an extendable linear connection $\nabla $\ for $E$.
\end{proposition}

\begin{proof}
Choose an extendable partition of unity subordinated to a suitable open
covering of $X$ and proceed as for the smooth case.
\end{proof}

Next, let $X$, $E$ and $\nabla $ be as above. $\nabla $ induces the $\mathbb{%
K}$-linear map $\mathbf{\nabla }:\Gamma _{ext}(T^{\mathbb{K\ast }}X^{\prime
}\otimes E|_{X^{\prime }})\rightarrow \Gamma _{ext}(\Lambda ^{2}T^{\mathbb{%
K\ast }}X^{\prime }\otimes E|_{X^{\prime }})$ which to any $\omega \otimes
s\in \Gamma _{ext}(T^{\mathbb{K\ast }}X^{\prime }\otimes E|_{X^{\prime }})$
associates $\mathbf{\nabla }(\omega \otimes s)=d\omega \otimes s-\omega
\wedge \nabla s$ and then it is linearly extended. The map $K^{\nabla
}:\Gamma _{ext}(E|_{X^{\prime }})\rightarrow \Gamma _{ext}(\Lambda ^{2}T^{%
\mathbb{K\ast }}X^{\prime }\otimes E|_{X^{\prime }})$ defined by $K^{\nabla
}=\mathbf{\nabla }\circ \nabla $ is \emph{the curvature of }$\nabla $. As a
note, for any $f\in \Gamma _{ext}((X\times \mathbb{K)}|_{X^{\prime }})$ and
for any $s\in \Gamma _{ext}(E|_{X^{\prime }})$ it results $K^{\nabla
}(fs)=fK^{\nabla }(s)$. So $K^{\nabla }\in \Gamma _{ext}(\Lambda ^{2}T^{%
\mathbb{K\ast }}X^{\prime }\otimes E^{\ast }|_{X^{\prime }}\otimes
E|_{X^{\prime }})$. In order to prove that $\mathbf{\nabla }$ and $K^{\nabla
}$ are well defined, use the classical, explicit, local expressions of $%
\nabla $ and $K^{\nabla }$, noting that the differential forms involved in
such expressions with respect to any local extendable frame of $%
E|_{X^{\prime }}$ and $T^{\mathbb{K}}X^{\prime }$ are extendable (cp. \cite%
{SuwaLibro}, Ch. II, Sec. 7).

Let $n\in 
\mathbb{N}
$. For each $q\in \{1,...,n\}$ denote by $\Sigma _{q}\in 
\mathbb{C}
\lbrack t_{1},...,t_{n}]$ the $q^{th}$ elementary symmetric function in the $%
n$ variables $t_{1},...,t_{n}$ and recall that $\Sigma _{q}$ is a polynomial
of degree $q$.

Let $X$ be an abstract complex analytic variety of complex dimension $n$, $%
E\rightarrow X$ a differentiable (holomorphic) complex vector bundle of rank 
$e$ over $X$ and $\nabla $ an extendable connection for $E$. Proceeding as
for the smooth case, by using the \emph{curvature forms of }$\nabla $, that
is the extendable forms locally representing $K^{\nabla }$, we can associate
a global closed extendable form with any elementary symmetric function (cp. 
\cite{Kobayashi}, Ch. II, Sec. 2). Namely, let $q\in \{1,...,n\}$ and argue
as in \cite{Kobayashi}, Ch. II, Sec. 2 (see also \cite{SuwaLibro}, Ch. II,
Sec. 7). The closed extendable $2q$-form associated with $\nabla $ and $%
\Sigma _{q}$ and achieved in such a way is denoted by $\Sigma _{q}(\nabla
)\in \Gamma _{ext}(\Lambda ^{2q}T^{\mathbb{%
\mathbb{C}
\ast }}X^{\prime })$.

Next, we define \emph{the }$q^{\text{\emph{th}}}$\emph{\ extendable Chern
form associated with }$\nabla $ as the extendable differential $2q$-form $%
c_{ext}^{q}(\nabla )\in \Gamma _{ext}(\Lambda ^{2q}T^{\mathbb{%
\mathbb{C}
\ast }}X^{\prime })$ given by $c_{ext}^{q}(\nabla )=(\sqrt{-1}/2\pi
)^{q}\Sigma _{q}(\nabla )$. For a local expression of $c_{ext}^{q}(\nabla )$
involving curvature forms of $\nabla $, see \cite{Kobayashi}, (2.13). Note
that $c_{ext}^{q}(\nabla )$ is real. Moreover, since it is closed, $%
c_{ext}^{q}(\nabla )$ represents an extendable cohomology class $%
[c_{ext}^{q}(\nabla )]\in H_{ext}^{2q}(X)$.

As in the smooth case, such a class is independent of the choice of the
connection $\nabla $ of $E$. Indeed, if $\nabla _{I},$ $\nabla _{II}$\ are
extendable linear connections for $E$, then there exists $c_{ext}^{q}(\nabla
_{I},\nabla _{II})\in \Gamma _{ext}(\Lambda ^{2q-1}T^{\ast }X)$ such that $%
c_{ext}^{q}(\nabla _{II},\nabla _{I})=-c_{ext}^{q}(\nabla _{I},\nabla _{II})$
and $dc_{ext}^{q}(\nabla _{I},\nabla _{II})=c_{ext}^{q}(\nabla
_{II})-c_{ext}^{q}(\nabla _{I})$. In fact, for the existence of $%
c_{ext}^{q}(\nabla _{I},\nabla _{II})$, proceed as for the classical case.
The extensibility of $c_{ext}^{q}(\nabla _{I},\nabla _{II})$ follows from
its explicit, local expression where only extendable differentiable forms
are involved (cp. \cite{Kobayashi}, Ch. II, Sec. 2, p. 38). $%
c_{ext}^{q}(\nabla _{I},\nabla _{II})$ is \emph{the extendable Bott
difference form with respect to }$c_{ext}^{q}(\nabla _{I})$\emph{\ and }$%
c_{ext}^{q}(\nabla _{II})$.

The extendable cohomology class $c_{ext}^{q}(E)\in H_{ext}^{2q}(X)$ defined
by%
\begin{equation}
c_{ext}^{q}(E)=[c_{ext}^{q}(\nabla )]  \label{Classe di Chern}
\end{equation}%
is \emph{the }$q^{\text{\emph{th}}}$\emph{\ extendable Chern class\ of }$E$.
Set $c_{ext}^{0}(E)=1$. The class $c_{ext}(E)\in \oplus
_{q=1}^{n}H_{ext}^{2q}(X)$ defined by $c_{ext}(E)=\tsum%
\nolimits_{q=0}^{n}c_{ext}^{q}(E)$ is \emph{the total extendable Chern class
of }$E$.

\begin{remark}
\label{Osserv. sulla Forma Differenza Estend.}Let $X$ be a complex analytic
variety of complex dimension $n$, $E\rightarrow X$ a differentiable
(holomorphic) complex vector bundle of rank $e$ and $\mathcal{C}%
=\{A_{l}\}_{l\in L}$ an atlas of $X$ associated with $\Lambda ^{p}T^{\mathbb{%
\mathbb{C}
\ast }}X^{\prime }\otimes E|_{X^{\prime }}$ for any $p\in 
\mathbb{N}
$. Fix $l\in L$ and $q\in \{0,...,n\}$.

\begin{enumerate}
\item Let $\nabla $ be an extendable linear connection for $E$ and $%
^{l}\nabla $ a linear connection for $E_{l}$ such that $\nabla
|_{A_{l}^{\prime }}=(F_{l}|_{A_{l}^{\prime }})^{\ast }(^{l}\nabla )$. Then $%
c_{ext}^{q}(\nabla )\in \Gamma _{ext}(\Lambda ^{2q}T^{\mathbb{%
\mathbb{C}
\ast }}X^{\prime })$ is completely extended on $A_{l}$ by the $q^{th}$
differentiable Chern form $c_{dif}^{q}(^{l}\nabla )\in \Gamma (\Lambda
^{2q}T^{\mathbb{%
\mathbb{C}
\ast }}U_{l})$.

\item Let $\nabla _{I},$ $\nabla _{II}$\ be extendable linear connections
for $E$ and $^{l}\nabla _{I},$ $^{l}\nabla _{II}$ linear connections for $%
E_{l}$ such that $\nabla _{I}|_{A_{l}^{\prime }}=(F_{l}|_{A_{l}^{\prime
}})^{\ast }(^{l}\nabla _{I})$ and $\nabla _{II}|_{A_{l}^{\prime
}}=(F_{l}|_{A_{l}^{\prime }})^{\ast }(^{l}\nabla _{II})$. Then $%
c_{ext}^{q}(\nabla _{I},\nabla _{II})\in \Gamma _{ext}(\Lambda ^{2q-1}T^{%
\mathbb{%
\mathbb{C}
\ast }}X^{\prime })$ is completely extended on $A_{l}$\ by the Bott
difference form $c_{dif}^{q}(^{l}\nabla _{I},^{l}\nabla _{II})\in \Gamma
(\Lambda ^{2q-1}T^{\mathbb{%
\mathbb{C}
\ast }}U_{l})$.
\end{enumerate}
\end{remark}

The following observation is a consequence of Theorem \ref{Isomorfismo tra
Coomologie} and Remark \ref{Osserv. sulla Forma Differenza Estend.}.

\begin{remark}
Let $X$ be an abstract $n$-dimensional complex analytic variety, $%
E\rightarrow X$ a differentiable (holomorphic) complex vector bundle of rank 
$e$ and $\mathcal{V}=\{V_{0},V_{1}\}$ an open covering of $X$. Let $\nabla
_{0}$ and $\nabla _{1}$ be extendable linear connections for $%
E|_{V_{0}}\rightarrow V_{0}$ and, respectively, $E|_{V_{1}}\rightarrow V_{1}$%
. Fix $q\in \{0,...,n\}$ with $q\leq e$ and consider the class%
\begin{equation}
\check{c}_{ext}^{q}(E)=P^{2q\ast }(c_{ext}^{q}(E))\in \check{H}_{ext}^{2q}(X,%
\mathcal{V}).  \label{Classe di Chern-Cech}
\end{equation}%
(see Section \ref{Ext. Cech cohomology groups (Sottosezione)}). The
definition of $P^{2q\ast }:H_{ext}^{2q}(X)\rightarrow \check{H}_{ext}^{2q}(X,%
\mathcal{V})$ implies that $\check{c}_{ext}^{q}(E)\in \check{H}_{ext}^{2q}(X,%
\mathcal{V})$ is represented by the cocycle%
\begin{equation}
\check{c}_{ext}^{q}(\nabla _{\ast })=(c_{ext}^{q}(\nabla
_{0}),c_{ext}^{q}(\nabla _{1}),c_{ext}^{q}(\nabla _{0},\nabla _{1})).
\label{Cociclo di Chern}
\end{equation}
\end{remark}

The following remark is in order.

\begin{remark}
The construction of extendable Chern classes can be extended to more general
cases we do not deal with in this paper. As an example, let $X$ be an
abstract finite dimensional complex analytic variety and $E^{\prime
}\rightarrow X^{\prime }$ an $\mathcal{S}_{E^{\prime }}$-extendable
differentiable (holomorphic) complex vector bundle over the regular part of $%
X$. If the coherent sheaf $\mathcal{S}_{E^{\prime }}$ associated with $%
E^{\prime }$ admits a finite resolution%
\begin{equation*}
0\rightarrow \mathcal{E}_{m}\rightarrow ...\rightarrow \mathcal{E}%
_{0}\rightarrow \mathcal{S}_{E^{\prime }}\rightarrow 0
\end{equation*}%
by locally free sheaves $\mathcal{E}_{0},$ $...,$ $\mathcal{E}_{m}$ of $%
\mathcal{C}_{X}^{\infty }$-modules over $X$, then we can define the total
extendable Chern class $c_{ext}(E^{\prime })$ of $E^{\prime }$ by setting $%
c_{ext}(E^{\prime })=\tprod\nolimits_{\iota \in \{0,...,m\}}c_{ext}(E_{\iota
})^{(-1)^{\iota }}$, with $E_{\iota }$ the vector bundle associated with $%
\mathcal{E}_{\iota }$.
\end{remark}

Now, we consider the relative case.

Let $X$ be a complex analytic variety of complex dimension $n$, $Z$ the
closure of a non empty open set that is also a polyhedron of $X$ and $%
E\rightarrow X$ be a differentiable (holomorphic) complex vector bundle of
rank $e$. Let $\mathbb{T}$\ be a triangulation of $X$ compatible with $%
Sing\left( X\right) $ and $Z$ and $q\in \{0,...,n\}$. Set $r=e-q+1$ and
consider the $(2q-1)$-skeleton $Skel^{2q-1}(X,\mathbb{T})$ of $(X,\mathbb{T}%
) $. Let $Z_{2q-1}$ be a polyhedron of $X$ that is the closure of an open
set and that is such that $Z_{2q-1}\supseteq Skel^{2q-1}\left( X,\mathbb{T}%
\right) $ and $Z_{2q-1}\sim Skel^{2q-1}\left( X,\mathbb{T}\right) $. Set $%
Z_{\bullet }=Z_{2q-1}\cup Z$ and let $V_{Z}$ be an open neighborhood of $%
Z_{\bullet }$ such that $V_{Z}\sim Z_{\bullet }$ .

Let $s^{(r)}$ be a differentiable $r$-section of $E$ whose restriction at $%
V_{Z}$ is an $r$-frame. If $\nabla $ is an extendable linear connection for $%
E$ that is $s^{(r)}$-trivial on $V_{Z}$, then the $q^{th}$\ extendable Chern
form $c_{ext}^{q}(\nabla )$\ vanishes on $V_{Z}$. Namely, $%
c_{ext}^{q}(\nabla )\in \Gamma _{ext}(\Lambda ^{2q}T^{\mathbb{%
\mathbb{C}
\ast }}X^{\prime })_{Z}$. To indicate this fact, we write $%
c_{ext}^{q}(\nabla ,s^{(r)})$ instead of $c_{ext}^{q}(\nabla )$. So, $%
[c_{ext}^{q}(\nabla ,s^{(r)})]\in H_{ext}^{2q}(X)_{Z}$. Moreover, if $\nabla
^{\prime }$ is another extendable linear connection for $E$ that is $s^{(r)}$%
-trivial on $V_{Z}$, then $[c_{ext}^{q}(\nabla ^{\prime
},s^{(r)})]=[c_{ext}^{q}(\nabla ,s^{(r)})]$ as classes in $%
H_{ext}^{2q}(X)_{Z}$. So, $[c_{ext}^{q}(\nabla ,s^{(r)})]\in
H_{ext}^{2q}(X)_{Z}$\ does not depend on the choice of the extendable linear
connection $\nabla $ that is $s^{(r)}$-trivial on $V_{Z}$. To prove this,
proceed as in \cite{SuwaLibro}, Ch. III, Sec. 3. The class $%
c_{ext}^{q}(E,s^{(r)})=[c_{ext}^{q}(\nabla ,s^{(r)})]$ is \emph{the
localization\ outside }$Z$\emph{\ with respect to }$s^{(r)}$\emph{\ of }$%
c_{ext}^{q}(E)$. Note that, as a relative class, $[c_{ext}^{q}(\nabla
,s^{(r)})]$ depends on the frame $s^{\left( r\right) }$ (cp. \cite{SuwaLibro}%
, Ch. III, Sec. 3).

Let $X$, $Z$, $E\rightarrow X$, $q$, $r$, $\mathbb{T}$, $Z_{2q-1}$, $%
Z_{\bullet }$, $V_{Z}$ and $s^{(r)}$ be as above. Set $V_{0}=V_{Z}$\ and let 
$V_{1}$ be an open set in $X$ such that $V_{1}\cap Z_{\bullet }=\emptyset $
and $V_{0}\cup V_{1}=X$. Then the open covering $\mathcal{V}=\{V_{0},V_{1}\}$
of $X$ is adapted to $Z_{\bullet }$ .

Let $\nabla _{0}$ be an $s^{(r)}$-trivial extendable linear connection for $%
E|_{V_{0}}\rightarrow V_{0}$ and $\nabla _{1}$ any extendable linear
connection for $E|_{V_{1}}\rightarrow V_{1}$. By the $s^{(r)}$-triviality of 
$\nabla _{0}$ , we get $c_{ext}^{q}(\nabla _{0})=0$. So, $\check{c}%
_{ext}^{q}(\nabla _{\ast })=(c_{ext}^{q}(\nabla _{0}),$ $c_{ext}^{q}(\nabla
_{1}),$ $c_{ext}^{q}(\nabla _{0},\nabla _{1}))$ lies in $K^{q}(X,\mathcal{V}%
,V_{0})$ (cp. Subsection \ref{Ext. Cech cohomology groups (Sottosezione)}).
In order to indicate this, write $\check{c}_{ext}^{q}(\nabla _{\ast
},s^{(r)})$ instead of $\check{c}_{ext}^{q}(\nabla _{\ast })$. Then $[\check{%
c}_{ext}^{q}(\nabla _{\ast },s^{(r)})]\in \check{H}_{ext}^{2q}\left( X,%
\mathcal{V},V_{0}\right) $. Moreover, it can be proved that $[\check{c}%
_{ext}^{q}(\nabla _{\ast },s^{(r)})]$ does not depend on both the choices of
the $s^{(r)}$-trivial connection $\nabla _{0}$ and of the connection $\nabla
_{1}$ (cp. \cite{SuwaLibro}, Ch. III, Sec. 3).

The class $\check{c}_{ext}^{q}(E,s^{(r)})=[\check{c}_{ext}^{q}(\nabla _{\ast
},s^{(r)})]\in \check{H}_{ext}^{2q}(X,\mathcal{V},V_{0})$ is called \emph{%
the localization\ outside }$Z$\emph{\ with respect to }$s^{(r)}$\emph{\ of }$%
c_{ext}^{q}(E)$. As a reason for this name, note that $\check{c}%
_{ext}^{q}(E,s^{(r)})$ is the image of $c_{ext}^{q}(E,s^{(r)})$ via $%
P^{2q\ast }:H_{ext}^{2q}\left( X,Z\right) \rightarrow \check{H}%
_{ext}^{2q}\left( X,\mathcal{V},V_{0}\right) $. Note that, as a relative
class, $[\check{c}_{ext}^{q}(\nabla _{\ast },s^{(r)})]$ depends on the frame 
$s^{\left( r\right) }$.

Next, we wish to study the following case.

Let $X$ be a complex analytic variety of complex dimension $n$ and $%
E\rightarrow X$ be a holomorphic complex vector bundle of rank $e$. Take $%
q\in \{0,...,n\}$ and set $r=e-q+1$. Let be $s^{\left( r\right) }$ a
holomorphic $r$-section of $E$ and denote by $S$ its singular locus. Since $%
s^{\left( r\right) }$ is holomorphic, $S$ is a closed complex analytic
subvariety of $X$. Furthermore, $S$ is a polyhedron of $X$.

Now, suppose that $S$ is compact and let $U_{S}$ be an open neighborhood of $%
S$ in $X$ enjoying (1) and (2) of Lemma \ref{Sottovarieta' Compatte}. We can
assume without loss of generality that the closure $\overline{U_{S}}$ of $%
U_{S}$ is also compact. Set $Z=X\setminus U_{S}$ and note that the
restriction $s^{(r)}|_{X\setminus U_{S}}$ of $s^{(r)}$\ at $X\setminus U_{S}$
is an $r$-frame, because $X\setminus U_{S}\subseteq X\setminus S$. In the
above situation, both $c_{ext}^{q}(E,s^{(r)})\in H_{ext}^{2q}(X)_{Z}$ and $%
\check{c}_{ext}^{q}(E,s^{(r)})\in \check{H}_{ext}^{2q}(X,\mathcal{V},V_{Z})$%
, the localizations outside $Z$ with respect to $s^{(r)}$\ of $%
c_{ext}^{q}(E) $, are called \emph{localization\ at }$S$\emph{\ with respect
to }$s^{(r)}$\emph{\ of }$c_{ext}^{q}(E)$.

\subsection{The homomorphisms $\boldsymbol{P}_{k}^{\ast }$ and $\boldsymbol{A%
}_{S,k}^{\ast }$}

Let $X$ be a compact irreducible complex analytic variety of complex
dimension $n$ and $[X]$ the fundamental class of $X$ (cp. Remark \ref%
{Topologia Algebrica (Osservazione)}). For each $k\in \{0,...,2n\}$ the cap
product with $[X]$ induces a homomorphism $\boldsymbol{P}_{k}^{\ast
}:H^{k}\left( X\right) \rightarrow H_{2n-k}\left( X\right) $ called \emph{%
the }$k$\emph{-Poincar\'{e} homomorphism}. Namely, $\boldsymbol{P}_{k}^{\ast
}\left( [c]\right) =[c]\cap \lbrack X]$ for any $[c]\in H^{k}\left( X\right) 
$.

Let $X$ be a compact irreducible complex analytic variety of complex
dimension $n$. If $\mathbb{T}$\ is a finite triangulation of $X$ compatible
with $Sing(X)$ and if an orientation of $(X,\mathbb{T)}$ is already given,
then it can be proved that $\boldsymbol{P}_{k}^{\ast }$ is induced by%
\begin{equation}
\begin{array}[b]{rccl}
\boldsymbol{P}_{k}: & C_{\mathbb{T}^{\prime }}^{k}\left( X\right) & 
\rightarrow & C_{2n-k}^{\mathbb{T}}\left( X\right) \\ 
& c & \mapsto & \boldsymbol{P}_{k}\left( c\right) =\tsum\nolimits_{\Delta
\in \mathbb{T}_{2n-k}}c(\hat{\Delta})\text{ }\Delta%
\end{array}%
,  \label{Poincare' Somma}
\end{equation}%
where $\mathbb{T}^{\prime }$ is the first barycentric subdivision of $%
\mathbb{T}$, $\mathbb{T}^{\ast }$ is the dual block decomposition of $(X,%
\mathbb{T)}$ and $\hat{\Delta}$ denotes the dual block of $\Delta \in 
\mathbb{T}_{2n-k}$ (cp. \cite{Munkres}, Ch. 8, Sec. 64).

Let $X$ be a irreducible complex analytic variety of complex dimension $n$, $%
S$ a closed compact analytic subvariety of $X$ and $\mathbb{T}$\ a finite
triangulation of $X$ compatible with $Sing\left( X\right) $ and $S$. If $%
c\in C_{\mathbb{T}^{\prime }}^{k}\left( X,X\setminus S\right) $ is a $k$%
-cochain relative to $\left( X,X\setminus S\right) $, then in a sum as in (%
\ref{Poincare' Somma}) only appear the simplices in $S$. Namely,%
\begin{equation}
\tsum\nolimits_{\Delta \in \mathbb{T}_{2n-k}:\Delta \subseteq S}c(\hat{\Delta%
})\text{ }\Delta  \label{Alexander-Lefschetz Somma}
\end{equation}%
Such a finite sum induces a homomorphism $\boldsymbol{A}_{S,k}^{\ast
}:H^{k}\left( X,X\setminus S\right) \rightarrow H_{2n-k}\left( S\right) $
called \emph{the }$(S,k)$\emph{-Alexander-Lefschetz homomorphism}.

The following proposition is a direct consequence of the constructions of $%
\boldsymbol{P}_{k}$\ and $\boldsymbol{A}_{S,k}^{\ast }$ (cp. \cite{Munkres}).

\begin{proposition}
\label{Poincare' e Alexander Commutano}Let $X$ be a compact irriducible
complex analytic variety of dimension $n$ and $S$ a compact analytic
subvariety of $X$. Let $k\in \left\{ 0,...,2n\right\} $. Then the following
diagram is commutative%
\begin{equation*}
\begin{array}{ccccc}
H^{k}\left( X,X\setminus S\right) &  & \overset{j^{\ast }}{\rightarrow } & 
& H^{k}\left( X\right) \\ 
&  &  &  &  \\ 
\downarrow _{\text{ }^{\boldsymbol{A}_{S,k}^{\ast }}} &  &  &  & \downarrow
_{\text{ }^{\boldsymbol{P}_{k}^{\ast }}} \\ 
&  &  &  &  \\ 
H_{2n-k}\left( S\right) &  & \overset{i_{\ast }}{\rightarrow } &  & 
H_{2n-k}\left( X\right)%
\end{array}%
\end{equation*}
\end{proposition}

\subsection{Topological Chern classes\label{Topological Chern classe for
extendable v.b. via obstruction theory}}

We briefly recall the construction of topological Chern classes for
continuous complex vector bundles (cp. \cite{Steenrod}).

Take $e\in 
\mathbb{N}
$, let $r\in 
\mathbb{N}
$ be such that $r\leq e$ and denote by $W_{r}(%
\mathbb{C}
^{e})$ the complex Stiefel manifold of $r$-frames. Recall that $W_{r}(%
\mathbb{C}
^{e})$ is $(2e-2r)$-connected and that $\pi _{2e-2r+1}(W_{r}(%
\mathbb{C}
^{e}))\simeq 
\mathbb{Z}
$ has a canonical generator $\mathbf{\varsigma }$ (cfr. \cite{Steenrod},
25.7).

\begin{definition}
\label{Def. Sezioni "Topologiche"}Let $X$ be a finite dimensional complex
analytic variety, $Y\subseteq X$ a subset of $X$ and $E\rightarrow X$ a
continuous complex vector bundle of rank $e$.

An $r$\emph{-section of }$E\rightarrow X$ \emph{on }$Y$ is an ordered family 
$s^{(r)}=(s_{1},...,s_{r})$ of $r$ continuous sections of $E\rightarrow X$
over $Y$. A \emph{singular point of }$s^{\left( r\right) }=(s_{1},...,s_{r})$
is a point where the sections $s_{1},$ $...,$ $s_{r}$ fail to be linearly
independent over $%
\mathbb{C}
$. An $r$\emph{-frame (frame) of }$E\rightarrow X$ \emph{on }$Y$ is a
non-singular $r$-section ($e$-section) over $Y$.
\end{definition}

Let $X$ be a finite dimensional complex analytic variety, $E\rightarrow X$ a
continuous complex vector bundle of rank $e$ and $W_{r}(E)$ the continuous
bundle of complex $r$-frames of $E$ over $X$. Recall that $W_{r}(E)$ is a
bundle associated with $E$ whose fibre $W_{r}(E_{x})$ over a point $x\in X$
is homeomorphic to $W_{r}(%
\mathbb{C}
^{e})$.

Set $q=e-r+1$ and note that $2e-2r+1=2q-1$. The primary obstruction to
constructing a section of $W_{r}(E)$ is \emph{the }$q^{th}$\emph{\
topological Chern class }$c^{q}\left( E\right) $\emph{\ of }$E$. A priori
this class is defined in the cohomology with local coefficients system
formed by $\pi _{2q-1}(W_{r}(E_{x}))$. But, since the group of
transformations of $E\rightarrow X$ is the unitary group, that is arcwise
connected, this system is in fact trivial (cp. \cite{Steenrod}, 30.4). Then
the obstruction class $c^{q}\left( E\right) $ lies in the integral
cohomology group $H^{2q}\left( X\right) $.

Next, we introduce the notion of index of an $r$-section at its singular
point.

Let $U$ be an open set in $X$ such that $E|_{U}$ is trivial and $\kappa
:U\times W_{r}\left( 
\mathbb{C}
^{e}\right) \rightarrow W_{r}\left( 
\mathbb{C}
^{e}\right) $ the projection on the second factor. Let $B^{2q}\subset U$ be
a suitably oriented ball of dimension $2q$ in $U$, $p\in B^{2q}$ be a point
in the interior of $B^{2q}$ and $S^{2q-1}$ the boundary of $B^{2q}$. As a
note, $S^{2q-1}$ is an oriented $\left( 2q-1\right) $-sphere. If $s^{\left(
r\right) }$ is an $r$-section of $E$ on a neighborhood of $B^{2q}$ in $U\ $%
and $p$ is an isolated singularity of $s^{\left( r\right) }$, then the map%
\begin{equation*}
\psi =\psi _{p}:S^{2q-1}\overset{s^{\left( r\right) }}{\rightarrow }%
W_{r}(E)|_{U}\simeq U\times W_{r}(%
\mathbb{C}
^{e})\overset{\kappa }{\rightarrow }W_{r}(%
\mathbb{C}
^{e})
\end{equation*}%
given by $\psi =\kappa \circ s^{\left( r\right) }$ is well defined.

By the very definition of homotopy groups, the map $\psi $\ determines an
element $\mathbf{\psi }$ in $\pi _{2e-2r+1}(W_{r}(%
\mathbb{C}
^{e}))$. Furthermore, since $\pi _{2e-2r+1}(W_{r}(%
\mathbb{C}
^{e}))\simeq 
\mathbb{Z}
$ $\mathbf{\varsigma }$, there exists an integer $I(E|_{B^{2q}},s^{\left(
r\right) },p)\in 
\mathbb{Z}
$\ such that%
\begin{equation}
\mathbf{\psi }=I(E|_{B^{2q}},s^{\left( r\right) },p)\text{ }\mathbf{%
\varsigma }.  \label{Def. Indice}
\end{equation}%
Such an integer is \emph{the index of the }$r$\emph{-section }$s^{\left(
r\right) }$\emph{\ of }$E\rightarrow X$\emph{\ at the point }$p$.

Let $X$ be a finite dimensional complex analytic variety, $E\rightarrow X$ a
continuous complex vector bundle and $\mathbb{T}$ a triangulation of $X$. We
try to construct an $r$-frame of $E$, that is a section of $W_{r}\left(
E\right) $, on each skeleton of $\mathbb{T}$\ inductively from the skeleton
of dimension $0$.

First of all, it is always possible to construct a section $s^{\left(
r\right) }$ of $W_{r}\left( E\right) $\ over $Skel^{0}\left( X,\mathbb{T}%
\right) $.

Next, if\ $\Delta _{h}\in \mathbb{T}_{h}$ is contained in an open set $U$ of 
$X$ trivializing $W_{r}(E)$ and an $r$-section $s^{\left( r\right) }$ is
given on $\partial \Delta _{h}$, then there is a well defined map $\kappa
\circ s^{\left( r\right) }:S^{2h-1}\rightarrow W_{r}(E)|_{U}\simeq U\times
W_{r}(%
\mathbb{C}
^{e})\rightarrow W_{r}(%
\mathbb{C}
^{e})$ that determines an element in $\pi _{h-1}(W_{r}(%
\mathbb{C}
^{e}))$, because $\partial \Delta _{h}$ is homeomorphic to a $\left(
2h-1\right) $-sphere $S^{2h-1}$. Then, by the very definition of homotopy
groups, if $h\leq 2q-1=2e-2r+1$, the section $s^{\left( r\right) }$ can be
extended to a continuous $r$-section defined on the interior of $\Delta _{h}$%
, because for any $h\in 
\mathbb{N}
$ such that $h\leq 2e-2r+1$ it results $\pi _{h-1}(W_{r}(%
\mathbb{C}
^{e}))=0$.

Iterating this process, we reach an obstruction for $h=2q$. Namely, the $r$%
-frame $s^{\left( r\right) }$ can be extended to an $r$-section on a
suitably small $\Delta _{2q}\in \mathbb{T}_{2q}$ with at most a singularity
at an interior point $p$ of $\Delta $ and the measure of such an obstruction
is given by the index $I(E|_{\Delta },s^{\left( r\right) },p)$ of $s^{\left(
r\right) }$ at $p$. So, up to choose a $p_{\Delta }$ in the interior of any $%
\Delta \in \mathbb{T}_{2q}$, we may define a cochain $\Gamma _{\mathbb{T}%
}^{2q}\in C_{\mathbb{T}}^{2q}(X)$ associated with $s^{(r)}$. Namely, $\Gamma
_{\mathbb{T}}^{2q}$ is the cochain whose value at $\Delta \in \mathbb{T}%
_{2q} $ is%
\begin{equation}
\Gamma _{\mathbb{T}}^{2q}(\Delta )=I(E|_{\Delta },s^{\left( r\right)
},p_{\Delta })  \label{Def. Cociclo}
\end{equation}%
and then it is extended by linearity. It can be proved that $\Gamma _{%
\mathbb{T}}^{2q}$ is a cocycle and that the cohomology class represented by $%
\Gamma _{\mathbb{T}}^{2q}$ is independent of all the choices made in the
definition (cp. \cite{Steenrod}).

The class identified by $\Gamma _{\mathbb{T}}^{2q}$ is denoted by $%
c_{top}^{q}(E)\in H^{2q}(X)$ and it is called \emph{the }$q^{th}$\emph{\
topological Chern class of }$E\rightarrow X$. The class $c_{top}(E)=\tsum%
\nolimits_{q=0}^{e}c_{top}^{q}(E)$, with $c_{top}^{0}(E)=1$, is \emph{the
total topological Chern class }$c_{top}(E)$\emph{\ of }$E\rightarrow X$. As
a note, $c_{top}(E)$ is an invertible element of the cohomology ring $\oplus
_{r\in 
\mathbb{N}
}H^{r}\left( X\right) $.

\begin{remark}
Let $X$ be a compact irreducible complex analytic variety of complex
dimension $n$, $\mathbb{T}$ a finite triangulation of $X$ and $q\in
\{0,...,n\}$. Then, by (\ref{Poincare' Somma}), the class $\boldsymbol{P}%
_{2q}^{\ast }(c_{top}^{q}(E))\in H_{2n-2q}\left( X\right) $ is represented
by the cycle%
\begin{equation}
\tsum\nolimits_{\Delta \in \mathbb{T}_{2n-2q}}I(E|_{\Delta },s^{\left(
r\right) },p_{\Delta })\text{ }\Delta .  \label{P(c(E))}
\end{equation}
\end{remark}

We need observations about the localization of topological Chern classes.

Let $X$ be a complex analytic variety of complex dimension $n$, $Z$ a
polyhedron of $X$ that is the closure of a non empty open set and $%
E\rightarrow X$ a continuous complex vector bundle of rank $e$. Let $\mathbb{%
T}$ be a triangulation of $X$ compatible with $Z$ and $Sing(X)$ and consider 
$Z$ as a closed subcomplex of $(X,\mathbb{T})$.

Let $q\in \{0,...,n\}$ and set $r=e-q+1$. If $s^{(r)}$ is a continuous $r$%
-section of $E\rightarrow X$ whose restriction at $Z$ is an $r$-frame, then $%
I(E|_{\Delta },s^{\left( r\right) },p_{\Delta })=0$ for any $\Delta \in 
\mathbb{T}_{2q}$ contained in $Z$. So, in this case, the cocycle $\Gamma _{%
\mathbb{T}}^{2q}\in C_{\mathbb{T}}^{2q}\left( X\right) $ representing $%
c_{top}^{q}(E)$\ lies in $C_{\mathbb{T}}^{2q}\left( X,Z\right) $ and it
represents a class $c_{top}^{q}(E,s^{\left( r\right) })\in H^{2q}\left(
X,Z\right) $. The class $c_{top}^{q}(E,s^{\left( r\right) })$ is \emph{the
localization\ outside }$Z$ \emph{with respect to }$s^{(r)}$ \emph{of }$%
c_{top}^{q}(E)$. Indeed, the image\ of $c_{top}^{q}(E,s^{\left( r\right) })$%
\ via the map $j^{\ast }:H^{2q}\left( X,Z\right) \rightarrow H^{2q}\left(
X\right) $ induced by the inclusion $j:(X,\emptyset )\rightarrow (X,Z)$ is $%
c_{top}^{q}\left( E\right) $. Note that, as a relative class, $%
c_{top}^{q}(E,s^{\left( r\right) })$ depends on the frame $s^{\left(
r\right) }$ (cp. \cite{Steenrod}).

Next, we study the following case (see the end of Subsection \ref{Ext. Chern
classes (Sottosezione)}).

Let $X$ be a complex analytic variety of complex dimension $n$ and $%
E\rightarrow X$ be a continuous complex vector bundle of rank $e$. Take $%
q\in \{0,...,n\}$ and set $r=e-q+1$. Let be $s^{\left( r\right) }$ a
continuous $r$-section of $E$ and denote by $S$ its singular locus. Suppose
that $S$ is a closed complex analytic subvariety of $X$ and note that, in
this case, $S$ is a polyhedron of $X$.

Now, suppose that $S$ is compact. Let $U_{S}$ be an open neighborhood of $S$
enjoying (1) and (2) of Remark \ref{Sottovarieta' Compatte} and such that
its closure $\overline{U_{S}}$ is also compact. Set $Z=X\setminus U_{S}$ and
note that the restriction of $s^{(r)}$ at $Z$ is an $r$-frame, because $%
Z\subset X\setminus S$. In the above situation, $c_{top}^{q}(E,s^{(r)})\in
H^{2q}(X,Z)$, the localization outside $Z$ with respect to $s^{(r)}$\ of $%
c_{top}^{q}(E)$, is called \emph{localization\ at }$S$\emph{\ with respect
to }$s^{(r)}$\emph{\ of }$c_{top}^{q}(E)$. Furthermore, since $(X,Z)\sim
(X,X\setminus S)$, $c_{top}^{q}(E,s^{\left( r\right) })$ corresponds to a
class in $H^{2q}(X,X\setminus S)$, still denoted by $c_{top}^{q}(E,s^{\left(
r\right) })$.

\begin{definition}
Let $X$ be an irreducible complex analytic variety of complex dimension $n$
and $E\rightarrow X$ a continuous complex vector bundle of rank $e$. Take $%
q\in \{0,...,n\}$ with $q\leq e$ and set $r=e-q+1$. Let $s^{\left( r\right)
} $ be a continuous $r$-section of $E\rightarrow X$ and $S$ its singular
locus. Suppose that $S$ is a closed compact complex analytic subvariety of $%
X $. \emph{The topological residue of }$s^{\left( r\right) }$ \emph{for }$%
c_{top}^{q}(E)$\emph{\ at }$S$ is the homology class $%
TopRes_{c_{top}^{q}}(E,s^{\left( r\right) },S)\in H_{2n-2q}(S)$ defined by%
\begin{equation}
TopRes_{c_{top}^{q}}(E,s^{\left( r\right) },S)=\boldsymbol{A}_{S,2q}^{\ast
}(c_{top}^{q}(E,s^{\left( r\right) }))
\label{Residuo Topologico (Definizione)}
\end{equation}
\end{definition}

We have the following remark (cp. Proposition \ref{Poincare' e Alexander
Commutano}).

\begin{remark}
\label{Poincare' e Alexander Commutano (per Residui)}Let $X$ be a compact
irreducible complex analytic variety of complex dimension $n$ and $%
E\rightarrow X$ a continuous complex vector bundle of rank $e$. Take $q\in
\{0,...,n\}$ with $q\leq e$ and set $r=e-q+1$. Let $s^{\left( r\right) }$ be
a continuous $r$-section of $E\rightarrow X$ and $S$ its singular locus.
Suppose that $S$ is a closed compact complex analytic subvariety of $X$.
Then, by Proposition \ref{Poincare' e Alexander Commutano},%
\begin{equation}
i_{\ast }(TopRes_{c_{top}^{q}}(E,s^{\left( r\right) },S))=\boldsymbol{P}%
_{2q}^{\ast }\left( c_{top}^{q}(E)\right)  \label{TopRes e P(c(E))}
\end{equation}
\end{remark}

By (\ref{P(c(E))})\ and (\ref{Alexander-Lefschetz Somma}), $%
TopRes_{c_{top}^{q}}(E,s^{\left( r\right) },S)$ is represented by the cycle%
\begin{equation}
\tsum\nolimits_{\Delta \in \mathbb{T}_{2n-2q}:\Delta \subset S}I(E|_{\Delta
},s^{\left( r\right) },p_{\Delta })\text{ }\Delta
\label{Formula Esplicita per Loc. Omol. (Topologica)}
\end{equation}%
If $S$ has a finite number of connected components $\{S_{\nu }\}_{\nu \in
\{1,...,N\}}$, then, correspondingly to the decomposition $%
H_{2n-2q}(S)=\oplus _{\nu \in \{1,...,N\}}H_{2n-2q}(S_{\nu })$, for any $\nu
\in \{1,...,N\}$ we get the residue $TopRes_{c_{top}^{q}}(E,s^{\left(
r\right) },S_{\nu })$.

\subsection{Localization of differentiable Chern classes\label{SEZIONE
Localized Diff. Chern Classes}}

We consider the localization of differentiable Chern classes of a smooth
complex vector bundle defined over a manifold, assuming that the
construction of such classes is already familiar to the reader. For the
necessary background of differential geometry, refer to \cite{SuwaLibro}.

Let $M$ be an oriented differentiable manifold of real dimension $m$. Then $%
M $ is a triangulable space. From now on, suppose that a triangulation $%
\mathbb{T}$ of $M$ is already given and still denote by $M$ the simplicial
complex associated with $M$\ and $\mathbb{T}$. As a matter of notations, let 
$S$ be a subcomplex of $M$ and $k\in \{0,...,m\}$. The $k^{th}$ de Rham
cohomology group of $M$ and the $k^{th}$ de Rham cohomology group relative
to $(M,M\setminus S)$ are denoted by $H_{dR}^{k}(M)$ and, respectively, $%
H_{dR}^{k}(M,M\setminus S)$. Recall that $H_{dR}^{2q}(M)\simeq H^{2q}(M)$
and $H_{dR}^{2q}(M,M\setminus S)\simeq H^{2q}(M,M\setminus S)$.

Let $E\rightarrow M$ be a differentiable complex vector bundle of rank $e$
over $M$ and for each $q\in \{0,...,[\frac{m}{2}]\}$ denote by $%
c_{dif}^{q}(E)\in H_{dR}^{2q}(M)$ the $q^{th}$ differentiable Chern class of 
$E$.

Next, take $r\in \{0,...,e\}$ and let $q\in \{0,...,[\frac{m}{2}]\}$ be such
that $q\geq e-r+1$. Let $S$ be a subcomplex of $M$. If $s^{(r)}$ is a
differentiable $r$-section of $E\rightarrow M$ whose restriction at $%
M\setminus S$ is an $r$-frame, then\ there exists a cohomology class $%
c_{dif}^{q}(E,s^{(r)})\in H_{dR}^{2q}(M,M\setminus S)$ called \emph{the
localization at }$S$ \emph{with respect to }$s^{(r)}$\emph{\ of }$%
c_{dif}^{q}(E)$.

To prove this, set $D_{0}=M\setminus S$ and let $D_{1}$ be an open
neighborhood of $S$ such that $D_{1}\sim S$. Then $\mathcal{D}%
=\{D_{0},D_{1}\}$ is an open covering\ of $M$. Let $\nabla _{0},$ $\nabla
_{1}$ be differentiable linear connections for $E|_{D_{0}}\rightarrow D_{0}$
and, respectively, $E|_{D_{1}}\rightarrow D_{1}$ and consider the \v{C}%
ech-de Rham cocycle $\check{c}_{dif}^{q}(\nabla _{\ast
})=(c_{dif}^{q}(\nabla _{0}),$ $c_{dif}^{q}(\nabla _{1}),$ $%
c_{dif}^{q}(\nabla _{0},\nabla _{1}))$ representing $c_{dif}^{q}(E)$. Now,
if $\nabla _{0}$\ is $s^{(r)}$-trivial, then $c_{dif}^{q}(\nabla _{0})=0$
and $\check{c}_{dif}^{q}(\nabla _{\ast })$ represents a cohomology class $[%
\check{c}_{dif}^{q}(\nabla _{\ast },s^{(r)})]\in H_{dR}^{2q}(M,M\setminus S)$%
. It can be proved that $[\check{c}_{dif}^{q}(\nabla _{\ast },s^{(r)})]$
does not depend on the choices of both $\nabla _{1}$ and $\nabla _{0}$ that
is $s^{(r)}$-trivial (cp. \cite{SuwaLibro}, Ch. III, Sec. 3). The class $%
c_{dif}^{q}(E,s^{(r)})\in H_{dR}^{2q}(M,M\setminus S)$ is defined by $%
c_{dif}^{q}(E,s^{(r)})=[\check{c}_{dif}^{q}(\nabla _{\ast },s^{(r)})]$. As a
note, $j^{\ast }(c_{dif}^{q}(E,s^{(r)}))=c_{dif}^{q}(E)$, with $j^{\ast
}:H_{dR}^{2q}(M,M\setminus S)\rightarrow H_{dR}^{2q}(M)$ induced by the
inclusion $j:(M,\emptyset )\rightarrow (M,M\setminus S)$. Notice that, as a
relative class, $c_{dif}^{q}(E,s^{(r)})$ depends on the frame $s^{\left(
r\right) }$ (cp. \cite{SuwaLibro}, Ch. III, Sec. 3).

Now, suppose that $S$ is also compact. Let $\{S_{\nu }\}_{\nu \in
\{1,...,N\}}$ be the finite set of connected components of $S$ and consider
the $(S,2q)$-Alexander-Lefschetz duality $\boldsymbol{A}_{S,2q}^{\ast
}:H^{2q}(M,M\setminus S)\overset{\simeq }{\rightarrow }\oplus _{\nu \in
\{1,...,N\}}H_{m-2q}(S_{\nu })$. For each $\nu \in \{1,...,N\}$ \emph{the
differential geometric residue of }$s^{(r)}$\emph{\ for }$c_{dif}^{q}$\emph{%
\ at }$S_{\nu }$ is the homology class $DifRes_{c_{dif}^{q}}(s^{(r)},E,S_{%
\nu })\in H_{m-2q}(S_{\nu })$ defined by%
\begin{equation}
DifRes_{c_{dif}^{q}}(s^{(r)},E,S_{\nu })=\boldsymbol{A}_{S,2q}^{\ast
}(c_{dif}^{q}(E,s^{\left( r\right) })).  \label{DifRes}
\end{equation}

In order to give a description of such a residue, let $D_{0},$ $D_{1}$ be as
above. For any $\nu \in \{1,...,N\}$ let $\breve{D}_{\nu }$ be an open
neighborhood of $S_{\nu }$ in $D_{1}$ and $R_{\nu }$ an $m$-dimensional
manifold with differentiable boundary such that $S_{\nu }\subset \mathring{R}%
_{\nu }\subset R_{\nu }\subset \breve{D}_{\nu }$. Assume that for any $\nu
_{1},\nu _{2}\in \{1,...,N\}$ such that $\nu _{1}\neq \nu _{2}$ it holds $%
\breve{D}_{\nu _{1}}\cap \breve{D}_{\nu _{2}}=\emptyset $. Then the residue $%
DifRes_{c_{dif}^{q}}(s^{(r)},E,S_{\nu })$ is represented by an $(m-2q)$%
-cycle $C_{\nu }$ in $S_{\nu }$ such that for any closed differentiable $%
(m-2q)$-form $\tau \in \Gamma (\Lambda ^{m-2q}T^{\ast }\breve{D}_{\nu })$ it
holds%
\begin{equation}
\tint\nolimits_{C_{\nu }}\tau =\tint\nolimits_{R_{\nu }}(c_{dif}^{q}(\nabla
_{1})\wedge \tau )+\tint\nolimits_{-\partial R_{\nu }}(c_{dif}^{q}(\nabla
_{0},\nabla _{1})\wedge \tau )  \label{DifRes e Integrali}
\end{equation}%
If $2q=m$, then the differential geometric residue can be identified with
the complex number given by%
\begin{eqnarray}
DifRes_{c_{dif}^{m}}(s^{(r)},E,S_{\nu }) &=&\tint\nolimits_{R_{\nu
}}c_{dif}^{m}(\nabla _{1})+\tint\nolimits_{-\partial R_{\nu
}}c_{dif}^{m}(\nabla _{0},\nabla _{1})
\label{Espressione del Res. Classico in p} \\
&=&\tint\nolimits_{R_{\nu }}c_{dif}^{m}(\nabla
_{1})-\tint\nolimits_{\partial R_{\nu }}c_{dif}^{m}(\nabla _{0},\nabla _{1})
\notag
\end{eqnarray}

Next, let $M$ and $E\rightarrow M$ be as above. We give a topological
interpretation of the residue in case the compact subset $S$ of $M$ is just
a point $p$. Namely, $S=\{p\}$. In order to proceed, we need some
definitions. Let $z_{1},...,z_{m}$ denote the complex coordinates on $%
\mathbb{C}
^{m}$. Write $\Theta (z)=dz_{1}\wedge ...\wedge dz_{m}$ and for each $h\in
\{1,...,m\}$ set $\Theta _{h}(z)=(-1)^{h}z_{h}(dz_{1}\wedge ...\wedge 
\widehat{dz_{h}}\wedge ...\wedge dz_{m})$. \emph{The Bochner-Martinelli
kernel on }$%
\mathbb{C}
^{m}$ is $(2m-1)$-differential form $\beta _{m}:%
\mathbb{C}
^{m}\rightarrow \Lambda ^{2m-1}T^{\ast }%
\mathbb{C}
^{m}$ defined by $\beta _{m}(z)=\frac{(m-1)!}{(2\pi \sqrt{-1})^{m}}\cdot 
\frac{(-1)^{\frac{m(m-1)}{2}}}{\left\Vert z\right\Vert ^{2m}}\left(
\tsum\nolimits_{h=1}^{m}\overline{\Theta _{h}(z)}\wedge \Theta (z)\right) $.

Let $B^{2m}\subset 
\mathbb{R}
^{2m}\simeq 
\mathbb{C}
^{m}$ be a ball of real dimension $2m$ and $E\rightarrow D$\ a
differentiable complex vector bundle of rank $m$ defined on a open
neighborhood $D$ of $B^{2m}$. Suppose there is a non vanishing
differentiable section $s$ of $E\rightarrow D$ defined on an open
neighborhood of the boundary $S^{2m-1}$ of $B^{2m}$. Suppose also that the $%
s $ is defined on the whole of $D$ with at most an isolated singularity at a
point $p$ in the interior of $B^{2m}$. Indeed, this can be done, by
dimensional reasons. Then, on one hand, we have the index $%
I(E|_{B^{2q}},s,p)\in 
\mathbb{Z}
$ (cp. (\ref{Def. Indice})) and, on the other hand, we have $%
DifRes_{c_{dif}^{m}}(s,E,p)\in H_{0}(p)\simeq 
\mathbb{C}
$, the differential geometric residue of $s$ for $c_{dif}^{m}$\ at $\{p\}$.

We claim that they are in fact the same number. To prove this, suppose that $%
E\rightarrow D$ is trivial and let $\mathbf{e}^{(m)}=(\mathbf{e}_{1},...,%
\mathbf{e}_{m})$ be a frame of $E\rightarrow D$\ on $D$. Then $%
s=\tsum\nolimits_{h=1}^{m}f_{h}\mathbf{e}_{h}$, with $f_{h}:D\rightarrow 
\mathbb{C}
$ a differentiable function for any $h\in \{1,...,m\}$. Then we have a
differentiable map $f=(f_{1},...,f_{m}):D\rightarrow 
\mathbb{C}
^{m}$ that takes the value $(0,...,0)\in 
\mathbb{C}
^{m}$ only at $p\in B^{2m}$. In the above situation it results $%
DifRes_{c_{dif}^{m}}(s,E,p)=\tint\nolimits_{S^{2m-1}}f^{\ast }(\beta _{m})$.
So, in particular%
\begin{equation}
DifRes_{c_{dif}^{m}}=I(E|_{B^{2m}},s,p)
\label{Res. Classico e Indice Semplice}
\end{equation}%
(cp. \cite{SuwaLibro}, Ch. III, Sec. 4).

Now, we consider a more general case. Let $B^{2m}\subset 
\mathbb{R}
^{2m}\simeq 
\mathbb{C}
^{m}$ be a ball of real dimension $2m$ and $D$ a open neighborhood of $%
B^{2m} $. Denote by $S^{2m-1}$ the boundary of $B^{2m}$ and let $D^{\prime }$%
\ be an open neighborhood of $S^{2m-1}$. Let $E\rightarrow D$ be\ a
differentiable complex vector bundle of rank $e$. Set $r=e-m+1$ and let $%
s^{(r)}=(s_{1},...,s_{r})$ be a differentiable $r$-section of $E\rightarrow
D $ with at most an isolated singularity at a point $p$ in the interior of $%
B^{2m}$ and such that its restriction at $D^{\prime }$ is an $r$-frame. In
fact, as before, this can be done for dimensional reasons. Let $t\in 
\mathbb{N}
$ be such that $t<r$, write $s^{(t)}=(s_{1},...,s_{t})$ and set $%
s^{(r-t)}=(s_{t+1},...,s_{r})$. Suppose that $s^{(t)}$\ is non singular on
the whole of $D$. So, $s^{(t)}$ generates a trivial complex vector bundle $%
E^{t}\rightarrow D$ of rank $t$ that is a subbundle of $E\rightarrow D$.
Suppose that $E\rightarrow D$ is the trivial vector bundle. Then we have the
following exact sequence%
\begin{equation*}
0\rightarrow E^{t}\rightarrow E\rightarrow E_{\bullet }\rightarrow 0,
\end{equation*}%
with $E_{\bullet }\rightarrow D$\ a complex vector bundle of rank $e-t$.
Denote by $s_{\bullet }^{(r-t)}$ the $(r-t)$-section of $E_{\bullet
}\rightarrow D$ induced by $s^{(r-t)}$. \ Then $s_{\bullet }^{(r-t)}$ has at
most an isolated singularity at $p\in B^{2m}$. In the above situation it
results%
\begin{equation}
DifRes_{c_{dif}^{m}}(s^{(r)},E,p)=DifRes_{c_{dif}^{m}}(s_{\bullet
}^{(r-t)},E_{\bullet },p)  \label{Res. Classico e Indice Complice}
\end{equation}%
and%
\begin{equation}
DifRes_{c_{dif}^{m}}(s^{(r)},E,p)=I(E|_{B^{2m}},s^{(r)},p)
\label{Res. Classico e Indice}
\end{equation}%
Indeed, if $t=r-1$, then (\ref{Res. Classico e Indice}) follows from (\ref%
{Res. Classico e Indice Semplice}) and (\ref{Res. Classico e Indice Complice}%
). For a proof of (\ref{Res. Classico e Indice Complice}), refer to \cite%
{SuwaLibro}, Ch. III, Sec. 4.

\subsection{Topological and extendable Chern classes}

We describe the topological Chern classes by means of the extendable Chern
classes.

\begin{theorem}
\label{Classi Estendibili e Topologiche coincidono copy(2)}Let $X$ be a
complex analytic variety of complex dimension $n$ and $Z$ either the empty
set or the closure of a non empty open set that is a polyhedron of $X$. Let $%
E\rightarrow X$ be a differentiable (holomorphic) complex vector bundle of
rank $e$. Take $q\in \{0,...,n\}$ with $q\leq e$ and set $r=e-q+1$.

Let $\mathcal{C}=\{A_{l}\}_{l\in L}$\ be an atlas of $X$ associated with $%
\Lambda ^{p}T^{\mathbb{%
\mathbb{C}
\ast }}X^{\prime }\otimes E|_{X^{\prime }}$ for any $p\in 
\mathbb{N}
$. Let $\mathbb{T}$\ be a $\mathcal{C}$-small triangulation of $X$
compatible with $Sing\left( X\right) $ and $Z$. Let $Z_{2q-1}$ be a
polyhedron of $X$ which is the closure of an open set so that $%
Z_{2q-1}\supset Skel^{2q-1}\left( X,\mathbb{T}\right) $ and $Z_{2q-1}\sim
Skel^{2q-1}\left( X,\mathbb{T}\right) $. Set $Z_{\bullet }=Z_{2q-1}\cup Z$
and let $V_{0}$ be an open neighborhood of $Z_{\bullet }$ such that $%
V_{0}\sim Z_{\bullet }$ and such that for any $\Delta \in \mathbb{T}_{2q}$
for which $\Delta \nsubseteq Z$\ it holds $\Delta \cap V_{0}=\Delta
\setminus \{p_{\Delta }\}$, where $p_{\Delta }$ is the barycentre of $\Delta 
$.

Let $s^{(r)}$ be a differentiable (holomorphic) $r$-section of $E$ whose
restriction at $V_{0}$ is an $r$-frame. Consider the localizations $%
c_{top}^{q}(E,s^{(r)})$\ and $c_{ext}^{q}(E,s^{(r)})$\ outside $Z\subset
V_{0}$\ with respect to $s^{(r)}$\ of $c_{top}^{q}(E)$ and, respectively, $%
c_{ext}^{q}(E)$. Then%
\begin{equation*}
c_{top}^{q}(E,s^{(r)})=H_{Z}^{2q}(c_{ext}^{q}(E,s^{(r)})).
\end{equation*}
\end{theorem}

\begin{proof}
Let $V_{1}$ be an open set of $X$ such that $V_{1}\cap Z_{\bullet
}=\emptyset $ and $V_{0}\cup V_{1}=X$. Then the open covering $\mathcal{V}%
=\{V_{0},V_{1}\}$ is adapted to $Z$, because $Z\subseteq Z_{\bullet }$ .

Let $\nabla _{0}$ be an $s^{(r)}|_{V_{0}}$-trivial extendable linear
conection for $E|_{V_{0}}$ and $\nabla _{1}$ an extendable linear connection
for $E|_{V_{1}}$. Since $\nabla _{0}$\ is $s^{(r)}$-trivial and $r=e-q+1$,
it results $c_{ext}^{q}(\nabla _{0})=0$. Then, by $Z\subseteq Z_{\bullet }$,
the cocycle $\check{c}_{ext}^{q}(\nabla _{\ast })=(c_{ext}^{q}(\nabla _{0}),$
$c_{ext}^{q}(\nabla _{1}),$ $c_{ext}^{q}(\nabla _{0},\nabla _{1}))$
represents both $\check{c}_{ext}^{q}(E,s^{(r)})\in \check{H}_{ext}^{2q}(X,%
\mathcal{V},V_{0})$ and $c_{ext}^{q}(E,s^{(r)})\in H_{ext}^{2q}(X,Z)$.

Let $\{\rho _{0},$ $\rho _{1}:X\rightarrow 
\mathbb{R}
\}$ be an extendable partition of unity subordinated to $\mathcal{V}$. Then
for any $\Delta \in \mathbb{T}_{2q}$ we have $\rho _{1}|_{\partial \Delta
}\equiv 0$, because $\partial \Delta \subseteq Skel^{2q-1}\left( X,\mathbb{T}%
\right) \subset Z_{2q-1}\subseteq Z_{\bullet }$ . Let $\mathcal{R}%
=\{R_{0},R_{1}\}$ be a honeycomb cell system associated with $\mathcal{V}$
and suppose that $R_{0}\supset Z_{\bullet }$, that $R_{0}\sim Z_{\bullet }$
and that for any $\Delta \in \mathbb{T}_{2q}$ it holds $R_{(1,0)}\cap \Delta
\sim \partial \Delta $. Then the inclusions $Z_{\bullet }\subset
R_{0}\subset V_{0}$ are homotopic equivalences.

Let $\Gamma _{\mathbb{T}}^{2q}\in C_{\mathbb{T}}^{2q}(X)$ be the cocycle
associated with $s^{(r)}$ and representing the class $c_{top}^{q}(E)\in
H^{2q}(X)$. Such a cocycle is defined by assigning to each $\Delta \in 
\mathbb{T}_{2q}$ the value $\Gamma _{\mathbb{T}}^{2q}(\Delta )=I(E|_{\Delta
},s^{\left( r\right) },p_{\Delta })$ (cp. (\ref{Def. Cociclo}) and (\ref%
{Def. Indice})). Furthermore, $\Gamma _{\mathbb{T}}^{2q}$ belongs to $C_{%
\mathbb{T}}^{2q}(X,Z)$ and it represents the class $c_{top}^{q}(E,s^{\left(
r\right) })\in H^{2q}(X,Z)$. Indeed, $s^{\left( r\right) }|_{Z}$ is an $r$%
-frame. So, $I(E|_{\Delta },s^{\left( r\right) },p_{\Delta })=0$ for any $%
\Delta \in \mathbb{T}_{2q}$ such that $\Delta \subseteq Z$ (cp. Subsection %
\ref{Topological Chern classe for extendable v.b. via obstruction theory}).

As a matter of notations, for any $C\in C_{2q}^{\mathbb{T}}\left( X\right) $
set $\mathbf{C}=C+C_{2q}^{\mathbb{T}}\left( Z\right) $. Consider the
operator of integration $\check{\eta}_{Z}^{2q}:\check{Z}_{ext}^{2q}(X,%
\mathcal{V},V_{0})\rightarrow Z_{\mathbb{T}}^{2q}(X,Z)$ (cp. (\ref{Cech
Integrale Relativo}) and (\ref{Relaz. tra Oper. d'Integr. Relat. Piu'
Esplicita})). We claim that for any $\Delta \in \mathbb{T}_{2q}$ it results $%
\check{\eta}_{Z}^{2q}(\check{c}_{ext}^{q}(\nabla _{\ast }))(\mathbf{\Delta }%
)=\Gamma _{\mathbb{T}}^{2q}(\mathbf{\Delta })$.

First of all, if $\Upsilon \in C_{2q}^{\mathbb{T}}\left( Z\right) $, then $%
\check{\eta}_{Z}^{2q}(\check{c}_{ext}^{q}(\nabla _{\ast }))(\mathbf{\Upsilon 
})=0$ (cp. Subsection \ref{Integrazione delle Classi di Coom. di Cech Est.
(Sottosezione)}). So, for any $\Upsilon \in C_{2q}^{\mathbb{T}}\left(
Z\right) $ it results%
\begin{equation}
\check{\eta}_{Z}^{2q}(\check{c}_{ext}^{q}(\nabla _{\ast }))(\mathbf{\Upsilon 
})=\Gamma _{\mathbb{T}}^{2q}(\mathbf{\Upsilon }),  \label{Uguaglianza 1}
\end{equation}%
because $\Gamma _{\mathbb{T}}^{2q}|_{C_{2q}^{\mathbb{T}}\left( Z\right)
}\equiv 0$.

Next, let $\Delta $ be any simplex in $\mathbb{T}_{2q}$. We claim that%
\begin{equation}
\tint\nolimits_{\Delta \cap R_{1}}c_{ext}^{q}(\nabla
_{1})-\tint\nolimits_{\Delta \cap R_{(1,0)}}c_{ext}^{q}(\nabla _{0},\nabla
_{1})=I(E|_{\Delta },s^{\left( r\right) },p_{\Delta })  \label{Uguaglianza 2}
\end{equation}%
and this will be enough to conclude, because of the hypotheses on $\mathcal{R%
}$\ (see (\ref{Cech Integrazione Relativo}) and, more generally, Subsection %
\ref{Integrazione delle Classi di Coom. di Cech Est. (Sottosezione)}).

To prove (\ref{Uguaglianza 2}), we proceed locally. Take $\Delta \in \mathbb{%
T}_{2q}$ and let $l\in L$ be such that $A_{l}\supset \Delta $. Let $%
s_{l}^{\left( r\right) }$ be a differentiable $r$-section of $%
E_{l}\rightarrow U_{l}$ extending $s^{(r)}|_{A_{l}}$ on $A_{l}$. Let $%
^{l}\nabla _{0}$ be an $s_{l}^{\left( r\right) }$-trivial differentiable
linear connection for $E_{l}$ extending $\nabla _{0}$ and $^{l}\nabla _{1}$
a differentiable linear connection for $E_{l}$ extending $\nabla _{1}$.

As a matter of notations, write $\Lambda $, $\Lambda \cap R_{1}$, $\Lambda
\cap R_{(1,0)}$ and $p_{\Lambda }$ instead of $F_{l}(\Delta )$, $%
F_{l}(\Delta \cap R_{1})$, $F_{l}(\Delta \cap R_{(1,0)})$ and $%
F_{l}(p_{\Delta })$. Then, by Proposition \ref{Integrale per
Triangolazione(Omega)} and Remark \ref{Osserv. sulla Forma Differenza
Estend.}, it suffices to prove that%
\begin{equation}
\tint\nolimits_{\Lambda \cap R_{1}}c_{dif}^{q}(^{l}\nabla
_{1})-\tint\nolimits_{\Lambda \cap R_{(1,0)}}c_{dif}^{q}(^{l}\nabla
_{0},^{l}\nabla _{1})=I(E_{l}|_{\Lambda },s_{l}^{\left( r\right)
},p_{\Lambda })  \label{Uguaglianza Decisiva}
\end{equation}%
Actually, (\ref{Uguaglianza Decisiva}) follows from (\ref{Res. Classico e
Indice}) and (\ref{Espressione del Res. Classico in p}), because $U_{l}$ is
a differentiable complex manifold. Then, by (\ref{Uguaglianza 1}) and (\ref%
{Uguaglianza 2}), the cocycles $\check{\eta}_{Z}^{2q}(\check{c}%
_{ext}^{q}(\nabla _{\ast }))\in C_{\mathbb{T}}^{2q}(X,Z)\ $and $\Gamma _{%
\mathbb{T}}^{2q}\in C_{\mathbb{T}}^{2q}(X.Z)$ coincide and we are done.
\end{proof}

If $Z$ is empty, then the following theorem holds.

\begin{theorem}
\label{Classi Estendibili e Topologiche coincidono copy(1)}Let $X$ be an
abstract complex analytic variety of complex dimension $n$ and $E\rightarrow
X$ a differentiable (holomorphic) complex vector bundle of rank $e$. Take $%
q\in \{1,...,n\}$ with $q\leq e$. Then%
\begin{equation*}
c_{top}^{q}(E)=H^{2q}(c_{ext}^{q}(E)).
\end{equation*}
\end{theorem}

As an application of Theorem \ref{Classi Estendibili e Topologiche
coincidono copy(2)}, we prove an abstract residue theorem. For the necessary
background on topological Chern classes and their residues, see \cite{Suwa
Sao Carlos Versione II}, Ch. 1.

Let $X$ be a complex analytic variety of complex dimension $n$ and $%
E\rightarrow X$ a holomorphic complex vector bundle of rank $e$. Take $q\in
\{0,...,n\}$ with $q\leq e$ and set $r=e-q+1$. Let $s^{\left( r\right) }$ be
a holomorphic $r$-section of $E$ and denote by $S$ its singular locus. Then $%
S$ is a closed complex analytic subvariety of $X$ that is also a polyhedron
of $X$.

Suppose that $S$ is compact and let $U_{S}$ be an open neighborhood of $S$
in $X$ enjoying (1) and (2) of Lemma \ref{Sottovarieta' Compatte} such that
its closure $\overline{U_{S}}$ is also compact. Set $Z=X\setminus U_{S}$ and
let $c_{top}^{q}(E,s^{(r)})\in H^{2q}(X,Z)$ and $\check{c}%
_{ext}^{q}(E,s^{(r)})\in \check{H}_{ext}^{2q}(X,Z)$ be the localization\ at $%
S$\ with respect to $s^{(r)}$\ of $c_{top}^{q}(E)$ and, respectively, $%
c_{ext}^{q}(E)$. Then, by Theorem \ref{Classi Estendibili e Topologiche
coincidono copy(2)}, $c_{top}^{q}(E,s^{(r)})=H_{Z}^{2q}\left(
c_{ext}^{q}(E,s^{(r)})\right) $ Furthermore, since $(X,Z)\sim \left(
X,X\setminus S\right) $, $c_{top}^{q}(E,s^{\left( r\right) })\in H^{2q}(X,Z)$
corresponds to a class in $H^{2q}(X,X\setminus S)$ still denoted by $%
c_{top}^{q}(E,s^{\left( r\right) })$.

If $X$ is compact and irriducible, then the following abstract residue
theorem holds (see Proposition \ref{Poincare' e Alexander Commutano} and
Remark \ref{Poincare' e Alexander Commutano (per Residui)}).

\begin{theorem}
\label{Teor. dei Res.} \textbf{(Residue theorem)}\ Let $X$ be a compact and
irreducible complex analytic variety of complex dimension $n$ and $%
E\rightarrow X$ a holomorphic complex vector bundle of rank $e$. Take $q\in
\{0,...,n\}$ with $q\leq e$ and set $r=e-q+1$. Let $s^{\left( r\right) }$ be
a holomorphic $r$-section of $E$ and $S$ the singular locus of $s^{(r)}$.
Then $i_{\ast }(TopRes_{c_{top}^{q}}(E,s^{\left( r\right) },S))=\boldsymbol{P%
}_{2q}^{\ast }\circ H^{2q}(c_{ext}^{q}\left( E\right) )$. If $q=n$, then%
\begin{equation}
i_{\ast }(TopRes_{c_{top}^{n}}(E,s^{\left( r\right)
},S))=\tint\nolimits_{[X]}c_{ext}^{n}(E)  \label{Formula dei Res.}
\end{equation}
\end{theorem}

As a note, under the hypotheses of Theorem \ref{Teor. dei Res.}, using the
notations employed at the end of the proof of Theorem \ref{Classi
Estendibili e Topologiche coincidono copy(2)}, the right hand side of (\ref%
{Formula dei Res.}) can be written as%
\begin{equation}
\tint\nolimits_{\lbrack X]}c_{ext}^{n}(E)=\left[ \tsum\nolimits_{\Delta \in 
\mathbb{T}_{2n}}\left( \tint\nolimits_{\Lambda \cap
R_{1}}c_{dif}^{n}(^{l}\nabla _{1})-\tint\nolimits_{\Lambda \cap
R_{(1,0)}}c_{dif}^{n}(^{l}\nabla _{0},^{l}\nabla _{1})\right) \text{ }%
p_{\Delta }\right]  \label{Residuo Astratto Esplicito}
\end{equation}%
In fact, (\ref{Residuo Astratto Esplicito}) follows from (\ref%
{Alexander-Lefschetz Somma}), because of (\ref{Formula Esplicita per Loc.
Omol. (Topologica)}), (\ref{Uguaglianza 2}) and (\ref{Uguaglianza Decisiva})
(cp. the proof of Theorem \ref{Classi Estendibili e Topologiche coincidono
copy(2)}).

\subsection{Generalized Camacho-Sad index theorem\label{Generalized
Camacho-Sad theorem (Sottosezione)}}

As a matter of notations, the stalk at $a$ of a sheaf $\mathcal{S}%
\rightarrow A$\ is denoted by $\mathcal{S}_{a}$. Let $X$ be an abstract
finite dimensional complex analytic variety. From now on, the sheaf of germs
of holomorphic vector fields on $X$ is denoted by $\mathcal{T}X$ instead of $%
\mathcal{O}_{X}(\mathbf{T}X)$. For the necessary background about
foliations, refer to \cite{Bracci} and \cite{SuwaLibro}, Ch. VI, Sec. 6.

Let $X$ be an abstract complex analytic variety of complex dimension $n$ and 
$Y$ a complex analytic subvariety of $X$ of complex dimension $m\lneq n$
such that $Y\nsubseteq Sing(X)$. Set $Y^{\prime }=Y\setminus ((Sing(X)\cap
Y)\cup Sing(Y))$.

Let $\mathcal{F}$\ be a holomorphic foliation of $X$ of rank $k\leq m$ and
write $Sing(\mathcal{F})=\{x\in X^{Reg}:(\mathcal{T}X/\mathcal{F})_{x}$ $is$ 
$locally$ $free\}\cup Sing(X)$. Then $\mathcal{F}|_{X\setminus Sing(\mathcal{%
F})}$ is the sheaf of holomorphic sections of a holomorphic vector bundle $F$
over $X\setminus Sing(\mathcal{F})$.

Suppose that $Y$\ is $\mathcal{F}$-invariant. Then the image of the sheaf
homomorphism $\mathcal{F}\otimes \mathcal{O}_{Y}\rightarrow \mathcal{T}%
X\otimes \mathcal{O}_{Y}$, still denoted by $\mathcal{F}\otimes \mathcal{O}%
_{Y}$, is a holomorphic foliation of $Y$ of rank $k$ and $(\mathcal{F}%
\otimes \mathcal{O}_{Y})|_{Y^{Reg}}$\ is a possibly singular foliation of
the manifold $Y^{Reg}$. Consider the following exact sequence of sheaves $%
0\rightarrow \mathcal{F}\otimes _{\mathcal{O}_{X}}\mathcal{O}_{Y}\rightarrow 
\mathcal{T}X\otimes _{\mathcal{O}_{X}}\mathcal{O}_{Y}\rightarrow \mathcal{Q}%
\otimes _{\mathcal{O}_{X}}\mathcal{O}_{Y}\rightarrow 0$, set $S=(Sing(%
\mathcal{F})\cap Y)\cup Sing(Y)$ and write $Y^{\prime \prime }=Y\setminus S$%
. Then $(\mathcal{F}\otimes _{\mathcal{O}_{X}}\mathcal{O}_{Y})|_{Y^{\prime
\prime }}$\ is the sheaf of holomorphic sections of a holomorphic vector
bundle. By $Y^{\prime \prime }=Y\setminus S\subseteq Y^{\prime }$, we have
the following diagram%
\begin{equation}
\begin{array}{ccccccccc}
0 & \rightarrow & (\mathcal{F}\otimes _{\mathcal{O}_{X}}\mathcal{O}%
_{Y})|_{Y^{\prime \prime }} & \rightarrow & (\mathcal{T}X\otimes _{\mathcal{O%
}_{X}}\mathcal{O}_{Y})|_{Y^{\prime \prime }} & \rightarrow & (\mathcal{Q}%
\otimes _{\mathcal{O}_{X}}\mathcal{O}_{Y})|_{Y^{\prime \prime }} & 
\rightarrow & 0 \\ 
&  &  &  & \downarrow &  &  &  &  \\ 
0 & \rightarrow & \mathcal{T}Y|_{Y^{\prime \prime }} & \rightarrow & (%
\mathcal{T}X\otimes _{\mathcal{O}_{X}}\mathcal{O}_{Y})|_{Y^{\prime \prime }}
& \overset{\pi }{\rightarrow } & \mathcal{N}_{Y}|_{Y^{\prime \prime }} & 
\rightarrow & 0%
\end{array}
\label{Diagramma Foliazione}
\end{equation}

\begin{remark}
Let $M$ be a complex manifold. For definitions and general results
concerning partial connections for a complex vector bundle $E\rightarrow M$,
see \cite{Abate Bracci Tovena 2} and \cite{Baum-Bott}. Let $H\subseteq 
\mathbf{T}M$ be an involutive holomorphic bundle. For definitions and
results about $H$-bundles and (flat) holomorphic actions of $H$ on a given
holomorphic vector bundle over $M$, refer to \cite{Abate Bracci Tovena 2}
and \cite{Baum-Bott}.
\end{remark}

Denote by $N_{Y^{\prime }}\rightarrow Y^{\prime }$\ the complex vector
bundle associated with $\mathcal{N}_{Y}|_{Y^{\prime }}$. $N_{Y^{\prime }}$
is the normal bundle of $Y^{\prime }$. It is known that $N_{Y^{\prime
}}|_{Y^{\prime \prime }}\rightarrow Y^{\prime \prime }\ $is an $%
(F|_{Y^{\prime \prime }})$-vector bundle with respect to the map%
\begin{equation}
\begin{array}{cccl}
\tau : & \Gamma (F|_{Y^{\prime \prime }})\times \Gamma (N_{Y^{\prime
}}|_{Y^{\prime \prime }}) & \rightarrow & \Gamma (N_{Y^{\prime
}}|_{Y^{\prime \prime }}) \\ 
& (f,s) & \mapsto & \tau (f,s)=\pi ([\tilde{f},\text{ }\tilde{s}%
]|_{Y^{\prime \prime }})%
\end{array}
\label{Form. Azione Olom. Tau}
\end{equation}%
where $\tilde{f}$ and $\tilde{s}$ are sections of $\Gamma (\mathbf{T}%
X|_{X\setminus (Sing(\mathcal{F})\cup Sing(Y))})$ such that $\tilde{f}%
|_{Y^{\prime \prime }}=f$ and, respectively, $\pi (\tilde{s}|_{Y^{\prime
\prime }})=s$. Furthermore, $\tau $ is a flat holomorphic action of $%
F|_{Y^{\prime \prime }}$\ on $N_{Y^{\prime }}|_{Y^{\prime \prime }}$ (cp. 
\cite{Abate Bracci Tovena 2}).

Let $\nabla $ be a linear connection of type $(1,0)$ for $N_{Y^{\prime
}}|_{Y^{\prime \prime }}\rightarrow Y^{\prime \prime }$ extending the
partial connection $(F_{Y^{\prime \prime }}\oplus \mathbf{\bar{T}}%
Y|_{Y^{\prime \prime }},$ $\tau \oplus \bar{\partial})$. Denote by $K$ the
curvature of $\nabla $. Then for any symmetric homogeneous polynomial $\Phi
\in 
\mathbb{C}
\lbrack t_{1},...,t_{n}]$ of degree $q\in \{m-k+1,...,m\}$ it results $\Phi
(K)=0$. In particular, $c_{ext}^{q}(\nabla )=c_{dif}^{q}(\nabla )=0$ for any 
$q\in \{m-k+1,...,m\}$ (cp. \cite{Abate Bracci Tovena 2}, Theorem 6.1).

Suppose that $N_{Y^{\prime }}|_{Y^{\prime \prime }}\rightarrow Y^{\prime
\prime }$ is the restriction at $Y^{\prime \prime }=Y\setminus S$ of a
holomorphic vector bundle $N_{Y}\rightarrow Y$ defined over the whole of $Y$%
. In this case, $N_{Y}|_{Y^{\prime }}=N_{Y^{\prime }}=(\mathbf{T}%
X^{Reg}|_{Y^{\prime }})/\mathbf{T}Y^{\prime }$ and so $N_{Y}|_{Y^{\prime
\prime }}=N_{Y^{\prime }}|_{Y^{\prime \prime }}$. Suppose that $Y$ is
compact and globally irreducible. Then it is possible to localize some
extendable Chern classes of $Y$ around $S$.

To prove this, let $V_{0}$ be the open subset of $Y$ defined as $%
V_{0}=Y^{\prime \prime }=Y\setminus S$. Let $V_{1}$ be a neighborhood of $S$
open in $Y$, homotopically equivalent to $S$ and such that $\overline{V_{1}}$
is compact. Set $Z=Y\setminus V_{1}$ and suppose that both $\overline{V_{1}}$
and $Z$ are polytopes with respect to a triangulation of $X$. Then $\mathcal{%
V}=\{V_{0},V_{1}\}$ is an open covering of $Y$ adapted to $Z$. Let $\nabla
_{0}$ be a linear connection of type $(1,0)$ for $N_{Y^{\prime
}}|_{V_{0}}\rightarrow V_{0}$ extending the partial connection $%
(F_{Y^{\prime \prime }}\oplus \mathbf{\bar{T}}Y|_{Y^{\prime \prime }},$ $%
\tau \oplus \bar{\partial})$. Let $\nabla _{1}$ be an extendable linear
connection for $N_{Y}|_{V_{1}}\rightarrow V_{1}$. Take $q\in \{m-k+1,...,m\}$%
. Then $\check{c}_{ext}^{q}(\nabla _{\ast })=(c_{ext}^{q}(\nabla _{0}),$ $%
c_{ext}^{q}(\nabla _{1}),$ $c_{ext}^{q}(\nabla _{0},\nabla _{1}))=(0,$ $%
c_{ext}^{q}(\nabla _{1}),$ $c_{ext}^{q}(\nabla _{0},\nabla _{1}))$, because $%
c_{ext}^{q}(\nabla _{0})=c_{dif}^{q}(\nabla _{0})=0$ (cp. \cite{Abate Bracci
Tovena 2}). So $[\check{c}_{ext}^{q}(\nabla _{\ast })]\in \check{H}%
_{ext}^{r}\left( Y,\mathcal{V},V_{0}\right) $.

Consider the following commutative diagram%
\begin{equation}
\begin{array}{ccccccccc}
H_{ext}^{2q}\left( Y,Z\right) &  & \overset{H_{Z}^{2q}}{\longrightarrow } & 
& H^{2q}(Y,Y\setminus S) &  & \overset{\boldsymbol{A}_{S,2q}^{\ast }}{%
\longrightarrow } &  & H_{2n-2q}(S) \\ 
\downarrow &  &  &  & \downarrow &  &  &  & \downarrow _{\text{ }^{i_{\ast
}}} \\ 
H_{ext}^{2q}\left( Y\right) &  & \overset{H^{2q}}{\longrightarrow } &  & 
H^{2q}(Y) &  & \overset{\boldsymbol{P}_{2q}^{\ast }}{\longrightarrow } &  & 
H_{2n-2q}(Y)%
\end{array}
\label{Diagramma Residuo-Foliazione}
\end{equation}%
denote by $c_{ext}^{q}(N_{Y},\mathcal{F},Y\setminus Z)\in H_{ext}^{2q}\left(
Y,Z\right) $ the cohomology class corresponding to $[\check{c}%
_{ext}^{q}(\nabla _{\ast })]\in \check{H}_{ext}^{r}\left( Y,\mathcal{V}%
,V_{0}\right) $ and set $Res_{c_{ext}^{q}}(N_{Y},\mathcal{F},S)=\boldsymbol{A%
}_{S,2q}^{\ast }\circ H_{Z}^{2q}(c_{ext}^{q}(N_{Y},\mathcal{F},Y\setminus
Z)) $.

\begin{theorem}
\label{Teorema Camacho-Sad}Let $X$ be an abstract complex analytic variety
of complex dimension $n$ and $Y$ a compact and globally irreducible complex
analytic subvariety of $X$ of complex dimension $m\lneq n$ such that $%
Y\nsubseteq Sing(X)$. Let $\mathcal{F}$\ be a holomorphic foliation of $X$
of rank $k\leq m$ and suppose that $Y$ is $\mathcal{F}$-invariant. Set $%
S=(Sing(\mathcal{F})\cap Y)\cup Sing(Y)$ and write $Y^{\prime \prime
}=Y\setminus S$. Let $N_{Y}\rightarrow Y$ be a holomorphic vector bundle on $%
Y$ whose restriction at $Y^{\prime \prime }$ coincides with the normal
bundle of $Y^{\prime \prime }$. Then $\boldsymbol{P}_{2q}^{\ast }\circ
H^{2q}(c_{ext}^{q}\left( N_{Y}\right) )=i_{\ast }(Res_{c_{ext}^{q}}(N_{Y},%
\mathcal{F},S))$. If $q=m$, then%
\begin{equation}
\tint\nolimits_{\lbrack Y]}c_{ext}^{m}(N_{Y})=i_{\ast
}(Res_{c_{ext}^{m}}(N_{Y},\mathcal{F},S))  \label{Res. Camacho-Sad}
\end{equation}
\end{theorem}

We have the following remark.

\begin{remark}
As an example of a complex vector bundle $N_{Y}\rightarrow Y$\ enjoying the
above hypotheses, consider the restriction at $Y$ of the line bundle $%
L_{Y}\rightarrow X$ canonically associated with a Cartier divisor $Y$ of the
ambient variety $X$. Another class of examples is given by subvarieties $Y$
of $X\ $which are defined as the zero locus of a section of a holomorphic
vector bundle defined over the ambient variety $X$. The last example is
similar to the construction for subvarieties of complex manifolds which are
also strongly local complete intersection (cp. \cite{Lehmann-Suwa}).
\end{remark}

Next, we give an explicit expression of $i_{\ast }(Res_{c_{ext}^{q}}(N_{Y},%
\mathcal{F},S))$\ in a simple but fundamental case. Let $X$ be an abstract
complex analytic variety of complex dimension $2$ and $Y$ be a compact and
globally irreducible Cartier divisor of $X$ such that $Y\nsubseteq Sing(X)$.
Consider be the line bundle $L_{Y}\rightarrow X$ canonically associated with 
$Y$ and let $N_{Y}\rightarrow Y$ be the restriction at $Y$ of $L_{Y}$.

Let $\mathcal{F}$\ be a holomorphic foliation of rank $1$ of $X$ and suppose
that $Y$ is $\mathcal{F}$-invariant. Suppose that $S=(Sing(\mathcal{F})\cap
Y)\cup Sing(Y)$ is an isolated singular point $p\in Sing(Y)\cap Sing(%
\mathcal{F})\cap Sing(X)$ and that the stalk $\mathcal{F}_{p}$ is generated
on $\mathcal{O}_{X,p}$ by a single element of $\mathcal{T}X_{p}$. Write $%
Y^{\prime \prime }=Y\setminus S$ and recall that $N_{Y^{\prime
}}|_{Y^{\prime \prime }}\rightarrow Y^{\prime \prime }$ is an $%
(F|_{Y^{\prime \prime }})$-bundle with respect to the action $\tau $
described in (\ref{Form. Azione Olom. Tau}).

Let $W_{1}$ be a neighborhood of $p$ open in $X$ such that $W_{1}\cap
Y=\{x\in X:h(x)=0\}$, where $h$ is a local holomorphic definition function
for $Y$ defined on $W_{1}$. Denote by $\mathfrak{h}$\ the non vanishing
section of $L_{Y}|_{W_{1}}\rightarrow X|_{W_{1}}$ associated with $h$. Let $%
V_{1}$ be the neighborhood of $p$ open in $Y$ defined by $V_{1}=W_{1}\cap Y$%
. Shrinking $W_{1}$, if necessary, we can assume without loss of generality
that $V_{1}$ is topologically contractible and that its closure $\overline{%
V_{1}}$ is compact. Set $Z=Y\setminus V_{1}$ and suppose that both $%
\overline{V_{1}}$ and $Z$ are polyhedra with respect to a triangulation $%
\mathbb{T}$\ of $X$ compatible with $Sing(X)\cup Y$ and such that $p$ is in
the interior of some $2$-simplex of $\mathbb{T}$. Shrinking $W_{1}$, if
necessary, we can also assume that on $W_{1}$ the foliation $\mathcal{F}$\
is generated by one holomorphic vector field $\digamma \in \mathcal{T}X$.
Indeed, the sheaf $\mathcal{F}$\ is coherent and $\mathcal{F}_{p}$\ is
generated by only one element of $\mathcal{T}X_{p}$.

Write $V_{0}=Y\setminus \{p\}$, set $\mathcal{V}=\{V_{0},V_{1}\}$ and note
that $Z\subset V_{0}$. Let $\nabla _{0}$ be a linear connection of type $%
(1,0)$ for $N_{Y^{\prime }}|_{V_{0}}\rightarrow V_{0}$ extending the partial
connection $(F_{Y^{\prime \prime }}\oplus \mathbf{\bar{T}}Y|_{Y^{\prime
\prime }},$ $\tau \oplus \bar{\partial})$. Let $\nabla _{1}$ be an $%
\mathfrak{h}|_{V_{1}}$-trivial extendable linear connection for $%
N_{Y}|_{V_{1}}\rightarrow V_{1}$. Then $\check{c}_{ext}^{1}(\nabla _{\ast
})=(c_{ext}^{1}(\nabla _{0}),$ $c_{ext}^{1}(\nabla _{1}),$ $%
c_{ext}^{1}(\nabla _{0},\nabla _{1}))=(0,$ $0,$ $c_{ext}^{1}(\nabla
_{0},\nabla _{1}))$. Indeed, $c_{ext}^{1}(\nabla _{0})=c_{dif}^{1}(\nabla
_{0})=0$, because of Theorem 6.1 of \cite{Abate Bracci Tovena 2}, and $%
c_{ext}^{1}(\nabla _{1})=0$, because $\nabla _{1}$\ is $\mathfrak{h}%
|_{V_{1}} $-trivial. Let $\left\{ \rho _{0},\rho _{1}:Y\rightarrow 
\mathbb{R}
\right\} $ be an extendable partition of unity subordinated to $\mathcal{V}$
and set $\nabla =\rho _{0}\nabla _{0}+\rho _{1}\nabla _{1}$. Then, by (\ref%
{Res. Camacho-Sad}), (\ref{Form. Integr. Omolog. per Var. Comp.}), (\ref%
{Form. Integr. per Var. Comp.}) and (\ref{Cech Integrazione Relativo}), it
results%
\begin{equation}
i_{\ast }(Res_{c_{ext}^{1}}(N_{Y},\mathcal{F},p))=\tint%
\nolimits_{[Y]}c_{ext}^{1}(N_{Y})=\tint\nolimits_{Y}c_{ext}^{1}(\nabla
)=-\tint\nolimits_{Lk(p)}c_{ext}^{1}(\nabla _{0},\nabla _{1}),
\label{Calcolo Res. Camacho-Sad}
\end{equation}%
with $Lk(p)\subset V_{0}\cap V_{1}$ the link of $p$ in $Y$ with respect to a
triangulation $\mathbb{T}_{\bullet }$\ of $X$ compatible with $Sing(X)\cup
Y\cup \{p\}$. As a note, $\mathbb{T}_{\bullet }\mathbb{\neq T}$, because $p$
is not a vertex of $\mathbb{T}$.

So, we only have to explicitly compute the extendable Bott difference form $%
c_{ext}^{1}(\nabla _{0},\nabla _{1})$. To do this, observe that $%
c_{ext}^{1}(\nabla _{0},\nabla _{1})$ is defined on the differentiable
complex manifold $V_{(0,1)}=V_{0}\cap V_{1}$. Consider the differentiable
vector bundle $E=N_{Y}|_{V_{(0,1)}}\times 
\mathbb{R}
\rightarrow V_{(0,1)}\times 
\mathbb{R}
$ and let $\bar{\nabla}$ be the linear connection for $E$ defined by $\bar{%
\nabla}=(1-\varsigma )\nabla _{0}+\varsigma \nabla _{1}$, with $\varsigma
\in 
\mathbb{R}
$. Let $\Xi _{\ast }$ denote the integration along the fibres of the
projection $\Xi :V_{(0,1)}\times \lbrack 0,1]\rightarrow V_{(0,1)}$. Then,
by its very definition, we have $c_{ext}^{1}(\nabla _{0},\nabla _{1})=\Xi
_{\ast }(c_{dif}^{1}(\bar{\nabla}))$.

Let $h$ and $\mathfrak{h}$\ be as above. By the parametrization theorem (cp. 
\cite{Gunning}, Vol. II, Ch. D), we may find a holomorphic function $%
y:W_{1}\rightarrow 
\mathbb{C}
$ defined on $W_{1}$ such that $(dh\wedge dy)|_{Y^{\prime }}$ does not
vanish on a neighborhood $V$ of $Y^{\prime }\setminus \{p\}$ that, without
loss of generality, we can assume to contain $V_{(0,1)}$. Then $(h,y)$ are
local coordinates on $X^{Reg}$ near $p$ and $y$ is a local coordinate on $%
Y^{\prime }$ near each point of $Y^{\prime }\setminus \{p\}$. In particular, 
$y$ is a local coordinate on $V_{(0,1)}\subset V_{1}\setminus \{p\}$. Since $%
Y$ is $\mathcal{F}$-invariant, using the coordinates $(h,y)$ and with slight
abuses of notation, we can write the holomorphic vector field $\digamma \in 
\mathcal{T}X$ generating $\mathcal{F}$\ on $W_{1}$ as $\digamma =a(h,y)h%
\frac{\partial }{\partial h}+b(h,y)\frac{\partial }{\partial y}$, with $a$
and $b$ holomorphic functions defined on $W_{1}$ such that $b(0,y)$ is not
identically equal to zero.

Let $\bar{\theta}$ be the connection form of $\bar{\nabla}$ and $\theta _{0}$%
\ the connection form of $\nabla _{0}$. Since $\nabla _{1}$ is $\mathfrak{h}%
|_{V_{1}}$-trivial, the connection form $\theta _{1}$ with respect to $%
\mathfrak{h}|_{V_{1}}$\ is zero. Then $\bar{\theta}=(1-\varsigma )\theta
_{0} $ and%
\begin{equation}
c_{ext}^{1}(\nabla _{0},\nabla _{1})=\tfrac{\sqrt{-1}}{2\pi }\Xi _{\ast }(d%
\bar{\theta})=\tfrac{1}{2\pi \sqrt{-1}}\theta _{0}
\label{Forma Differ. Camacho-Sad}
\end{equation}%
Now, to compute $\theta _{0}|_{Y^{\prime }\setminus \{p\}}$, look at the
very definition of $\theta _{0}$\ and $\tau $ (cp. (\ref{Form. Azione Olom.
Tau})). In fact, since $\frac{\partial }{\partial h}\in \Gamma (V_{(0,1)},%
\mathbf{T}X|_{V_{(0,1)}})$ is an extension of $\mathfrak{h}|_{V_{(0,1)}}\in
\Gamma (V_{(0,1)},N_{Y}|_{V_{(0,1)}})$ on $V_{(0,1)}$, it results $\theta
_{0}(\frac{\partial }{\partial y})\mathfrak{h}|_{V_{(0,1)}}=(\nabla _{0})_{%
\frac{\partial }{\partial y}}(\mathfrak{h}|_{V_{(0,1)}})=\tau (\frac{%
\partial }{\partial y}|_{Y^{\prime }\setminus \{p\}},$ $\mathfrak{h}%
|_{V_{(0,1)}})=-\frac{a(0,y)}{b(0,y)}\mathfrak{h}|_{V_{(0,1)}}$. Then $%
\theta _{0}=-\frac{a(0,y)}{b(0,y)}dy$. So, by (\ref{Calcolo Res. Camacho-Sad}%
) and (\ref{Forma Differ. Camacho-Sad}), we get the following formula for
the residue%
\begin{equation}
i_{\ast }(Res_{c_{ext}^{1}}(N_{Y},\mathcal{F},p))=\tfrac{1}{2\pi \sqrt{-1}}%
\tint\nolimits_{Lk(p)}\tfrac{a(0,y)}{b(0,y)}dy
\label{Residuo-Foliazione (Espressione Finale)}
\end{equation}

\end{document}